# Existence for low-regularity McKean–Vlasov dynamics via emergence of regularity

By **Robert Alexander Crowell** *of ETH Zürich*

**Abstract:** We establish the existence of solutions to common noise McKean–Vlasov martingale problems for coefficients with low regularity. Our approach is able to handle the key challenge posed by drift coefficients that are discontinuous with respect to the narrow convergence of measures. This case arises for e.g. singular interactions. Our proof strategy proceeds via a two-step approximation using smoothed McKean–Vlasov $n$-particle systems: We first pass to the large system limit by taking $n \to \infty$, and subsequently remove the smoothing. A novel aspect of our work is the use of a crucial emergence of regularity property. It ensures that after the first limit, we obtain a process of measures that are absolutely continuous with respect to the Lebesgue measure and provides quantitative integrability bounds on their densities. We use this regularity to establish a tightness result in a stronger topology than is typically considered. In this way we obtain a sufficiently strong mode of convergence that lets us subsequently remove the smoothing and solve the McKean–Vlasov martingale problem via the particle system approximations.

**Robert Alexander Crowell**

rc@math.ethz.ch

Department of Mathematics
ETH Zürich
Rämistrasse 101
CH-8092 Zürich

**Acknowledgements:** This work grew out of my PhD thesis written at ETH Zürich. I am deeply indebted to my advisor Martin Schweizer who gave me the freedom to pursue this line of research, and at the same time guided me with his immense experience throughout. I benefited greatly from his detailed feedback, and the arguments and the exposition in this manuscript improved markedly from his careful comments.

**CONTENTS**

## INTRODUCTION

This paper is the third in a series studying McKean–Vlasov systems with common noise when the coefficients have only low regularity. Here, we aim to prove the existence of solutions of a *common noise McKean–Vlasov martingale problem* $\mathrm{MP}(\nu, \mathscr{A}, \mathscr{Q})$, where $\mathscr{A}$ and $\mathscr{Q}$ are certain operators on $\mathbf{M}_1^+(\mathbb{R}^d) \subseteq \mathcal{S}'(\mathbb{R}^d)$ specified below, and where $\mathcal{S}'(\mathbb{R}^d)$ is the space of tempered distributions and $\mathbf{M}_1^+(\mathbb{R}^d)$ is the subset of Borel probability measures. More explicitly, we seek a measure $\mathbf{P}$ on $\Omega = \mathbf{C}([0,T]; \mathcal{S}')$ with canonical process $\Lambda = (\Lambda_t)_{t \in [0,T]}$ such that $\Lambda$ can be written under $\mathbf{P}$ as

$$\Lambda_t = \Lambda_0 + \int_0^t \mathscr{A}_s(\Lambda_s)\,\mathrm{d}s + M_t, \qquad \text{for } t \in [0,T]\,, \tag{1}$$

where $\Lambda_0 \overset{d}{=} \nu$ and $M = (M_t)_{t \in [0,T]}$ is under $\mathbf{P}$ an $\mathcal{S}'$-valued martingale with continuous trajectories, whose quadratic variation is given by

$$\langle\!\langle M \rangle\!\rangle_t = \int_0^t \mathscr{Q}_s(\Lambda_s)\,\mathrm{d}s\,, \qquad \text{for } t \in [0,T]\,. \tag{2}$$

In (1), $\mathscr{A}$ is map $[0,T] \times \mathbf{M}_1^+ \to \mathcal{S}'$ given by a second-order differential operator specified in terms of coefficients $b$, $\sigma$, $\bar{\sigma}$. We are interested in the case where these coefficients possess only low regularity properties. Motivating examples include drifts $b : [0,T] \times \mathbb{R}^d \times \mathbf{M}_1^+ \to \mathbb{R}^d$ that are not narrowly continuous in the measure argument, e.g., convolutional drifts $b_t(x,\mu) = (k_t * \mu)(x)$ with an irregular kernel $k$, or singular drifts $b_t(x,\mu) = B_t(x, \delta_x[\mu])$ for an appropriate function $B$. The latter example is defined only if $f_\mu := \mathrm{d}\mu/\mathrm{d}x$ exists so that $\delta_x[\mu] = f_\mu(x)$ has meaning almost everywhere. The low regularity of the coefficients translates to poor narrow continuity properties of $\mathbf{M}_1^+ \ni \lambda \mapsto \mathscr{A}_t(\lambda) \in \mathcal{S}'$, which makes solving $\mathrm{MP}(\nu, \mathscr{A}, \mathscr{Q})$ difficult.

Our main result in Theorem 1.9 shows under rather weak regularity assumptions on $b$, $\sigma$, $\bar{\sigma}$ that *we can solve the martingale problem* $\mathrm{MP}(\nu, \mathscr{A}, \mathscr{Q})$. This result requires ellipticity of $\sigma$ and a type of local continuity in the measure argument $\mu$ of the drift $b$; see Assumptions 1.3 and 1.7 below. Importantly, the continuity requirement is not with respect to narrow convergence but rather with respect to $\mathbb{L}^1(\mathrm{d}x)$-convergence of the densities for absolutely continuous measures. These assumptions are indispensable for solutions of $\mathrm{MP}(\nu, \mathscr{A}, \mathscr{Q})$ to exist and at the same time generous enough to cover interesting cases that appear prominently in the literature such as those mentioned above; see also Section 6 for examples and Section 7 for a discussion.

The possibly non-unique solutions of $\mathrm{MP}(\nu, \mathscr{A}, \mathscr{Q})$ that we produce can be obtained through an approximation by a sequence of *smoothed McKean–Vlasov $n$-particle systems* given for $k \in \mathbb{N}$, $i \in \{1, \ldots, n\}$ and $t \in [0,T]$ by

$$\begin{cases} X_t^{i,n} = X_0^{i,n} + \int_0^t b_s^k(X_s^{i,n}, \mu_s^{n,k})\,\mathrm{d}s + \int_0^t \sigma_s(X_s^{i,n}, \mu_s^n)\,\mathrm{d}B_s^{i,n,k} + \int_0^t \bar{\sigma}_s(X_s^{i,n}, \mu_s^n)\,\mathrm{d}Z_s^n\,, \\ \mu_t^n = \dfrac{1}{n}\sum_{i=1}^n \delta_{X_t^{i,n}}\,. \end{cases} \tag{3}$$

Above $b^k$ and $\mu^{n,k}$ are obtained from $b$ and $\mu^n$ by mollification via an approximate identity $(h_k)_{k \in \mathbb{N}}$. The system (3) possesses a unique weak solution $\mathbb{P}_{x_0}^{n,k}$ under very general assumptions on the coefficients; it is well-defined even for a singular drift $b$ as above. For each $n \in \mathbb{N}$ and $k \in \mathbb{N}$, the law $\mathbf{P}_{x_0}^{n,k} := \mathbb{P}_{x_0}^{n,k} \circ (\mu^n)^{-1}$ solves a smoothed $n$-particle martingale problem $\mathrm{MP}(\nu_n, \mathscr{A}^k, \mathscr{Q} + \frac{1}{n}\mathscr{C})$. Formally, as $n \to \infty$ and $k \to \infty$, and if $\nu_n \to \nu$, we obtain our original martingale problem $\mathrm{MP}(\nu, \mathscr{A}, \mathscr{Q})$ as the limit.

To establish Theorem 1.9, we first consider for a fixed smoothing $k \in \mathbb{N}$ the sequence $(\mathbf{P}_{x_0}^{n,k})_{n \in \mathbb{N}}$. We prove that it is tight and use the results in the second paper [10] of this series to show that its cluster points $\mathbf{P}_{x_0}^{\infty,k}$ solve a smoothed martingale problem $\mathrm{MP}(\nu, \mathscr{A}^k, \mathscr{Q})$. Second, we use another tightness argument to obtain a cluster point $\mathbf{P}_{x_0}^{\infty,\infty}$ of the sequence $(\mathbf{P}_{x_0}^{\infty,k})_{k \in \mathbb{N}}$ of cluster points. This is our solution candidate for $\mathrm{MP}(\nu, \mathscr{A}, \mathscr{Q})$. While natural, implementing this strategy is rather intricate. Key for its success is a careful choice of the topologies with respect to which tightness is considered. Since $\Lambda_t$ is $\mathbf{P}_{x_0}^{n,k}$-a.s. an atomic measure, we require a topology that is not too restrictive; so for $(\mathbf{P}_{x_0}^{n,k})_{n \in \mathbb{N}}$ we prove tightness in $\Omega$. The situation is quite different for the sequence $(\mathbf{P}_{x_0}^{\infty,k})_{k \in \mathbb{N}}$ of cluster points. Tightness in $\Omega$ is not enough to handle coefficients that are not narrowly continuous. We therefore show that $(\mathbf{P}_{x_0}^{\infty,k})_{k \in \mathbb{N}}$, in addition to being tight in $\Omega$, is also tight in $\Omega_{p,1} = \mathbb{L}^p([0,T]; \mathbb{L}^1(\mathrm{d}x))$ for any $1 \le p < \infty$. For this to be possible, we must have that $\Lambda_t \ll \mathrm{d}x$ for almost every $t \in [0,T]$ under $\mathbf{P}_{x_0}^{\infty,k}$. This property, which we termed the *emergence of regularity* in the first paper [8] in this series is vital for our study. It comes with a quantitative bound $\big\| \|\Lambda_t\|_{\mathbb{L}^r(\mathrm{d}x)} \big\|_{\mathbb{L}^q(\mathbf{P}_{x_0}^{\infty,k})} \le c(1 \wedge t)^{-\gamma}$ for certain $r > 1$ and $0 < \gamma < 1$, which we use in various ways. In a first step, we trade off regularity of $\Lambda_t$ in space against regularity in time to obtain for a Skorohod representation $(\hat{\Lambda}^m)_{m \in \mathbb{N} \cup \{\infty\}}$ of a narrowly convergent subsequence $(\mathbf{P}_{x_0}^{\infty,k_m})_{m \in \mathbb{N}}$ with limit $\mathbf{P}_{x_0}^{\infty,\infty}$ an improved mode of convergence. Specifically, we show that for $1 \le \xi < r$ we have almost surely that $(\hat{\Lambda}_t^m)_{t \in (0,T]} \to (\hat{\Lambda}_t^\infty)_{t \in (0,T]}$ in $\mathbf{C}((0,T]; \mathbb{L}^\xi(\mathrm{d}x))$ with uniform convergence on compacts as $m \to \infty$. This mode of convergence is sufficiently strong for even the singular coefficients mentioned as examples above to exhibit continuity properties in mean. These we exploit by using the emergence of regularity a second time, but in a different way. The specific form of the functional $\mathscr{A}_t(\Lambda_t)$ lets us compensate the low regularity of the drift $b_t$ with the good regularity of $\Lambda_t$ to show that $\mathbf{P}_{x_0}^{\infty,\infty}$ solves the martingale problem $\mathrm{MP}(\nu, \mathscr{A}, \mathscr{Q})$.

These arguments and their delicate synthesis are novelties of our work. Together, they provide a systematic approach to overcome the difficulties created by $\mathscr{A}$ being discontinuous with respect to narrow convergence. This is enabled by the emergence of regularity property, which is paired with the semimartingale property of $\Lambda$ under $\mathbf{P}_{x_0}^{\infty,k}$ to deduce and use tightness in a comparatively strong topology. We demonstrate the virtues of our approach by implementing it under Assumptions 1.3 and 1.7. However, we emphasize that the conceptual aspects carry over to more general assumptions and also to other settings. We discuss this in Section 7.

We point out that the solutions of $\mathrm{MP}(\nu, \mathscr{A}, \mathscr{Q})$ yield solutions to the associated stochastic Fokker–Planck equation, a certain SPDE with coefficients given in terms of $b$, $\sigma$, $\bar{\sigma}$. Via superposition theory, solutions of this SPDE are related to solutions of the associated common noise McKean–Vlasov SDE, linking the existence theory of both. Therefore, our result in Theorem 1.9 provides a way to obtain weak solutions of common noise McKean–Vlasov SDEs with irregular or even singular drifts of a certain type. This is true even if the diffusion coefficients $\sigma$, $\bar{\sigma}$ depend on the measure argument.

**Organization**

This paper is structured as follows. The probabilistic setup, the main result in Theorem 1.9 and the outline of its proof are contained in Section 1. Section 2 develops the basic technical terminology and studies the smoothed $n$-particle martingale problem $\mathrm{MP}(\nu_n, \mathscr{A}^{\delta_k}, \mathscr{Q} + \frac{1}{n}\mathscr{C})$. Section 3 develops the necessary regularity results, and Section 4 the key tightness result. Section 5 combines all these ingredients into the proof of the main result in Theorem 1.9. Finally, Sections 6 and 7 conclude this work with examples, and a discussion and an outlook, respectively. The Appendices A and B contain useful results from functional analysis and probability, as well as deferred proofs.

**Notation**

We collect the essential notation. Appendix A contains more background and references.

- Throughout $c(a_1, \ldots, a_k)$ denotes a generic finite constant which depends on the quantities $a_1, \ldots, a_k$, and which may change from line to line.

- The space of *Schwartz functions* is denoted by $\mathcal{S} = \mathcal{S}(\mathbb{R}^d)$ and endowed with its usual Fréchet topology. We write $\mathcal{S}' = \mathcal{S}'(\mathbb{R}^d)$ for the topological vector space of *(tempered) distributions*, and $\lambda[\phi] := \langle \lambda; \phi \rangle_{\mathcal{S}' \times \mathcal{S}}$ for the duality pairing of $\lambda \in \mathcal{S}'$ and $\phi \in \mathcal{S}$. The Fourier transform $\mathscr{F}$ is a linear automorphism in $\mathcal{S}'$ with inverse $\mathscr{F}^{-1}$.

- For a measure-space $(\mathrm{X}, \mathcal{X}, \chi)$ and a Banach space $(B, \|\cdot\|_B)$, we denote the associated *Bochner space* by $\mathbb{L}^p(\chi; B) = \mathbb{L}^p((\mathrm{X}, \mathcal{X}, \chi); B)$. If $B = \mathbb{R}$, we write $\mathbb{L}^p(\chi) = \mathbb{L}^p(\mathrm{X}, \mathcal{X}, \chi)$ for the usual *Lebesgue space*.

- The *Bessel potential space* with integrability $r \in (1, \infty)$ and regularity $s \in \mathbb{R}$ is defined by $\mathsf{H}^s_r = \mathsf{H}^s_r(\mathbb{R}^d) := \{f \in \mathcal{S}'(\mathbb{R}^d) : \|f\|_{\mathsf{H}^s_r(\mathbb{R}^d)} < \infty\}$, where $\|f\|_{\mathsf{H}^s_r(\mathbb{R}^d)} := \|J^s f\|_{\mathbb{L}^r(\mathrm{d}x)}$ with $J^s f = \mathscr{F}^{-1}(h^s(\mathscr{F}f))$ for $f \in \mathcal{S}'$ and $h(\xi) := \langle \xi \rangle := (1 + |\xi|^2)^{1/2}$.

- The *Hermite–Fourier space* with regularity $p \in \mathbb{R}$ is the Hilbert space obtained by endowing $\mathscr{H}_p = \mathscr{H}_p(\mathbb{R}^d) := \{f \in \mathcal{S}'(\mathbb{R}^d) : \|f\|_{\mathscr{H}_p(\mathbb{R}^d)} < \infty\}$ with the (semi-)norm $\|f\|_{\mathscr{H}_p(\mathbb{R}^d)} := (\sum_{k \in \mathbb{N}^d} (\langle k \rangle^p f_k^{\#})^2)^{1/2}$, where $f_k^{\#}$ is the $k$'th Hermite coefficient of $f \in \mathcal{S}'$.

- If $(\mathrm{X}, \mathcal{X})$ is a measurable space then the set of *probability measures* is denoted by $\mathbf{M}_1^+(\mathrm{X}) = \mathbf{M}_1^+(\mathrm{X}, \mathcal{X})$. For $\mu, \nu \in \mathbf{M}_1^+(\mathrm{X})$ with $\int_{\mathrm{X}} |x| \mathrm{d}\mu, \int_{\mathrm{X}} |x| \mathrm{d}\nu < \infty$, we define the *Kantorovich–Rubinstein* distance by $\mathsf{d}_{\mathbb{W}_1}(\mu, \nu) := \sup\{\int f \, \mathrm{d}(\mu - \nu) : \|f\|_{\mathsf{Lip}} \leq 1\}$. We write $\mathcal{P}_{\mathrm{wk}^*} = \mathcal{P}_{\mathrm{wk}^*}(\mathbb{R}^d)$ for the set $\mathbf{M}_1^+ = \mathbf{M}_1^+(\mathbb{R}^d)$ with the *narrow topology*.

- Note that $\mathbf{M}_1^+ \subseteq \mathcal{S}'$ and that $\mathcal{P}_{\mathrm{wk}^*} \hookrightarrow \mathcal{S}'$, which is to say that the embedding is continuous. In particular, if $\mu \in \mathbf{M}_1^+$ and $\phi \in \mathcal{S}$, then $\mu[\phi] = \langle \mu; \phi \rangle_{\mathcal{S}' \times \mathcal{S}} = \int_{\mathbb{R}^d} \phi(x) \, \mu(\mathrm{d}x)$.

# 1 PROBABILISTIC SETUP

## 1.1 Probabilistic setup for the particle systems

Fix $T < \infty$ as the terminal time and $d \in \mathbb{N}$ for the dimension. The basis for our study is a family of $\delta$-*smoothed $n$-particle systems* of McKean–Vlasov type. These are defined for each $\delta \geq 0$, $n \in \mathbb{N}$, $i = 1, \ldots, n$ and $t \in [0, T]$ as the interacting $n$-particle system

$$X_t^{i,n} = X_0^{i,n} + \int_0^t b_s^\delta(X_s^{i,n}, \mu_s^{n,\delta}) \, \mathrm{d}s + \int_0^t \sigma_s(X_s^{i,n}, \mu_s^n) \, \mathrm{d}B_s^{i,n,\delta} + \int_0^t \bar{\sigma}_s(X_s^{i,n}, \mu_s^n) \, \mathrm{d}Z_s^n \,, \quad (1.1)$$

$$\mu_t^n = \frac{1}{n} \sum_{i=1}^n \delta_{X_t^{i,n}} \,. \tag{1.2}$$

For each $\delta \geq 0$, the particle system (1.1) consists of $n$ coupled stochastic differential equations (SDEs), each with values in $\mathbb{R}^d$. We refer to $X^{i,n}$, the $i$'th component of the system, as a *particle*, to $B^{1,n,\delta}, \ldots, B^{n,n,\delta}$ as the *idiosyncratic noise*, to $Z^n$ as the *common noise*, to $\mu_t^n$ as the *empirical measure* and to $\mu_t^{n,\delta}$ as the *smoothed empirical measure*.

For $\delta = 0$, we set $b^0 := b$, $\mu^{n,0} := \mu^n$ and $B^{i,n} := B^{i,n,0}$. For $\delta > 0$, the drift coefficient $b^\delta$ and the smoothed empirical measure $\mu^{n,\delta}$ is obtained from $b$, respectively $\mu^n$, by an appropriate mollification procedure that we make precise momentarily. To state the assumptions on the coefficients $b$, $\sigma$, $\bar{\sigma}$, we need on $[0, T] \times \mathbb{R}^d \times \mathbf{M}_1^+(\mathbb{R}^d)$ the product $\sigma$-algebra $\mathcal{E}$ generated by the Lebesgue-measurable sets of $[0, T] \times \mathbb{R}^d$ and the Borel-measurable sets of $\mathcal{P}_{\mathrm{wk}^*}(\mathbb{R}^d)$.

**Assumption 1.1** | The random variables $(X_0^{i,n})_{i=1,\ldots,n}$ are i.i.d. Moreover the law of $X_0^{1,n}$ is absolutely continuous with respect to Lebesgue measure, with density $\nu$.

**Assumption 1.2** | The drift coefficient $b : [0, T] \times \mathbb{R}^d \times \mathbf{M}_1^+(\mathbb{R}^d) \to \mathbb{R}^d$ is bounded and measurable relative to $\mathcal{E}$.

**Assumption 1.3** | The diffusion coefficients $\sigma : [0, T] \times \mathbb{R}^d \times \mathbf{M}_1^+(\mathbb{R}^d) \to \mathbb{R}^{d \times d}$ as well as $\bar{\sigma} : [0, T] \times \mathbb{R}^d \times \mathbf{M}_1^+(\mathbb{R}^d) \to \mathbb{R}^{d \times m}$ are bounded and measurable relative to $\mathcal{E}$. In addition,

(E) there is $\kappa > 0$ such that the uniform ellipticity condition

$$z^\top \Big( (\sigma_t \sigma_t^\top)(x, \mu) \Big) z \geq \kappa \, |z|^2$$

holds for all $z \in \mathbb{R}^d$, $t \in [0, T]$ and $(x, \mu) \in \mathbb{R}^d \times \mathbf{M}_1^+(\mathbb{R}^d)$;

(H) there exist $\beta > 0$ and a constant $C < \infty$ such that the Hölder-regularity conditions

$$|\sigma_t(x, \mu) - \sigma_t(x', \mu')| \leq C \Big( |x - x'|^\beta + \mathsf{d}_{\mathbb{W}_1}(\mu, \mu')^\beta \Big) \,,$$

$$|\bar{\sigma}_t(x, \mu) - \bar{\sigma}_{t'}(x', \mu')| \leq C \Big( |t - t'|^\beta + |x - x'|^\beta + \mathsf{d}_{\mathbb{W}_1}(\mu, \mu')^\beta \Big) \,,$$

hold for all $t, t' \in [0, T]$ and $(x, \mu), (x', \mu') \in \mathbb{R}^d \times \mathbf{M}_1^+(\mathbb{R}^d)$.

Let $h \in \mathbf{C}_c^\infty(\mathbb{R}^d)$ be a symmetric function such that $h \geq 0$ and $\int_{\mathbb{R}^d} h(x) \, \mathrm{d}x = 1$, and let $h_\delta := \delta^d h(\delta x)$ for $\delta > 0$. Then $(h_\delta)_{\delta > 0}$ is a compactly supported approximate identity; see Rudin [35, Def. 6.31]. For $\delta > 0$ and $\lambda \in \mathbf{M}_1^+(\mathbb{R}^d)$, we first define the smooth probability density

$$\lambda^\delta(x) := (\lambda * h_\delta)(x) := \int_{\mathbb{R}^d} h_\delta(y - x) \, \lambda(\mathrm{d}y) \qquad \text{for } x \in \mathbb{R}^d, \tag{1.3}$$

which we identify with the associated measure $\lambda^\delta \in \mathbf{M}_1^+(\mathbb{R}^d)$. It is a standard result that $\lambda^\delta \to \lambda$ narrowly as $\delta \to 0$. In the case of the empirical measure (1.2), this gives the smoothed empirical measure $\mu_t^{n,\delta} \in \mathbf{M}_1^+(\mathbb{R}^d)$ associated with

$$\mu_t^{n,\delta}(x) := (\mu_t^n * h_\delta)(x) = \frac{1}{n} \sum_{i=1}^n h_\delta(X_t^{i,n} - x) \qquad \text{for } x \in \mathbb{R}^d \,. \tag{1.4}$$

Similarly, we define

$$b^\delta : [0,T] \times \mathbb{R}^d \times \mathbf{M}_1^+ \to \mathbb{R}^d$$

for each $\delta > 0$ by setting

$$b_t^\delta(x,\lambda) := \big(b_t(\lambda) * h_\delta\big)(x) := \int_{\mathbb{R}^d} b_t(y,\lambda) h_\delta(x-y) \, \mathrm{d}y \,, \tag{1.5}$$

where we denote by $b_t(\lambda)$ for $(t,\lambda) \in [0,T] \times \mathbf{M}_1^+$ the function $b_t(\,\cdot\,,\lambda) : \mathbb{R}^d \to \mathbb{R}^d$. The following standard result can be found in e.g. Makarov and Podkorytov [28, Cor. 9.3.2].

**Lemma 1.4** | *Under Assumption 1.2 and for each $\delta > 0$, the $\delta$-smoothed drift coefficient $b^\delta : [0,T] \times \mathbb{R}^d \times \mathbf{M}_1^+(\mathbb{R}^d) \to \mathbb{R}^d$ is well-defined and we have that $\sup_{\delta > 0} \|b^\delta\|_\infty \leq \|b\|_\infty$.*

The $n$-particle system (1.1) is solved in the usual weak sense of SDEs; cf. Stroock and Varadhan [40, Ch. 6]. To construct such a solution, we write

$$\Omega_{\mathrm{par}}^n := \mathbf{C}([0,T]; \mathbb{R}^{nd} \times \mathbb{R}^m)$$

with the coordinate process representing the particles and the common noise, i.e.,

$$(t,\omega) \mapsto \omega_t =: (X_t^n, Z_t^n)(\omega) = (X_t^{1,n}, \dots, X_t^{n,n}, Z_t^n)(\omega) \qquad \text{for } (t,\omega) \in [0,T] \times \Omega_{\mathrm{par}}^n \,.$$

We endow $\Omega_{\mathrm{par}}^n$ with the filtration $\mathbb{G}^n := (\mathcal{G}_t^n)_{0 \leq t \leq T}$ given by the right-continuous version of the canonical raw filtration, i.e.,

$$\mathcal{G}_t^n := \bigcap_{0 < \varepsilon < T-t} \sigma\big((X_s^n, Z_s^n) : 0 \leq s < t+\varepsilon\big) \qquad \text{for } t \in [0,T),$$

$$\mathcal{G}^n := \mathcal{G}_T^n := \sigma\big((X_s^n, Z_s^n) : 0 \leq s \leq T\big) \,.$$

**Lemma 1.5** | *Let $\delta > 0$. Under Assumptions 1.1, 1.2 and 1.3, there exist a unique measure $\mathbb{P}_{x_0}^{n,\delta}$ on $(\Omega_{\mathrm{par}}^n, \mathcal{G}^n)$ and independent standard $d$-dimensional $(\mathbb{G}^n, \mathbb{P}_{x_0}^{n,\delta})$-Brownian motions $B^{1,n,\delta}, \dots, B^{n,n,\delta}$ such that $Z^n$ is a standard $m$-dimensional $(\mathbb{G}^n, \mathbb{P}_{x_0}^{n,\delta})$-Brownian motion that is independent of $B^{n,\delta} := (B^{1,n,\delta}, \dots, B^{n,n,\delta})$, and $(X^n, B^{n,\delta}, Z^n)$ solves (1.1), (1.2). Moreover, the law of $X^n$ under $\mathbb{P}_{x_0}^{n,\delta}$ is exchangeable and we have*

$$\mathbf{E}^{\mathbb{P}_{x_0}^{n,\delta}} \left[ \sup_{t \in [0,T]} |X_t^{i,n}|^q \right] =: c_q < \infty \qquad \text{for any } q > 1, \tag{1.6}$$

*where the constant $c_q$ depends only on $q$, $T$ and the bounds for $b$, $\sigma$, $\bar\sigma$, but not on $n$ or $\delta$.*

**Proof**    With Lemma 1.4, the result can be proved just as [8, Lem. 1.3].    □

By exchangeability, we mean that for any permutation $\pi$ of $[n] := \{1, \ldots, n\}$, we have
$$\text{Law}_{\mathbb{P}_{x_0}^{n,\delta}}(X^{1,n}, \ldots, X^{n,n}) = \text{Law}_{\mathbb{P}_{x_0}^{n,\delta}}(X^{\pi(1),n}, \ldots, X^{\pi(n),n}) \, .$$
In the sequel, we denote the weak solution from Lemma 1.5 by
$$\boldsymbol{X}^{n,\delta} = \left( \Omega_{\text{par}}^n, \mathcal{G}^n, \mathbb{G}^n, \mathbb{P}_{x_0}^{n,\delta}, (X_t^n)_{t \in [0,T]}, (B_t^{n,\delta})_{t \in [0,T]}, (Z_t^n)_{t \in [0,T]} \right) \, . \tag{1.7}$$
Note that Assumption 1.3(E) also allows $\bar{\sigma} = 0$, in which case there is no common noise $Z^n$ in (1.1) and the above simplifies in the obvious way.

## 1.2    Probabilistic setup for the measure flow

From the $n$-particle system (1.2) we have the *empirical measure flow* $\mu^n := (\mu_t^n)_{t \in [0,T]}$. The sequence $(\mu^n)_{n \in \mathbb{N}}$ of these processes is a key object in this work. To focus on it, we transfer it to a canonical setup. For this, let
$$\Omega := \mathbf{C}\big( [0,T]; \, \mathcal{S}'(\mathbb{R}^d) \big)$$
denote the *canonical trajectory space*. Its topological vector space structure is induced by convergence uniformly on $[0,T]$ for bounded sets in $\mathcal{S}(\mathbb{R}^d)$; see e.g. Kallianpur and Xiong [24, Ch. 2.4] or Appendix A. We denote the canonical process on $\Omega$ by
$$(t, \omega) \mapsto \omega_t =: \Lambda_t(\omega) = \Lambda_t \qquad \text{for } (t, \omega) \in [0,T] \times \Omega \, ,$$
We endow $\Omega$ with the right-continuous version $\mathbb{F} := (\mathcal{F}_t)_{0 \leq t \leq T}$ of the the canonical filtration generated by $\Lambda$, i.e.,
$$\mathcal{F}_t := \bigcap_{0 < \varepsilon < T-t} \sigma(\Lambda_s : 0 \leq s < t + \varepsilon) \qquad \text{for } t \in [0,T),$$
$$\mathcal{F} := \mathcal{F}_T := \sigma(\Lambda_s : 0 \leq s \leq T) \, .$$

The set $\mathbf{M}_1^+(\mathbb{R}^d)$ is not closed in the topology of $\mathcal{S}'(\mathbb{R}^d)$. We therefore need to establish part 1) of the following technical result.

**Lemma 1.6** │ *Let $\delta \geq 0$. Under Assumptions 1.1, 1.2 and 1.3, let $\boldsymbol{X}^{n,\delta}$ be the weak solution from Lemma 1.5 of the $\delta$-smoothed $n$-particle system (1.1), (1.2). Then there exists a version of the process $\mu^n = (\mu_t^n)_{t \in [0,T]}$ which has continuous trajectories in the strong topology of $\mathcal{S}'(\mathbb{R}^d)$. In particular,*
$$\mathbf{P}_{x_0}^{n,\delta} := \mathbb{P}_{x_0}^{n,\delta} \circ (\mu^n)^{-1} \tag{1.8}$$
*is a well-defined probability measure on $(\Omega, \mathcal{F})$. Moreover:*

*1) The sequence $(\mathbf{P}_{x_0}^{n,\delta})_{n \in \mathbb{N}}$ is tight. Any narrow cluster point $\mathbf{P}_{x_0}^{\infty,\delta}$ of $(\mathbf{P}_{x_0}^{n,\delta})_{n \in \mathbb{N}}$ satisfies*
$$\mathbf{P}_{x_0}^{\infty,\delta}\big[ \, \mathbf{C}\big( [0,T]; \, \mathcal{S}'(\mathbb{R}^d) \cap \mathbf{M}_1^+(\mathbb{R}^d) \big) \, \big] = 1 \, , \tag{1.9}$$
*that is, $\mathbf{P}_{x_0}^{\infty,\delta}$ concentrates its mass on the set of probability-measure-valued processes.*

*2) For any $\varepsilon, K > 0$ and $q > 1$, we have for all $n \in \mathbb{N} \cup \{\infty\}$ the concentration bound*
$$\mathbf{P}_{x_0}^{n,\delta}\Big[ \omega \in \Omega : \Lambda_t\big[ [-K, K]^d \big] < 1 - \varepsilon \text{ for some } t \in [0,T] \Big] \leq \frac{c_q}{\varepsilon K^q} \, , \tag{1.10}$$
*where $c_q$ is the constant from (1.6), which depends on $q$ but not on $K$, $\varepsilon$ nor $\mathbf{P}_{x_0}^{n,\delta}$ or $\delta$.*

***Proof***  In view of Lemma 1.4, this can be proved in a similar way as [8, Lem. 1.4].  □

Lemma 1.6 lets us consider for each $\delta \geq 0$ and $n \in \mathbb{N}$ the $\mathcal{S}'$-valued canonical process $\Lambda$ under $\mathbf{P}_{x_0}^{n,\delta}$. Note the distinction between $\mathbb{P}_{x_0}^{n,\delta}$ and $\mathbf{P}_{x_0}^{n,\delta}$: the former denotes the joint law of $n$ particles and the common Brownian motion, the latter models the flow of the empirical measure of the particles. Accordingly, we term $\mathbf{P}_{x_0}^{n,\delta}$ the *empirical measure law* of the $\delta$-smoothed $n$-particle system. For $\delta = 0$, we simply write $\mathbf{P}_{x_0}^n := \mathbf{P}_{x_0}^{n,0}$ and note that with the notation we chose after (1.1) and (1.2) this is precisely the law of the $n$-particle system that is considered in [9].

## 1.3  Main result

For each $\delta \geq 0$ and $n \in \mathbb{N}$, the dynamics of $\Lambda$ under $\mathbf{P}_{x_0}^{n,\delta}$ can be described explicitly. Indeed, in Section 2.3 below, we recall that $\Lambda$ under $\mathbf{P}_{x_0}^{n,\delta}$ is an $\mathcal{S}'$-valued semimartingale with continuous trajectories, which means that there is a unique decomposition $\Lambda = A^\delta + M^\delta$, where $A^\delta$ is an $\mathcal{S}'$-valued process of finite variation and $M^\delta$ is an $\mathcal{S}'$-valued local martingale, both with continuous trajectories in the strong topology of $\mathcal{S}'$. In this decomposition, the tensor quadratic variation $\langle\!\langle M^\delta \rangle\!\rangle$ of $M^\delta$ depends on $n$. We refer to Section 2.1, below for a brief recap of these terms.

The processes $A^\delta$ and $M^\delta$ in the semimartingale decomposition can be represented efficiently using a *martingale problem*. We develop the necessary terminology and notation for this in Section 2.2, below. This description ultimately allows us to represent in our main result in Theorem 1.9 the limiting dynamics of $\Lambda$ under successive cluster points $\mathbf{P}_{x_0}^{\infty,\infty}$, obtained informally from the double-index sequence $(\mathbf{P}_{x_0}^{n,\delta})_{n \in \mathbb{N}, \delta > 0}$ by first letting $n \to \infty$ and then $\delta \to 0$.

To state our main result and discuss the strategy of proof, we preview the basic recurrent notation that we require. The data of the martingale problem is given in terms of maps defined for each $\delta \geq 0$ and $t \in [0, T]$ by

$$\mathscr{A}_t^\delta : \mathcal{S}' \to \mathcal{S}' \cup \{\infty\}, \tag{1.11}$$

$\lambda \mapsto \mathscr{A}_t^\delta(\lambda)$, and

$$\mathscr{Q}_t : \mathcal{S}' \to \mathcal{L}_{\mathrm{bl}}^+(\mathcal{S}), \tag{1.12}$$

$$\mathscr{C}_t : \mathcal{S}' \to \mathcal{L}_{\mathrm{bl}}^+(\mathcal{S}), \tag{1.13}$$

with $\lambda \mapsto \mathscr{Q}_t(\lambda)$ and $\lambda \mapsto \mathscr{C}_t(\lambda)$, respectively, where $\mathcal{L}_{\mathrm{bl}}^+(\mathcal{S}')$ is a space of bilinear maps $\mathcal{S} \times \mathcal{S} \to \mathbb{R}$. In (2.13)–(2.15), below, we make $\mathscr{A}^\delta$, $\mathscr{Q}$, $\mathscr{C}$ explicit terms of $b^\delta$, $\sigma$, $\bar{\sigma}$. For simplicity of notation, we set $\mathscr{A} := \mathscr{A}^0$, which makes sense because $b := b^0$.

As we recall in Proposition 2.6, the measure $\mathbf{P}_{x_0}^{n,\delta}$ solves for each $n \in \mathbb{N}$ the canonical $\mathcal{S}'$-valued martingale problem $\mathrm{MP}(\nu_n, \mathscr{A}^\delta, \mathscr{Q} + \frac{1}{n}\mathscr{C})$; see Definition 2.3 for what precisely this entails. As a consequence, the $\mathcal{S}'$-valued semimartingale decomposition of $\Lambda$ under $\mathbf{P}_{x_0}^{n,\delta}$ can be written explicitly as

$$\Lambda_t = \Lambda_0 + \int_0^t \mathscr{A}_s^\delta(\Lambda_s)\,\mathrm{d}s + M_t^\delta, \qquad t \in [0, T], \tag{1.14}$$

where $M^\delta$ is in fact an $\mathcal{S}'$-valued square-integrable martingale with continuous trajectories, i.e., $M^\delta \in \mathcal{M}_2^c(\mathbb{F}, \mathbf{P}_{x_0}^{n,\delta}; \mathcal{S}')$, with (tensor) quadratic variation

$$\langle\!\langle M^\delta \rangle\!\rangle_t = \int_0^t \mathcal{Q}_s(\Lambda_s) \, \mathrm{d}s + \frac{1}{n} \int_0^t \mathcal{C}_s(\Lambda_s) \, \mathrm{d}s \,, \qquad t \in [0, T] \,. \tag{1.15}$$

This is true already under Assumptions 1.2 and 1.3. Under the same assumptions, the sequence $(\mathbf{P}_{x_0}^{n,\delta})_{n \in \mathbb{N}}$ is tight and thus possesses at least one cluster point $\mathbf{P}_{x_0}^{\infty,\delta}$. However, cluster points need not solve the formally limiting martingale problem $\mathrm{MP}(\nu, \mathscr{A}^\delta, \mathscr{Q})$ obtained by letting $n \to \infty$; see, e.g., [8, Ex. 7.1]. But, using the result in [10], we can show that whenever $\delta > 0$, any cluster point $\mathbf{P}_{x_0}^{\infty,\delta}$ solves $\mathrm{MP}(\nu, \mathscr{A}^\delta, \mathscr{Q})$ if we strengthen Assumption 1.2 as follows.

**Assumption 1.7** | The drift coefficient $b : [0, T] \times \mathbb{R}^d \times \mathbf{M}_1^+(\mathbb{R}^d) \to \mathbb{R}^d$ is bounded and measurable relative to $\mathcal{E}$. In addition, $b$ satisfies:

($C_w$) If $(\mu_k)_{k \in \mathbb{N} \cup \{\infty\}}$ in $\mathbf{M}_1^+(\mathbb{R}^d)$ consists of measure that are absolutely continuous with respect to Lebesgue measures and $\mathrm{d}\mu_k/\mathrm{d}x \to \mathrm{d}\mu_\infty/\mathrm{d}x$ in $\mathbb{L}^1(\mathrm{d}x)$ as $k \to \infty$, then

$$\lim_{k \to \infty} \|b_t(\mu_k) - b_t(\mu_\infty)\|_{\mathbb{L}^1(K, \mathrm{d}x)} = 0$$

for any compact subset $K \subseteq \mathbb{R}^d$ and $t \in [0, T]$.

Note that Assumption 1.7 is weaker than the continuity property ($C_s$) that is imposed in [10]; see Lemma 2.7 below for this property. We should point out that the boundedness requirement and the continuity property ($C_w$) in Assumption 1.7 can be weakened; see Section 7 for a discussion. Specifically, in ($C_w$) we only need to consider sequences with $\mathrm{d}\mu_k/\mathrm{d}x \to \mathrm{d}\mu_\infty/\mathrm{d}x$ in $\mathbb{L}^\xi(\mathrm{d}x)$ as $k \to \infty$ for a certain $\xi > 1$ that depends on the dimension $d$ and the Hölder-regularity $\beta$ from Assumption 1.3.

We now introduce a two-step procedure to produce from the particle systems (1.1), (1.2) solution candidates for the martingale problem $\mathrm{MP}(\nu, \mathscr{A}, \mathscr{Q})$, which is obtained as the formal limit as $\delta \to 0$. In the first step, we

$$\text{fix a sequence } (\delta_k)_{k \in \mathbb{N}} \text{ in } (0, \infty) \text{ with } \delta_k \to 0 \text{ as } k \to \infty. \tag{1.16}$$

For all $k, n \in \mathbb{N}$, we then define

$$\mathbf{P}_{x_0}^{n,k} := \mathbf{P}_{x_0}^{n,\delta_k} := \mathbb{P}_{x_0}^{n,\delta_k} \circ (\mu^n)^{-1} \,. \tag{1.17}$$

To simplify notation, set $\delta_\infty = 0$, so that $b^{\delta_\infty} = b$. By Lemma 1.6, the sequence $(\mathbf{P}_{x_0}^{n,k})_{n \in \mathbb{N}}$ is for each fixed $k \in \mathbb{N}$ tight. In the second step, we can thus

$$\text{chose for each } k \in \mathbb{N} \text{ a cluster point } \mathbf{P}_{x_0}^{\infty,k} \text{ of } (\mathbf{P}_{x_0}^{n,k})_{n \in \mathbb{N}}. \tag{1.18}$$

Now let $\Omega_{p,1} := \mathbb{L}^p([0, T]; \mathbb{L}^1(\mathrm{d}x))$; see (4.1).

**Lemma 1.8** | Under Assumptions 1.1, 1.3 and 1.7, let $(\mathbf{P}_{x_0}^{\infty,k})_{k \in \mathbb{N}}$ be chosen as in (1.16), (1.18). For $1 \le p < \infty$, the sequence $(\mathbf{P}_{x_0}^{\infty,k})_{k \in \mathbb{N}}$ is tight in $\Omega$ and in $\Omega_{p,1}$. In particular, there exists a successive cluster point $\mathbf{P}_{x_0}^{\infty,\infty}$, for which we can find a subsequence $(\mathbf{P}_{x_0}^{\infty,k_m})_{m \in \mathbb{N}}$ which converges narrowly to $\mathbf{P}_{x_0}^{\infty,\infty}$ both in $\Omega$ and in $\Omega_{p,1}$.

With the notion of successive cluster points just given, we can state our main result.

**Theorem 1.9** | *Under Assumptions 1.1, 1.3 and 1.7, let $(\mathbf{P}_{x_0}^{\infty,k})_{k\in\mathbb{N}}$ be chosen as in (1.16), (1.18). For $p \geq 1$, let $\mathbf{P}_{x_0}^{\infty,\infty}$ be a successive cluster point as in Lemma 1.8. Then $\mathbf{P}_{x_0}^{\infty,\infty}$ solves the martingale problem $\mathrm{MP}(\nu, \mathscr{A}, \mathscr{Q})$. In particular, the coordinate process $\Lambda$ on $(\Omega, \mathcal{F}, \mathbf{P}_{x_0}^{\infty,\infty})$ possesses the $\mathcal{S}'$-valued semimartingale decomposition*

$$\Lambda_t = \Lambda_0 + \int_0^t \mathscr{A}_s(\Lambda_s)\,\mathrm{d}s + M_t, \qquad t \in [0,T]\,, \tag{1.19}$$

*where $\Lambda_0 \overset{d}{=} \nu$ under $\mathbf{P}_{x_0}^{\infty,\infty}$ and where $M = (M_t)_{t\in[0,T]}$ is in $\mathscr{M}_2^c(\mathbb{F}, \mathbf{P}_{x_0}^{\infty,\infty}; \mathcal{S}')$ with quadratic variation*

$$\langle\langle M \rangle\rangle_t = \int_0^t \mathscr{Q}_s(\Lambda_s)\,\mathrm{d}s, \qquad t \in [0,T]\,. \tag{1.20}$$

*Moreover, for all $1 \leq q < \infty$, there exist exponents $1 < r < \infty$, $0 < w < \infty$ as well as $0 < \gamma < 1$, and a constant $c = c_{\mathrm{Thm.\ 3.1}} < \infty$, all independent of $\mathbf{P}_{x_0}^{\infty,\infty}$, such that the function $t \mapsto (\omega \mapsto \mathrm{d}\Lambda_t(\omega)/\mathrm{d}x)$ defines a strongly measurable map*

$$p^\infty : (0,T] \to \mathbb{L}^q\big((\Omega, \mathcal{F}, \mathbf{P}_{x_0}^{\infty,\infty}); \mathsf{H}_r^w(\mathbb{R}^d)\big)$$

*which satisfies for some constant the bound*

$$\Big\|\|p^\infty(t)\|_{\mathsf{H}_r^w(\mathbb{R}^d)}\Big\|_{\mathbb{L}^q(\mathbf{P}_{x_0}^{\infty,\infty})} \leq c_{\mathrm{Thm.\ 3.1}}(1 \wedge t)^{-\gamma}$$

*for all $t \in (0,T]$. In particular, for each $t \in (0,T]$, we have $\mathbf{P}_{x_0}^{\infty,\infty}$-a.s. that $\Lambda_t \ll \mathrm{d}x$ and $\mathrm{d}\Lambda_t(\omega)/\mathrm{d}x = p^\infty(t)(\omega)$ in $\mathsf{H}_r^w(\mathbb{R}^d)$ $(\mathrm{d}t \otimes \mathrm{d}\mathbf{P}_{x_0}^{\infty,\infty})$-a.e.*

The integrals in (1.19) and (1.20) are understood in the Bochner sense and $\mathsf{H}_r^w(\mathbb{R}^d)$ denotes the Bessel potential space with integrability parameter $r$ and regularity $w$.

Note that the martingale problem $\mathrm{MP}(\nu_n, \mathscr{A}, \mathscr{Q})$ is specified entirely in terms of the unsmoothed coefficients $b$, $\sigma$, $\bar{\sigma}$. Theorem 1.9 shows that any measure $\mathbf{P}_{x_0}^{\infty,\infty}$ obtained via the two-step procedure (1.16), (1.18) solves $\mathrm{MP}(\nu, \mathscr{A}, \mathscr{Q})$, which we view as the formally limiting problem obtained from the sequence $(\mathrm{MP}(\nu, \mathscr{A}^{\delta_k}, \mathscr{Q} + \frac{1}{n}\mathscr{C}))_{n\in\mathbb{N}, k\in\mathbb{N}}$ by first letting $n \to \infty$ and then $\delta_k \to 0$. This should be contrasted to the setting in [10], where for the unsmoothed coefficients, i.e. for $\delta_\infty = 0$ the limit in $n \to \infty$ is considered. As we explain below, the behavior of $\Lambda$ under the sequence $(\mathbf{P}_{x_0}^n)_{n\in\mathbb{N}}$ is more difficult to control and thus requires a more regular drift than we ask for in Assumption 1.7. Assumptions 1.3 and 1.7 allow us to cover interesting and relevant classes of drifts, e.g. certain singular interactions, that do not satisfy the stricter assumptions in [10].

As a particular consequence of Theorem 1.9 we obtain the following result, which we already previewed in [10].

**Corollary 1.10** | *Under Assumptions 1.1, 1.3 and 1.7, there exists a measure $\mathbf{P}$ on $(\Omega, \mathcal{F})$ solving the martingale problem $\mathrm{MP}(\nu, \mathscr{A}, \mathscr{Q})$.*

## 1.4 Strategy of Proof

The main difficulty in the proof of Theorem 1.9 is to establish the $(\mathbf{P}_{x_0}^{\infty,\infty}, \mathbb{F})$-martingale property of $M$ in (1.19). We now discuss why this is challenging and how we overcome this challenge.

**Step 1 — Smoothed particle systems**    The starting point of our analysis is the $\delta$-smoothed $n$-particle systems defined in (1.1), (1.2). Because the measure $\mu_t^{n,\delta}$ possesses a smooth density with respect to the Lebesgue measure, this systems are rather well-behaved and defined even for highly irregular or singular drifts, e.g., when $b_t(x, \mu) = B_t(x, f_\mu(x))$ with $f_\mu := \mathrm{d}\mu/\mathrm{d}x$. The laws $\mathbf{P}_{x_0}^{n,k}$, obtained via (1.16) and (1.17), form the basis of our study.

**Step 2 — The martingale problem formulation and a basic tightness result**    As the second step, we make the sequence of laws $(\mathbf{P}_{x_0}^{k,n})_{k,n\in\mathbb{N}}$ akin to a systematic analysis that ultimately leads to Theorem 1.9. An efficient framework for this is provided by $\mathcal{S}'$-valued processes and martingale problems, for which we recall the basic definitions in Section 2.1. In Section 2.2, we introduce the objects needed to formulate our martingale problems, specifically, the maps $\mathscr{A}^{\delta_k}$, $\mathscr{Q}$, and $\mathscr{C}$, which are derived from the particle system (1.1), (1.2) using Itô's formula.

In Proposition 2.6, we show that for each $k, n \in \mathbb{N}$, the empirical measure flow law $\mathbf{P}_{x_0}^{n,k}$ solves the *martingale problem* $\mathrm{MP}(\nu_n, \mathscr{A}^{\delta_k}, \mathscr{Q} + \frac{1}{n}\mathscr{C})$. This implies that the process $M^{\delta_k}$ from (1.14) is a $(\mathbf{P}_{x_0}^{n,k}, \mathbb{F})$-martingale, so that

$$\mathbf{E}^{\mathbf{P}_{x_0}^{n,k}}[(M_t^{\delta_k} - M_s^{\delta_k})g] = 0 \qquad \text{in } \mathcal{S}' \tag{1.21}$$

for all $0 \le s \le t \le T$ and all continuous, bounded, and $\mathcal{F}_s$-measurable functions $g : \Omega \to \mathbb{R}$. The expectation in (1.21) is a Bochner integral, which allows us to equivalently state that $\mathbf{E}^{\mathbf{P}_{x_0}^{n,k}}[(M_t[\phi] - M_s[\phi])g] = 0$ for all test functions $\phi \in \mathcal{S}$.

In Lemma 2.7, we establish for each fixed $k \in \mathbb{N}$ that the sequence $(\mathbf{P}_{x_0}^{n,k})_{n\in\mathbb{N}}$ is tight in $\Omega = \mathbf{C}([0, T]; \mathcal{S}')$. Starting from (1.21), we then pass to the limit $n \to \infty$ and deduce that any narrow limit point $\mathbf{P}_{x_0}^{\infty,k}$ of $(\mathbf{P}_{x_0}^{n,k})_{n\in\mathbb{N}}$ satisfies

$$\mathbf{E}^{\mathbf{P}_{x_0}^{\infty,k}}[(M_t^{\delta_k} - M_s^{\delta_k})g] = 0 \qquad \text{in } \mathcal{S}'. \tag{1.22}$$

Since in (1.1), (1.2) the coefficients are regularized, this is relatively standard. However, one minor technicality must be addressed. The traditional method for taking the limit $n \to \infty$ in (1.21) requires the the $\delta_m$-smoothed coefficient $b_t^{\delta_k}(x, \lambda^{\delta_k})$ to be jointly continuous in $(x, \lambda)$ when the space of probability measures $\mathbf{M}_1^+(\mathbb{R}^d)$ is equipped with the narrow topology. At this stage, $b_t^{\delta_k}(x, \lambda^{\delta_k})$ is only continuous in $\lambda$ and in $x$ separately. To resolve this minor point, we prove a simple result in Lemma 2.7 that allows us to use the results from the second paper [10] in this series, pass to the limit $n \to \infty$, and deduce (1.22).

**The main difficulties**    Next, we want to use (1.22) to show that any narrow limit point $\mathbf{P}_{x_0}^{\infty,\infty}$ of $(\mathbf{P}_{x_0}^{\infty,k})_{k\in\mathbb{N}}$ satisfies

$$\mathbf{E}^{\mathbf{P}_{x_0}^{\infty,\infty}}[(M_t - M_s)g] = 0 \qquad \text{in } \mathcal{S}', \tag{1.23}$$

for all $0 \le s \le t \le T$ and for $g$ as above. If we can establish (1.23), then a standard monotone class argument shows that $M$ is a $(\mathbf{P}_{x_0}^{\infty,\infty}, \mathbb{F})$-local martingale. This is the crucial property needed to prove that $\mathbf{P}_{x_0}^{\infty,\infty}$ solves the martingale problem $\mathrm{MP}(\nu, \mathscr{A}, \mathscr{Q})$, which is the formal limit of the sequence $(\mathrm{MP}(\nu_n, \mathscr{A}^{\delta_k}, \mathscr{Q} + \frac{1}{n}\mathscr{C}))_{n,k\in\mathbb{N}}$ when $n \to \infty$ and $k \to \infty$.

However, this step is far more intricate. A natural approach is to show, as in Step 2, that the sequence $(\mathbf{P}_{x_0}^{\infty,k})_{k\in\mathbb{N}}$ is tight in $\Omega$, and then try to deduce (1.23) from (1.21). However, such a tightness result in $\Omega$ is not strong enough for this purpose. We must therefore obtain a refined tightness result.

**Step 3 — The emergence of regularity**    In Section 3.1, we make use of the *emergence of regularity* result from the first paper in this series [8]. We show that for any $q \geq 1$, there exist exponents $r > 1$, $s > 0$, and $0 < \gamma < 1$ such that for all $k \in \mathbb{N}$,

$$\mathbf{P}_{x_0}^{\infty,k}\Big[\big\{\omega \in \Omega \,:\, \Lambda_t(\omega) \text{ is in } \mathsf{H}_r^w(\mathbb{R}^d) \text{ for almost all } t \in (0,T]\big\}\Big] = 1\,. \tag{1.24}$$

This property extends to the narrow cluster points $\mathbf{P}_{x_0}^{\infty,\infty}$ of $(\mathbf{P}_{x_0}^{\infty,k})_{k\in\mathbb{N}}$ on $\Omega$ from Step 2. Moreover, the density process $p^k : t \mapsto \mathrm{d}\Lambda_t/\mathrm{d}x$ admits the quantitative estimate

$$\Big\|\|p^k(t)\|_{\mathsf{H}_r^w(\mathbb{R}^d)}\Big\|_{\mathbb{L}^q(\mathbf{P}_{x_0}^{\infty,k})} \leq c_{Thm.\ 3.1}(1 \wedge t)^{-\gamma} \qquad \text{for all } t \in (0,T]\,. \tag{1.25}$$

This estimate is instrumental in improving the basic tightness result on $\Omega$.

**Step 4 — A refined tightness result**    Using the $\mathcal{S}'$-semimartingale property of $\Lambda$ on $\Omega$ under $\mathbf{P}_{x_0}^{\infty,k}$ and the quantitative bound (1.25) we show that

$$\mathbf{P}_{x_0}^{\infty,k}\Big[\big\{\omega \in \Omega \,:\, (0,T] \ni t \mapsto \Lambda_t(\omega) \in \mathbb{L}^r(\mathrm{d}x) \text{ is in } \mathbf{C}\big((0,T];\mathbb{L}^r(\mathrm{d}x)\big)\big\}\Big] = 1\,, \tag{1.26}$$

where $\mathbf{C}((0,T];\mathbb{L}^r(\mathrm{d}x))$ is equipped with the topology of uniform convergence on compacts. This is achieved using an interpolation argument to trade spatial regularity of $\Lambda$ against regularity in time. In fact, we find for the $\mathbf{P}_{x_0}^{\infty,k}$-a.s. $\mathbb{L}^r(\mathrm{d}x)$-valued process $\Lambda$ a locally Hölder-continuous modification in the time-variable $t \in (0,T]$, which we use to deduce our second, improved tightness result: For all $p \geq 1$, the sequence $(\mathbf{P}_{x_0}^{\infty,k})_{k\in\mathbb{N}}$ is tight in $\Omega_{p,1} = \mathbb{L}^p([0,T];\mathbb{L}^1(\mathrm{d}x))$. This leads to Lemma 1.8, above.

**Step 5 — The limiting martingale problem**    The refined tightness result allows us to find $\mathbf{P}_{x_0}^{\infty,\infty}$ and a subsequence $(\mathbf{P}_{x_0}^{\infty,k_m})_{m\in\mathbb{N}}$ such that $\mathbf{P}_{x_0}^{\infty,k_m} \to \mathbf{P}_{x_0}^{\infty,\infty}$ narrowly as $m \to \infty$, both in $\Omega$ and in $\Omega_{p,1}$. By the Skorokhod representation theorem, we can construct random variables $(\hat{\Lambda}^m)_{m\in\mathbb{N}\cup\{\infty\}}$ valued in both $\Omega$ and $\Omega_{p,1}$ on a common probability space $(\hat{\Omega},\hat{\mathcal{F}},\hat{\mathbf{P}})$ such that:

(a) for all $m \in \mathbb{N}\cup\{\infty\}$, we have $\mathrm{Law}_{\mathbf{P}_{x_0}^{\infty,k_m}}(\Lambda) = \mathrm{Law}_{\hat{\mathbf{P}}}(\hat{\Lambda}^m)$.

In addition,

(b) $\hat{\mathbf{P}}$-a.s., $\hat{\Lambda}^m \to \hat{\Lambda}^\infty$ in $\mathbb{L}^p([0,T];\mathbb{L}^1(\mathrm{d}x))$ as $m \to \infty$;

(c) $\hat{\mathbf{P}}$-a.s., $\hat{\Lambda}^m \to \hat{\Lambda}^\infty$ in $\mathbf{C}([0,T];\mathcal{S}')$ as $m \to \infty$.

Moreover, in a key result in Proposition 4.5, we combine properties (a)–(c) and the local $\mathbb{L}^r(\mathrm{d}x)$-continuity of $\Lambda$ from Step 4 to deduce that

(d) for all $1 \leq \xi < r$, we have $\hat{\mathbf{P}}$-a.s., $\hat{\Lambda}^m \to \hat{\Lambda}^\infty$ in $\mathbf{C}((0,T];\mathbb{L}^\xi(\mathrm{d}x))$ as $m \to \infty$.

Property (d) is the type of strong convergence that allows us to considerably weaken the regularity required from the coefficients to still obtain (1.23).

With the representation $(\hat{\Lambda}^m)_{m \in \mathbb{N} \cup \{\infty\}}$, the martingale property (1.21) can, for each $m \in \mathbb{N}$ be recast as $\mathbf{E}^{\hat{\mathbf{P}}}[(\hat{M}_t^m - \hat{M}_s^m)g] = 0$. Moreover, the dynamics in (1.14) allow us to write

$$\hat{M}_t^m - \hat{M}_s^m = \hat{\Lambda}_t^m - \hat{\Lambda}_s^m - \int_s^t \mathscr{A}_u^{\delta_m}(\hat{\Lambda}_u^m)\, \mathrm{d}u \qquad (1.27)$$

for the increments. Let us now preview the definition of $\mathscr{A}_u^{\delta_m}$ given in (2.9) and (2.13), below. We have $\hat{\mathbf{P}}$.a.s. that $\mathscr{A}_u^{\delta_m}(\hat{\Lambda}_t^m)$ is given by

$$\mathcal{S} \ni \phi \mapsto \mathscr{A}_u^{\delta_k}(\hat{\Lambda}_t^m)[\phi] = \hat{\Lambda}_t^m[\mathscr{L}_u^{\delta_k}(\hat{\Lambda}_t^m)\phi] = \int_{\mathbb{R}^d} \left( \mathscr{L}_t^{\delta_k}(\hat{\Lambda}_t^m)\phi \right)(x)\, \hat{\Lambda}_t^m(\mathrm{d}x) \in \mathbb{R}, \qquad (1.28)$$

where

$$\left( \mathscr{L}_u^{\delta_m}(\hat{\Lambda}_t^m)\phi \right)(x) \coloneqq b_u^{\delta_m}\left( x, (\hat{\Lambda}_t^m)^{\delta_m} \right) \cdot \nabla \phi(x) + \frac{1}{2}(\sigma_u \sigma_u^\mathsf{T} + \bar{\sigma}_u \bar{\sigma}_u^\mathsf{T})(x, \hat{\Lambda}_t^m) : \nabla^2 \phi(x). \quad (1.29)$$

In addition, $\mathscr{A}^{\delta_\infty} = \mathscr{A}^0 = \mathscr{A}$. By (1.27) and property (a), we see that in this setting, (1.23) follows from (1.22) once we establish that

$$\lim_{m \to \infty} \left| \mathbf{E}^{\hat{\mathbf{P}}} \left[ \left( \int_s^t \left( \mathscr{A}_u^{\delta_m}(\hat{\Lambda}_u^m)g(\hat{\Lambda}^m) - \mathscr{A}_u(\hat{\Lambda}_u^\infty)g(\hat{\Lambda}^\infty) \right) \mathrm{d}u \right) \right] \right| = 0\, ; \qquad (1.30)$$

see (5.3), below. To prove (1.30), we establish in Section 5 that for $(\mathrm{d}u \otimes \mathrm{d}\hat{\mathbf{P}})$-almost every $(u, \hat{\omega})$, we have the convergence $\mathscr{A}_u^{\delta_{k_m}}(\hat{\Lambda}_u^m)[\phi] \longrightarrow \mathscr{A}_u(\hat{\Lambda}_u^\infty)[\phi]$ for every $\phi \in \mathcal{S}$. It is here that we require the second, refined tightness result, which leads to the random variables $(\hat{\Lambda}^m)_{m \in \mathbb{N} \cup \{\infty\}}$.

**The role of regularity and comparison with the unsmoothed case**     At the heart of our approach are two central ideas, both of which hinge on the emergence of regularity. By taking first the limit $n \to \infty$ and then $k \to \infty$, we manage to leverage this principle. Specifically, the emergence of regularity property established in [8] provides crucial, *uniform* analytic estimates, captured by the following key result from Lemma 5.3:

$$\sup_{k \in \mathbb{N} \cup \{\infty\}} \left\| \|\hat{p}^k(t)\|_{\mathsf{H}_p^w(\mathbb{R}^d)} \right\|_{\mathbb{L}^q(\hat{\mathbf{P}})} \le c(1 \wedge t)^{-\gamma} \qquad \text{for all } t \in (0, T]. \qquad (1.31)$$

This uniform bound enables the two main arguments of our proof.

The first central idea is to use this regularity to obtain a strengthened tightness result. As a consequence, the Skorokhod representation theorem provides not just narrow convergence of $\hat{\Lambda}^k \to \hat{\Lambda}_t^\infty$, as would be the case if we had only tightness in $\Omega$, but instead we are able to deduce the much stronger mode of convergence in (d), above. This strong convergence is a payoff of the uniform estimate in (1.31) and is unattainable from tightness in $\Omega$ alone.

The second central idea explains why this regularity is so effective: it exploits the structure of the map $\mathscr{A}_u(\lambda)[\phi] = \lambda[\mathscr{L}_u(\lambda)\phi]$. This form allows favorable regularity properties of the measure $\lambda$ (specifically, of its density) to compensate for the lack of regularity in the $x$-variable of the drift coefficient $b$, which appears as part of the operator $\mathscr{L}_u(\lambda)$.

This approach should be contrasted with an analysis of the original, unsmoothed particle systems and taking the limit $n \to \infty$ there. For those systems, the empirical measures are purely atomic. Consequently, neither a regularity estimate like (1.31) nor a tightness result in $\Omega_{p,1}$ can hold. Proving the existence of a solution in that setting requires stronger continuity properties of the drift $b$ than those assumed here.

## 2 SMOOTHED MARTINGALE PROBLEM

This section serves a preparatory purpose. In Section 2.1, we recall some of the basic terminology from $\mathcal{S}'$-valued stochastic analysis. In particular, the notion of an *canonical $\mathcal{S}'$-valued martingale problem*, which we recall in Definition 2.3 below, is used repeatedly in this paper. In Section 2.2, starting with the empirical measure flow $\mu^n$ from (1.2) under $\mathbf{P}_{x_0}^{n,k}$, we develop the basic notation used to define for each $n, k \in \mathbb{N}$ the canonical $\mathcal{S}'$-valued martingale problem $\mathrm{MP}(\nu_n, \mathscr{A}^{\delta_k}, \mathscr{Q} + \frac{1}{n}\mathscr{C})$. Finally, in Section 2.3, we show that for each $n, k \in \mathbb{N}$, the empirical measure law $\mathbf{P}_{x_0}^{n,k}$ from (1.8) solves $\mathrm{MP}(\nu_n, \mathscr{A}^{n_k}, \mathscr{Q} + \frac{1}{n}\mathscr{C})$.

### 2.1 Background on $\mathcal{S}'$-valued semimartingales and martingale problems

Martingale problems form a backbone of the weak solution theory for $\mathbb{R}^d$-valued SDEs. This formulation is well known and goes back to a series of works by Stroock and Varadhan [38, 39]. We refer to the classic treatment in the book by the same authors [40, Ch. 6], or the more modern account in Karatzas and Shreve [26, Ch. 5.3].

In this work, we require the martingale problem formulation for processes valued in $\mathcal{S}'$ which we recall below; see e.g. Mikulevicius and Rozovskii [7, Ch. 6.3] for more details. To efficiently develop this theory, we recall first some of the basic structural properties of $\mathcal{S}'$ and basic notions of $\mathcal{S}'$-valued semimartingales. A foundational treatment of $\mathcal{S}'$-valued stochastic calculus is found in Itô [23]. More modern expositions are available in e.g. Kallianpur and Xiong [24] or Dalang and Sanz-Solé [11].

**Structure of $\mathcal{S}'$** Let $(\,\cdot\,;\,\cdot\,)_0$ be the usual inner product on $\mathbb{L}^2(\mathrm{d}x)$. The map $f \mapsto (\,\cdot\,; f)_0$ defines a continuous embedding $\mathcal{S} \hookrightarrow \mathcal{S}'$, which the allows us to identify $\mathcal{S}$ and $\mathbb{L}^2(\mathrm{d}x)$ with subspaces of $\mathcal{S}'$. The Fréchet topology of $\mathcal{S}$ and the strong topology of $\mathcal{S}'$ can be described by a family of Hilbertian seminorms, or equivalently, in terms of a family of Hilbert spaces; see e.g. Simon [36] or Appendix A. In concrete terms, we can consider the family $(\mathscr{H}_p)_{p \in \mathbb{R}}$ of *Fourier–Hermite spaces*, defined as follows. Let $(h_k)_{k \in \mathbb{N}^d}$ be the Hermite functions on $\mathbb{R}^d$. Note that $h_k \in \mathcal{S}$ for each $k \in \mathbb{N}^d$. Hence for $\lambda \in \mathcal{S}'$, the $k$'th Hermite–Fourier coefficient $\lambda_k^{\#} := \lambda[h_k]$ is well defined. For $k \in \mathbb{N}$, let $\langle k \rangle := (1 + |k|^2)^{1/2}$. For each $\lambda \in \mathcal{S}'$, there exists $p \in \mathbb{R}$ such that

$$\|\lambda\|_{\mathscr{H}_p} := \left( \sum_{k \in \mathbb{N}} (\langle k \rangle^{p/d} \lambda_k^{\#})^2 \right)^{1/2} < \infty. \tag{2.1}$$

In fact, if $\lambda \in \mathcal{S} \subsetneq \mathcal{S}'$, then the above sum is finite for any $p \in \mathbb{R}$. With this in mind, we can define a Hilbert-subspace of $\mathcal{S}'$ by considering

$$\mathscr{H}_p := \{ \lambda \in \mathcal{S}' : \|\lambda\|_{\mathscr{H}_p} < \infty \}$$

with the norm $\|\cdot\|_{\mathscr{H}_p}$. For $p \geq 0$, $\mathscr{H}_p$ can be identified with a subset of $\mathbb{L}^2(\mathrm{d}x)$ and thus with a set of functions. Finally, the $\mathbb{L}^2$ duality pairing $f \mapsto (f; \cdot)_0$ extends for any $p \in \mathbb{R}$ to an isometric isomorphism $J : \mathscr{H}_{-p} \to \mathscr{H}_p'$ so that $\mathscr{H}_{-p} \simeq \mathscr{H}_p'$. In sum, the above shows that $\mathcal{S} \hookrightarrow \mathscr{H}_p \hookrightarrow \mathbb{L}^2 \hookrightarrow \mathscr{H}_{-p} \hookrightarrow \mathcal{S}'$ for all $p > 0$. It is a classical fact that $\mathcal{S} = \bigcap_{p \in \mathbb{R}} \mathscr{H}_p$ and $\mathcal{S}' = \bigcup_{p \in \mathbb{R}} \mathscr{H}_p$, where the right-hand sides are respectively endowed with the limit and colimit tolology. This makes precise the statement that the topology of $\mathcal{S}$ and $\mathcal{S}'$ can be described by a family of Hilbert spaces.

**Semimartingales valued in $\mathcal{S}'$**   An $\mathcal{S}'$-valued process $M = (M_t)_{t \in [0,T]}$ is called a *square-integrable, continuous $\mathcal{S}'$-valued martingale* relative to $(\mathbb{F}, \mathbf{P})$ if for each $\phi \in \mathcal{S}$, there is a version of $M[\phi] = (M_t[\phi])_{t \in [0,T]}$ in $\mathscr{M}_2^c(\mathbb{F}, \mathbf{P}; \mathbb{R})$; we write $M \in \mathscr{M}_2^c(\mathbb{F}, \mathbf{P}; \mathcal{S}')$. The process $M$ is called a *continuous $\mathcal{S}'$-valued local martingale* if for each $\phi \in \mathcal{S}$, there is a version of $M[\phi]$ in $\mathscr{M}_{\mathrm{loc}}^c(\mathbb{F}, \mathbf{P}; \mathbb{R})$; we write $M \in \mathscr{M}_{\mathrm{loc}}^c(\mathbb{F}, \mathbf{P}; \mathcal{S}')$.

We call $M = (M_t)_{t \in [0,T]}$ a *square-integrable, continuous $\mathscr{H}_{-p}$-valued martingale* if it is $\mathscr{H}_{-p}$-valued and a martingale with continuous trajectories relative to the $\mathscr{H}_{-p}$-norm, $\mathbf{P}$-a.s., and such that

$$\|M\|_{\mathscr{M}_2^c(\mathbb{F}, \mathbf{P}; \mathscr{H}_{-p})} := \mathbf{E}[\|M_T\|_{\mathscr{H}_{-p}}^2]^{\frac{1}{2}} < \infty \,. \tag{2.2}$$

In this case, we write $M \in \mathscr{M}_2^c(\mathbb{F}, \mathbf{P}; \mathscr{H}_{-p})$. We call $M$ a *continuous $H$-valued local martingale* if there exists a localizing sequence $(\tau^n)_{n \in \mathbb{N}}$ such that $M_{\cdot \wedge \tau^n} \in \mathscr{M}_2^c(\mathbb{F}, \mathbf{P}; \mathscr{H}_{-p})$. In this case we write $M \in \mathscr{M}_{\mathrm{loc}}^c(\mathbb{F}, \mathbf{P}; \mathscr{H}_{-p})$.

We should emphasize the asymmetry in the requirements placed on the localizing sequence. In the $\mathcal{S}'$-valued case, this sequence may depend on $\phi$, in the $\mathscr{H}_{-p}$-valued case, it must not.

**Definition 2.1** | Let $\Lambda$ be a continuous $\mathcal{S}'$-valued process. We call $\Lambda$ a *continuous $\mathcal{S}'$-valued $(\mathbb{F}, \mathbf{P})$-semimartingale* if there are a random variable $\Lambda_0$ and processes $A = (A_t)_{t \in [0,T]}$, $M = (M_t)_{t \in [0,T]}$, all $\mathcal{S}'$-valued, such that

$$\Lambda_t = \Lambda_0 + A_t + M_t \quad \text{in } \mathcal{S}' \text{ for all } t \in [0,T]$$

and the following conditions are satisfied:

*1)* $\Lambda_0$ is $\mathcal{F}_0$-measurable.

*2)* $A$ is $\mathbb{F}$-adapted, $A_0 = 0$ and for all $\phi \in \mathcal{S}$, the process given by $A[\phi] = (A_t[\phi])_{t \in [0,T]}$ has continuous paths of finite variation.

*3)* $M \in \mathscr{M}_{\mathrm{loc}}^c(\mathbb{F}, \mathbf{P}; \mathcal{S}')$ and $M_0 = 0$. Equivalently, $M[\phi] = (M_t[\phi])_{t \in [0,T]} \in \mathscr{M}_{\mathrm{loc}}^c(\mathbb{F}, \mathbf{P}; \mathbb{R})$ for all $\phi \in \mathcal{S}$ and $M_0 = 0$.

Let $\Lambda$ be an $\mathscr{H}_{-p}$-valued process with continuous trajectories, $\mathbf{P}$-a.s. We call $\Lambda$ a *continuous $\mathscr{H}_{-p}$-valued $(\mathbb{F}, \mathbf{P})$-semimartingale* if there are an $\mathscr{H}_{-p}$-valued random variable $\Lambda_0$ and $\mathscr{H}_{-p}$-valued processes $A$ and $M$ such that

$$\Lambda_t = \Lambda_0 + A_t + M_t \quad \text{in } \mathscr{H}_{-p} \text{ for all } t \in [0,T] \,,$$

and the following conditions are satisfied: (a) $\Lambda_0$ is $\mathcal{F}_0$-measurable, (b) $A$ is $\mathbb{F}$-adapted, $A_0 = 0$ and $A$ has continuous paths of finite variation relative to the $\mathscr{H}_{-p}$-norm, (c) $M \in \mathscr{M}_{\mathrm{loc}}^c(\mathbb{F}, \mathbf{P}; \mathscr{H}_{-p})$ and $M_0 = 0$.

**Lemma 2.2** | *Let* $\tilde{\Lambda} = (\tilde{\Lambda}_t)_{t \in [0,T]}$ *be a cylindrical process on* $\mathcal{S}$, *i.e.* $\tilde{\Lambda}_t : \mathcal{S} \to \mathbb{L}^0(\Omega, \mathcal{F}, \mathbf{P})$ *is a linear map for each* $t \in [0,T]$. *Suppose that* $\tilde{\Lambda}_t : \mathcal{S} \to \mathbb{L}^0(\Omega, \mathcal{F}, \mathbb{P})$ *is continuous for each* $t \in [0,T]$ *and that for each* $\phi \in \mathcal{S}$, *the process* $\tilde{\Lambda}[\phi] = (\tilde{\Lambda}_t[\phi])_{t \in [0,T]}$ *is a real-valued* $(\mathbb{F}, \mathbf{P})$-*semimartingale with continuous trajectories and canonical decomposition*

$$\tilde{\Lambda}[\phi] = \Lambda_0^\phi + A^\phi + M^\phi,$$

*where* $\Lambda_0^\phi$ *is an* $\mathcal{F}_0$-*measurable real-valued random variable,* $A^\phi$ *is a real-valued* $\mathbb{F}$-*adapted and continuous process of finite variation, and* $M^\phi \in \mathcal{M}_{\mathrm{loc}}^c(\mathbb{F}, \mathbf{P}; \mathbb{R})$. *Then:*

1) *There is a continuous* $\mathcal{S}'$-*valued* $(\mathbb{F}, \mathbf{P})$-*semimartingale* $\Lambda$ *satisfying* $\Lambda_0[\phi] = \Lambda_0^\phi$, $A[\phi] = A^\phi$, $M[\phi] = M^\phi$ *for all* $\phi \in \mathcal{S}$, $\mathbf{P}$-*a.s.*

2) *If for each* $\phi \in \mathcal{S}$, *the real-valued semimartingale* $\Lambda[\phi] = \Lambda_0[\phi] + A[\phi] + M[\phi]$ *satisfies in addition that*

$$\left\| \int_0^T |\mathrm{d}A_t[\phi]| \right\|_{\mathbb{L}^2(\Omega, \mathcal{F}, \mathbf{P})} + \|M[\phi]\|_{\mathcal{M}_2^c(\mathbb{F}, \mathbf{P}; \mathbb{R})} < \infty,$$

*then there exists* $p_0 > 0$ *such that for all* $p \geq p_0$, *we have* $\mathbf{P}$-*a.s. that* $\Lambda, A, M \in \mathbf{C}([0,T]; \mathscr{H}_{-p})$ *and* $M \in \mathcal{M}_2^c(\mathbb{F}, \mathbf{P}; \mathscr{H}_{-p})$.

**Proof** By [29, Thm. 1], we find an $\mathcal{S}'$-valued continuous version $\Lambda$ of $\tilde{\Lambda}$, and part (i) of the lemma then follows from Pérez-Abreu [30, Thm. 1]. Part (ii) follows from e.g. Fonseca-Mora [14, Cor. 3.8]. □

**Tensor quadratic variation** We call a map $Q : \mathcal{S} \times \mathcal{S} \to \mathbb{R}$ *symmetric* if $Q[\phi, \psi] = Q[\psi, \phi]$ for all $\phi, \psi \in \mathcal{S}$, and *nonnegative* if $Q[\phi, \phi] \geq 0$ for all $\phi \in \mathcal{S}$. Denote by $\mathcal{L}_{\mathrm{bl}}^+(\mathcal{S})$ the space of bilinear maps $\mathcal{S} \times \mathcal{S} \to \mathbb{R}$ which are symmetric and nonnegative. This space carries the topology of uniform convergence on products of bounded sets in $\mathcal{S}$; see Huang and Yan [21, Sec. 3.3] for details.

A process $Q : [0, T] \times \Omega \to \mathcal{L}_{\mathrm{bl}}^+(\mathcal{S})$ is called *weakly* $\mathbb{F}$-*progressive* if for all $\phi, \psi \in \mathcal{S}$, the $\mathbb{R}$-valued process $Q[\phi, \psi]$ is $\mathbb{F}$-progressive, and *strongly* $\mathbb{F}$-*progressive* if the process $Q$ is $\mathbb{F}$-progressive when $\mathcal{L}_{\mathrm{bl}}^+(\mathcal{S}')$ is equipped with its Borel-$\sigma$-algebra. This is analogous to the weak and strong measurability in known from infinite-dimensional analysis because an element of $\mathcal{L}_{\mathrm{bl}}^+(\mathcal{S})$ can equivalently be seen as a continuous linear map $\mathcal{S} \to \mathcal{S}'$; see Huang and Yan [21, Thm. 3.17] or Simon [36, Sec. 4].

We say that $M \in \mathcal{M}_2^c(\mathbb{F}, \mathbf{P}; \mathcal{S}')$ has the *(tensor) quadratic variation* process $\int_0^t Q_s \, \mathrm{d}s$ with *(tensor) covariance operator* $Q : [0, T] \times \Omega \to \mathcal{L}_{\mathrm{bl}}^+(\mathcal{S})$ if $Q$ is weakly $\mathbb{F}$-progressive and

$$\left( M_t[\phi] M_t[\psi] - \int_0^t Q_s[\phi, \psi] \, \mathrm{d}s \right)_{t \in [0,T]} \tag{2.3}$$

is in $\mathcal{M}^c(\mathbb{F}, \mathbf{P}; \mathbb{R})$ for all $\phi, \psi \in \mathcal{S}$. It is well known that if $M \in \mathcal{M}_2^c(\mathbb{F}, \mathbf{P}; \mathcal{S}')$, then the tensor quadratic variation exists and is uniquely characterized by (2.3), and we have

$$\left( \int_0^t Q_s \, \mathrm{d}s \right)[\phi, \psi] = \int_0^t Q_s[\phi, \psi] \, \mathrm{d}s \tag{2.4}$$

for all $t \in [0, T]$ and $\phi, \psi \in \mathcal{S}$, $\mathbf{P}$-a.s. Moreover, $Q$ can be taken to be strongly $\mathbb{F}$-progressive so that the integral $\int Q_t \, \mathrm{d}t$ above is understood in the Bochner sense. Indeed,

Lemma 2.2, we may assume that there is a version (still denoted by $M$) in $\mathscr{M}_2^c(\mathscr{H}_{-p})$. But in the setting of Hilbert spaces, the result is classical, see Gawarecki and Mandrekar [16, Lem. 2.1], and $(Q_t)_{t \in [0,T]}$ is Bochner-integrable.

We generally denote the *tensor variation process* by $\langle\!\langle M \rangle\!\rangle$, that is, $\langle\!\langle M \rangle\!\rangle$ is the unique $\mathcal{L}_{\mathrm{bl}}^+(\mathcal{S})$-valued, progressive process such that

$$M[\phi]M[\psi] - \langle\!\langle M \rangle\!\rangle[\phi, \psi] \in \mathscr{M}^c(\mathbb{R}) \qquad \text{for all } \phi, \psi \in \mathcal{S}.$$

**Canonical $\mathcal{S}'$-valued martingale problems**   We take as pre-specified the measurable space $(\Omega, \mathcal{F})$ together with the right-continuous version of the canonical filtration $\mathbb{F} = (\mathcal{F}_t)_{t \in [0,T]}$ and the $\mathcal{S}'$-valued canonical process $\Lambda = (\Lambda_t)_{t \in [0,T]}$.

Assume we are then given a triplet $(\nu_0, A, Q)$, consisting of

- a probability measure $\nu_0$ on $\mathcal{B}(\mathbb{R}^d)$;

- an $\mathbb{F}$-progressive process $A : [0, T] \times \Omega \to \mathcal{S}' \cup \{+\infty\}$, i.e. $A[\phi]$ is $\mathbb{F}$-progressive for each $\phi \in \mathcal{S}$, and $A_0 = 0$;

- an $\mathbb{F}$-progressive process $Q : [0, T] \times \Omega \to \mathcal{L}_{\mathrm{bl}}^+(\mathcal{S})$, i.e. $Q[\phi, \psi]$ is $\mathbb{F}$-progressive for each $\phi, \psi \in \mathcal{S}$.

We say that $(\nu_0, A, Q)$ defines the *data of a martingale problem*, and write $\mathrm{MP}(\nu_0, A, Q)$.

**Definition 2.3** | A probability measure $\mathbf{P}$ on $(\Omega, \mathcal{F})$ is said to be a *solution of* $\mathrm{MP}(\nu_0, A, Q)$ if the following hold:
   *1)* We have that $\mathrm{Law}_{\mathbf{P}}(\Lambda_0) = \nu_0$.
   *2)* For each $\phi \in \mathcal{S}$, we have that $\mathbf{P}[\int_0^T |A_s[\phi]| \, \mathrm{d}s < \infty]$.
   *3)* If $M = (M_t)_{t \in [0,T]}$ is defined by

$$M_t := \Lambda_t - \Lambda_0 - \int_0^t A_s \, \mathrm{d}s \quad \text{for } t \in [0, T],$$

then $M \in \mathscr{M}_{\mathrm{loc}}^c(\mathbb{F}, \mathbf{P}; \mathcal{S}')$.
   *4)* $\langle M[\phi], M[\psi] \rangle = \int Q_t[\phi, \psi] \, \mathrm{d}t$ $\mathbf{P}$-a.s. whenever $\phi, \psi \in \mathcal{S}$.

The reason we allow $A$ to take the value $+\infty$ is because this lets us construct martingale problems whose solutions, provided they exist, are prohibited of putting mass on certain measurable subsets of $\mathcal{S}'$ via condition 1) above; see Mikulevicius and Rozovskii [7, Def. 6.3.3] for an alternative but equivalent formulation.

**Remark 2.4** | If $\mathbf{P}$ solves $\mathrm{MP}(\nu_0, A, Q)$ and if $M \in \mathscr{M}_2^c(\mathbb{F}, \mathbf{P}; \mathcal{S}')$, then part *4)* in Definition 2.3 becomes
   *4')* $\langle\!\langle M \rangle\!\rangle = \int Q_t \, \mathrm{d}t$ $\mathbf{P}$-a.s.

**Remark 2.5** | There is an equivalent way to formulate martingale problems in terms of a second order differential operator acting on Fréchet differentiable functions $\mathcal{S}' \to \mathbb{R}$. We refer to e.g. Kallianpur et al. [25] for this approach. We avoid it here because it would require additional notation.

## 2.2 Data for the martingale problem

Consider again the particle system (1.1) and (1.2). A central object in this section is the *empirical measure flow*, i.e.,

$$\mu_t^n = \frac{1}{n}\sum_{i=1}^n \delta_{X_t^{i,n}} \quad \text{with } t \in [0,T].\tag{2.5}$$

By Lemma 1.6, the canonical process $\Lambda$ on $(\Omega, \mathcal{F}, \mathbf{P}_{x_0}^{n,k})$ represents the dynamics of the empirical measure in the sense that for all $k, n \in \mathbb{N}$, $\mathrm{Law}_{\mathbb{P}_{x_0}^{n,\delta_k}}(\mu^n) = \mathrm{Law}_{\mathbf{P}_{x_0}^{n,k}}(\Lambda)$. As a steppingstone to describe the dynamics of $\Lambda$ under $\mathbf{P}_{x_0}^{n,k}$, we now develop the essential notation based on the dynamics of $n$-particle systems. This notation is used repeatedly. We introduce the maps

$$\mathscr{L}_t^\delta : \mathbf{M}_1^+(\mathbb{R}^d) \times \mathbb{S}(\mathbb{R}^d) \to \mathbb{L}^1(\mathbb{R}^d;\mathbb{R}) \cap \mathbb{L}^\infty(\mathbb{R}^d;\mathbb{R}),\tag{2.6}$$

$$\mathscr{R}_t : \mathbf{M}_1^+(\mathbb{R}^d) \times \mathbb{S}(\mathbb{R}^d) \to \mathbb{L}^1(\mathbb{R}^d;\mathbb{R}^d) \cap \mathbf{C}_b(\mathbb{R}^d;\mathbb{R}^d),\tag{2.7}$$

$$\mathscr{U}_t : \mathbf{M}_1^+(\mathbb{R}^d) \times \mathbb{S}(\mathbb{R}^d) \to \mathbb{L}^1(\mathbb{R}^d;\mathbb{R}^d) \cap \mathbf{C}_b(\mathbb{R}^d;\mathbb{R}^d),\tag{2.8}$$

defined as follows. We set

$$(\lambda,\phi) \mapsto \mathscr{L}_t^\delta(\lambda)\phi, \qquad (\lambda,\phi) \mapsto \mathscr{R}_t(\lambda)\phi, \qquad (\lambda,\phi) \mapsto \mathscr{U}_t(\lambda)\phi,$$

where

$$\left(\mathscr{L}_t^\delta(\lambda)\phi\right)(x) \coloneqq b_t^\delta(x,\lambda^\delta) \cdot \nabla\phi(x) + \frac{1}{2}a_t(x,\lambda) : \nabla^2\phi(x),\tag{2.9}$$

$$\left(\mathscr{R}_t(\lambda)\phi\right)(x) \coloneqq \bar{\sigma}_t^\mathsf{T}(x,\lambda)\nabla\phi(x),\tag{2.10}$$

$$\left(\mathscr{U}_t(\lambda)\phi\right)(x) \coloneqq \sigma_t^\mathsf{T}(x,\lambda)\nabla\phi(x).\tag{2.11}$$

As usual, for $\delta = 0$, we simply write $\mathscr{L} \coloneqq \mathscr{L}^0$. With $\lambda^\delta$ from (1.3), the notation specifying the function (2.9) is explicitly given by

$$b_t^\delta(x,\lambda^\delta) \cdot \nabla\phi(x) = \sum_{j=1}^d (b_t^\delta)^j(x,\lambda^\delta)\partial_j\phi(x),$$

$$a_t(x,\lambda) : \nabla^2\phi(x) = \sum_{i,j=1}^d a_t^{ij}(x,\lambda)\partial_{ij}\phi(x),$$

with

$$a_t^{ij}(x,\lambda) \coloneqq (\sigma_t\sigma_t^\mathsf{T} + \bar{\sigma}_t\bar{\sigma}_t^\mathsf{T})^{ij}(x,\lambda) = \sum_{k=1}^d (\sigma_t^{ik}\sigma_t^{kj} + \bar{\sigma}_t^{ik}\bar{\sigma}_t^{kj})(x,\lambda),\tag{2.12}$$

the notation specifying (2.10) by

$$\bar{\sigma}_t^\mathsf{T}(x,\lambda)\nabla\phi(x) = \sum_{j=1}^d \bar{\sigma}_t^{j\cdot}(x,\lambda)\partial_j\phi(x),$$

and the notation specifying (2.11) by

$$\sigma_t^{\mathsf{T}}(x,\lambda)\nabla\phi(x) = \sum_{j=1}^{d} \sigma_t^{j\cdot}(x,\lambda)\partial_j\phi(x)\,.$$

The fact that for each $(t,\lambda,\phi) \in [0,T] \times \mathbf{M}_1^+ \times \mathcal{S}$, the function $x \mapsto (\mathscr{L}_t^\delta(\lambda)\phi)(x)$ is an element of $\mathbb{L}^1(\mathbb{R}^d;\mathbb{R}) \cap \mathbb{L}^\infty(\mathbb{R}^d;\mathbb{R})$ follows from (2.9) together with observing that $b_t^\delta(\,\cdot\,,\lambda^\delta)$, $\sigma_t(\,\cdot\,,\lambda)$ and $\bar\sigma_t(\,\cdot\,,\lambda)$ are bounded, and that $\partial_i\phi$ and $\partial_{ij}\phi$ are for any $i,j \in [n]$ in $\mathcal{S}$, and thus also elements of $\mathbb{L}^s(\mathbb{R}^d;\mathbb{R})$ for any $s \in [1,\infty]$. Similarly, the fact that $x \mapsto (\mathscr{R}_t(\lambda)\phi)(x)$ and $x \mapsto (\mathscr{U}_t(\lambda)\phi)(x)$ are in $\mathbb{L}^1(\mathbb{R}^d;\mathbb{R}^d) \cap \mathbf{C}_b(\mathbb{R}^d;\mathbb{R}^d)$ follows from (2.10) and (2.11) along with Assumption 1.3 by which $x \mapsto \sigma_t(x,\lambda)$ and $x \mapsto \bar\sigma_t(x,\lambda)$ are continuous and bounded, together with the fact that $\partial_i\phi$ is for any $i \in [n]$ in $\mathcal{S}$, and thus also in $\mathbf{C}_b(\mathbb{R}^d;\mathbb{R})$ and in $\mathbb{L}^1(\mathbb{R}^d;\mathbb{R})$.

We often find it convenient to simplify the notation and write

$$\mathscr{L}_t^\delta(\lambda)\phi(x) \coloneqq \big(\mathscr{L}_t^\delta(\lambda)\phi\big)(x)\,,$$
$$\mathscr{R}_t(\lambda)\phi(x) \coloneqq \big(\mathscr{R}_t(\lambda)\phi\big)(x)\,,$$
$$\mathscr{U}_t(\lambda)\phi(x) \coloneqq \big(\mathscr{U}_t(\lambda)\phi\big)(x)\,.$$

This notation emphasizes that for any $(t,\lambda) \in [0,T] \times \mathbf{M}_1^+$, $\mathscr{L}_t^\delta(\lambda)$, $\mathscr{R}_t(\lambda)$ and $\mathscr{U}_t(\lambda)$ are *differential operators* acting on functions in $\mathcal{S}$.

We point out that the maps $\mathscr{L}_t^\delta$, $\mathscr{R}_t$ and $\mathscr{U}_t$ arise naturally via Itô's formula; see e.g. [41] for a details, or [10, Sec. 2.2] for a derivation using the same notation as we use here.

## 2.3 Martingale problem for large-system limits of the smoothed particle systems

The notation introduced in Section 2.2 lets us formulate for each $n \in \mathbb{N}$ a martingale problem in the sense of Section 2.1.

We find it convenient to introduce more succinct notation for the key objects which appear as the *characteristics* of the martingale problem we formulate. For this, we now return to (1.11)–(1.13) which, to recall, are given for each $\delta \geq 0$ and $t \in [0,T]$ by

$$\mathscr{A}_t^\delta : \mathcal{S}' \to \mathcal{S}' \cup \{\infty\}\,,$$

with $\lambda \mapsto \mathscr{A}_t^\delta(\lambda)$, and

$$\mathscr{Q}_t : \mathcal{S}' \to \mathcal{L}_{\mathrm{bl}}^+(\mathcal{S})\,,$$
$$\mathscr{C}_t : \mathcal{S}' \to \mathcal{L}_{\mathrm{bl}}^+(\mathcal{S})\,,$$

with $\lambda \mapsto \mathscr{Q}_t(\lambda)$ and $\lambda \mapsto \mathscr{C}_t(\lambda)$, respectively, where we recall that $\mathcal{L}_{\mathrm{bl}}^+(\mathcal{S})$ denotes the space of continuous, nonnegative definite and symmetric bilinear forms $\mathcal{S} \times \mathcal{S} \to \mathbb{R}$; see Section 2.1 for the definition. These maps are explicitly given as follows. We first recall the notation $\lambda[f] = \int f(x)\,\lambda(\mathrm{d}x)$, where the integral is understood to be componentwise if $f$ is an $\mathbb{R}^d$-valued bounded function. For each $\lambda \in \mathbf{M}_1^+$ and $\phi,\psi \in \mathcal{S}$, we set

$$\mathscr{A}_t^\delta(\lambda)[\phi] \coloneqq \lambda[\mathscr{L}_t^\delta(\lambda)\phi]\,, \tag{2.13}$$

and $\mathscr{A}_t^\delta(\lambda) = \infty$ whenever $\lambda \in \mathcal{S}' \backslash \mathbf{M}_1^+$. Moreover, if $\lambda \in \mathbf{M}_1^+$, we set

$$\mathscr{Q}_t(\lambda)[\phi, \psi] := \lambda[\mathscr{R}_t(\lambda)\phi] \cdot \lambda[\mathscr{R}_t(\lambda)\psi], \tag{2.14}$$

$$\mathscr{C}_t(\lambda)[\phi, \psi] := \lambda\Big[\big(\mathscr{U}_t(\lambda)\phi\big) \cdot \big(\mathscr{U}_t(\lambda)\phi\big)\Big], \tag{2.15}$$

and $\mathscr{Q}_t(\lambda)[\phi, \psi] = 0$ and $\mathscr{C}_t(\lambda)[\phi, \psi] = 0$ if $\lambda \in \mathcal{S}' \backslash \mathbf{M}_1^+$. On the space $(\Omega, \mathcal{F})$ with the filtration $\mathbb{F} = (\mathcal{F}_t)_{t \in [0,T]}$ and the canonical process $\Lambda = (\Lambda_t)_{t \in [0,T]}$ introduced in Section 1.2, and with the fixed sequence $(\delta_k)_{k \in \mathbb{N}}$ from (1.16), this defines for each $k, n \in \mathbb{N}$ a martingale problem $\mathrm{MP}(\nu_n, \mathscr{A}^{\delta_k}, \mathscr{Q} + \frac{1}{n}\mathscr{C})$ in the sense of Section 2.1 with $\nu_0 = \nu$, $A = \mathscr{A}^{\delta_k}$ and $Q = \mathscr{Q} + \frac{1}{n}\mathscr{C}$.

**Proposition 2.6** | *Under Assumptions 1.1, 1.2 and 1.3, the measure $\mathbf{P}_{x_0}^{n,k}$ defined by (1.8) solves $\mathrm{MP}(\nu_n, \mathscr{A}^{\delta_k}, \mathscr{Q} + \frac{1}{n}\mathscr{C})$ in the sense of Definition 2.3. In particular, $\Lambda$ on $(\Omega, \mathcal{F}, \mathbf{P}_{x_0}^{n,k})$ possesses the canonical semimartingale decomposition*

$$\Lambda_t = \Lambda_0 + \int_0^t \mathscr{A}_s^{\delta_k}(\Lambda_s)\,\mathrm{d}s + M_t^{\delta_k}, \qquad t \in [0, T], \tag{2.16}$$

*with $M^{\delta_k}$ in $\mathscr{M}_2^c(\mathbb{F}, \mathbf{P}_{x_0}^{n,k}; \mathcal{S}')$ and with tensor quadratic variation*

$$\langle\!\langle M^{\delta_k} \rangle\!\rangle_t = \int_0^t \mathscr{Q}_s(\Lambda_s)\,\mathrm{d}s + \frac{1}{n}\int_0^t \mathscr{C}_s(\Lambda_s)\,\mathrm{d}s, \qquad t \in [0, T]. \tag{2.17}$$

**Proof**    This follows from Lemma 1.4 and [10, Prop. 2.6]. □

Starting point for our strategy to prove Theorem 1.9 is the following simple result.

**Lemma 2.7** | *If $\delta > 0$ and $b$ satisfies Assumption 1.7, then $b^\delta$ defined in (1.5) also satisfies Assumption 1.7. In addition, $b_t^\delta$ satisfies:*

$(\mathrm{C_s})$ *If $\mu_\infty \in \mathbf{M}_1^+$ is absolutely continuous with respect to Lebesgue measure and $(\mu_k)_{k \in \mathbb{N}}$ is an arbitrary sequence in $\mathbf{M}_1^+$ with $\mu_k \to \mu_\infty$ narrowly as $k \to \infty$, then*

$$\lim_{k \to \infty} \sup_{x \in K} |b_t^\delta(x, \mu_k^\delta) - b_t^\delta(x, \mu_\infty^\delta)| = 0$$

*for any compact subset $K \subseteq \mathbb{R}^d$ and $t \in [0, T]$.*

**Proof**    See Appendix B. □

This motivates the two-step procedure (1.16), (1.18). Indeed, note that the continuity property $(\mathrm{C_s})$ in Lemma 2.7 is stronger than the property $(\mathrm{C_w})$ that we ask for in Assumption 1.7. In fact, the property $(\mathrm{C_s})$ is the type of continuity imposed in [8, Thm. 1.6]. We can therefore apply the results of that paper to show that any cluster point $\mathbf{P}_{x_0}^{\infty,k}$ of $(\mathbf{P}_{x_0}^{n,k})_{k \in \mathbb{N}}$ chosen in (1.18) solves the martingale problem $\mathrm{MP}(\nu, \mathscr{A}^{\delta_k}, \mathscr{Q})$, which is obtained for each fixed $k \in \mathbb{N}$ as the formal limit of the sequence $(\mathrm{MP}(\nu_n, \mathscr{A}^{\delta_k}, \mathscr{Q} + \frac{1}{n}\mathscr{C}))_{n \in \mathbb{N}}$ as $n \to \infty$. This is summarized in the following result.

**Proposition 2.8** | *Under Assumptions 1.1, 1.3 and 1.7, any narrow cluster point $\mathbf{P}_{x_0}^{\infty,k}$ of the sequence $(\mathbf{P}_{x_0}^{n,k})_{n \in \mathbb{N}}$ solves the martingale problem $\mathrm{MP}(\nu, \mathscr{A}^{\delta_k}, \mathscr{Q})$. In particular, the coordinate process $\Lambda$ on $(\Omega, \mathscr{F}, \mathbf{P}_{x_0}^{\infty,k})$ possesses the $\mathcal{S}'$-valued semimartingale decomposition*

$$\Lambda_t = \Lambda_0 + \int_0^t \mathscr{A}_s^{\delta_k}(\Lambda_s)\,\mathrm{d}s + M_t^{\delta_k}, \qquad t \in [0, T], \tag{2.18}$$

*where $M^{\delta_k} = (M_t^{\delta_k})_{t \in [0,T]}$ is in $\mathscr{M}_2^c(\mathbb{F}, \mathbf{P}_{x_0}^{\infty}; \mathcal{S}')$ with quadratic variation*

$$\langle\!\langle M^{\delta_k} \rangle\!\rangle_t = \int_0^t \mathscr{Q}_s(\Lambda_s)\,\mathrm{d}s, \qquad t \in [0, T]. \tag{2.19}$$

**Proof** This follows from Lemma 2.7 and [8, Thm. 1.6]. □

Via the two-step procedure in (1.16) and (1.18), we thus not only produce a sequence of measures $(\mathbf{P}_{x_0}^{\infty,k})_{k \in \mathbb{N}}$, but, by the previous result, a sequence of solutions of the martingale problems $(\mathrm{MP}(\nu, \mathscr{A}^{\delta_k}, \mathscr{Q}))_{k \in \mathbb{N}}$.

## 2.4 A basic tightness result

We next show that $(\mathbf{P}_{x_0}^{\infty,k})_{k \in \mathbb{N}}$ is tight in $\Omega$.

**Lemma 2.9** | *Under Assumptions 1.1, 1.2 and 1.3, the sequence $(\mathbf{P}_{x_0}^{\infty,k})_{k \in \mathbb{N}}$ is tight in $\mathbf{M}_1^+(\Omega)$ and there exists $p_1 > 0$ such that $\mathbf{P}_{x_0}^{\infty,k}[\mathbf{C}([0,T]; \mathscr{H}_{-p})] = 1$ for all $k \in \mathbb{N}$ and $p > p_1$. Moreover,*

$$\mathbf{P}_{x_0}^{\infty,\infty}\Big[\mathbf{C}\big([0,T];\, \mathscr{H}_{-p} \cap \mathbf{M}_1^+(\mathbb{R}^d)\big)\Big] = 1 \tag{2.20}$$

*for any any cluster point $\mathbf{P}_{x_0}^{\infty,\infty}$ of $(\mathbf{P}_{x_0}^{\infty,k})_{k \in \mathbb{N}}$ and any $p > p_1$. Hence $\mathbf{P}_{x_0}^{\infty,\infty}$ concentrates its mass on the set of probability-measure-valued processes.*

**Proof** See Appendix B. □

The previous lemma gives us at least one cluster point $\mathbf{P}_{x_0}^{\infty,\infty}$, and this is our solution candidate for the martingale problem $\mathrm{MP}(\nu, \mathscr{A}, \mathscr{Q})$, which is obtained as the formal limit of the sequence $(\mathrm{MP}(\nu, \mathscr{A}^{\delta_k}, \mathscr{Q}))_{k \in \mathbb{N}}$ as $k \to \infty$. To establish that $\mathbf{P}_{x_0}^{\infty,\infty}$ is indeed a solution, we require certain regularity estimates that we establish in the next section. Before doing this, however, we collect here for later reference a useful estimate.

**Remark 2.10** | Let $\phi \in \mathcal{S}$ and set $\|\phi\|_m^* = \max_{|\alpha| \le m} \sup_{x \in \mathbb{R}^d} |(1 + |x|^2)^{m/2} \partial^\alpha \phi(x)|$ for $m \in \mathbb{N}_0$; see (A.2) for more on the definition of the family of seminorms $\|\cdot\|_m^*$ with $m \in \mathbb{N}$. Recalling the convention from Section 2.2 that $\mathscr{L}^0 = \mathscr{L}$ and in view of $\|b^\delta\|_\infty \le \|b\|_\infty$ by part 3) of Lemma B.2, we obtain under Assumptions 1.2 and 1.3, for all $\lambda, \lambda' \in \mathbf{M}_1^+$, $t \in [0, T]$ and $\delta \ge 0$ the bound

$$|\lambda[\mathscr{L}_t^\delta(\lambda')\phi]| \le \|b^\delta\|_\infty \|\phi\|_1^* + \frac{1}{2}(\|\sigma\|_\infty^2 + \|\bar\sigma\|_\infty^2)\|\phi\|_2^* \le c_{b,\sigma,\bar\sigma}\|\phi\|_2^* \tag{2.21}$$

with $c_{b,\sigma,\bar\sigma} := \|b\|_\infty + \frac{1}{2}(\|\sigma\|_\infty^2 + \|\bar\sigma\|_\infty^2) < \infty$. This is an estimate that we rely on repeatedly, for instance in the following form. By (1.9), $\mathbf{P}_{x_0}^{\infty,k}$ puts for each $k \in \mathbb{N}$ all mass on the set $\mathbf{C}([0,T];\mathcal{S}'\cap\mathbf{M}_1^+)$. Under Assumptions 1.3 and 1.7, we can use Proposition 2.8 and (2.13), we can evaluate (2.18) under $\mathbf{P}_{x_0}^{\infty,k}$ to get for each $t\in[0,T]$ that

$$M_t^{\delta_k}[\phi] = \Lambda_t[\phi] - \Lambda_0[\phi] - \left(\int_s^t \mathscr{A}_u^{\delta_k}(\Lambda_u)\right)[\phi]\,\mathrm{d}u = \Lambda_t[\phi] - \Lambda_0[\phi] - \int_s^t \Lambda_u[\mathscr{L}_u^{\delta_k}(\Lambda_u)\phi]\,\mathrm{d}u\,.$$

The estimate (2.21) now shows that we can bound the process $M^{\delta_k}[\phi] := (M_t^{\delta_k}[\phi])_{t\in[0,T]}$ uniformly on the set $[0,T]\times\mathbf{C}([0,T];\mathcal{S}'\cap\mathbf{M}_1^+)$ by

$$\sup_{t\in[0,T]}|M_t^{\delta_k}[\phi]| \le 2\|\phi\|_0^* + Tc_{b,\sigma,\bar{\sigma}}\|\phi\|_2^* \le (2+Tc_{b,\sigma,\bar{\sigma}})\|\phi\|_2^*. \tag{2.22}$$

In fact, using (2.2) and Lemma A.2 we show via (2.22) that for an appropriate $p>0$ we get that

$$\|M_T^{\delta_k}\|_{\mathscr{H}_{-p}}^2 \le c(2+Tc_{b,\sigma,\bar{\sigma}})\,. \tag{2.23}$$

Together with the regularization result in Lemma 2.2 this implies that $M\in\mathscr{M}_2^c(\mathbf{P}_{x_0}^{\infty,k};\mathscr{H}_{-p})$.

Since by (2.20) in Lemma 2.9 also any cluster point $\mathbf{P}_{x_0}^{\infty,\infty}$ of $(\mathbf{P}_{x_0}^{\infty,k})_{k\in\mathbb{N}}$, puts all mass on the set $\mathbf{C}([0,T];\mathcal{S}'\cap\mathbf{M}_1^+)$, the estimates in (2.21), (2.22) and (2.23) remain valid for $\mathbf{P}_{x_0}^{\infty,\infty}$ and with $\delta = \delta_\infty = 0$. Thus if $M[\phi] = M^0[\phi]$ happens to have the $(\mathbb{F},\mathbf{P}_{x_0}^{\infty,\infty})$-local martingale property, then $M[\phi]$ is also a bounded and thus a true martingale.

# 3 UNIFORM REGULARITY ESTIMATES

In (1.18), we fixed a sequence $(\mathbf{P}_{x_0}^{\infty,k})_{k\in\mathbb{N}}$ by choosing for each $k\in\mathbb{N}$ a cluster point of the empirical measure flow laws of large-system limits with smoothed coefficients, i.e. a cluster point of $(\mathbf{P}_{x_0}^{n,k})_{n\in\mathbb{N}}$. By Lemma 2.9, the sequence $(\mathbf{P}_{x_0}^{\infty,k})_{k\in\mathbb{N}}$ is tight.

In this section we show that the canonical process $\Lambda$ under $(\mathbf{P}_{x_0}^{\infty,k})_{k\in\mathbb{N}}$ and any cluster point $\mathbf{P}_{x_0}^{\infty,\infty}$ satisfies a uniform regularity estimate. To achieve this, we use the *emergence of regularity* property that we obtained in [8].

## 3.1 Revisiting the emergence of regularity

Let $\mathsf{H}_r^s(\mathbb{R}^d)$ for $r\in(1,\infty)$ and $s\in\mathbb{R}$ be the Bessel potential spaces; see e.g. Bergh and Löfström [5, Ch. 6] or [8, Appendix A]. For $u,s>0$ and $w\in(-u,s)$, we have the chain of continuous inclusions

$$\mathcal{S}(\mathbb{R}^d) \hookrightarrow \mathsf{H}_r^s(\mathbb{R}^d) \hookrightarrow \mathsf{H}_r^w(\mathbb{R}^d) \hookrightarrow \mathsf{H}_r^{-u}(\mathbb{R}^d) \hookrightarrow \mathcal{S}'(\mathbb{R}^d). \tag{3.1}$$

In this setting, $\mathsf{H}_r^w(\mathbb{R}^d)$ is called an *intermediate space* between $\mathsf{H}_r^s(\mathbb{R}^d)$ and $\mathsf{H}_r^{-u}(\mathbb{R}^d)$. Recall also that $\mathsf{H}_r^0(\mathbb{R}^d)$ equals $\mathbb{L}^r(\mathbb{R}^d,\mathrm{d}x)$ with equivalent norms; so for $w\ge0$, elements of $\mathsf{H}_r^w(\mathbb{R}^d)$ are functions, not merely distributions. In addition, for $r'$ conjugate to $r$ via $1 = 1/r + 1/r'$, we have that $\mathsf{H}_{r'}^{-w}(\mathbb{R}^d)$ is the dual space of $\mathsf{H}_r^w(\mathbb{R}^d)$.

Let us revisit the main result from [8], the first paper in this series, with notation adapted to the present setting.

**Theorem 3.1** | *Impose Assumptions 1.1, 1.2 and 1.3. Then there exist real numbers $w > 0$ and $r > 1$ such that for any cluster point $\mathbf{P}_{x_0}^{\infty,k}$ of the sequence $(\mathbf{P}_{x_0}^{n,k})_{n \in \mathbb{N}}$, we have*

$$\mathbf{P}_{x_0}^{\infty,k}\left[\left\{\omega \in \Omega \,:\, \frac{\mathrm{d}\Lambda_t(\omega)}{\mathrm{d}x} = \bar{p}(t, \omega) \text{ is in } \mathsf{H}_r^w(\mathbb{R}^d) \text{ for almost all } t \in (0, T]\right\}\right] = 1, \quad (3.2)$$

*where $\bar{p} : (0, T] \times \Omega \to \mathsf{H}_r^w(\mathbb{R}^d)$ is a strongly measurable function.*

*More precisely, for each $t \in (0, T]$, we have $\mathbf{P}_{x_0}^{\infty,k}$-a.s. that $\Lambda_t \ll \mathrm{d}x$, and for any $q \geq 1$, $w$ and $r$ can be chosen in such a way that there exists $1 > \gamma > 0$ such that the function $t \mapsto (\omega \mapsto \mathrm{d}\Lambda_t(\omega)/\mathrm{d}x)$ defines a strongly measurable map*

$$p^k : (0, T] \to \mathbb{L}^q\left((\Omega, \mathcal{F}, \mathbf{P}_{x_0}^{\infty,k}); \mathsf{H}_r^w(\mathbb{R}^d)\right)$$

*which satisfies for some constant $c_{\text{Thm. 3.1}} < \infty$ the bound*

$$\left\|\|p^k(t)\|_{\mathsf{H}_r^w(\mathbb{R}^d)}\right\|_{\mathbb{L}^q(\mathbf{P}_{x_0}^{\infty,k})} \leq c_{\text{Thm. 3.1}}(1 \wedge t)^{-\gamma} \quad (3.3)$$

*for all $t \in (0, T]$. In (3.3), $c_{\text{Thm. 3.1}}$ depends on $r$, $w$, $q$, $\gamma$, but not on $t$ and neither on $\mathbf{P}_{x_0}^{\infty,k}$.*

*Finally, $\bar{p}$ is unique up to $(\mathrm{d}t \otimes \mathrm{d}\mathbf{P}_{x_0}^{\infty,k})$-a.e. equality and we have that $\bar{p}(t, \cdot) = p^k(t)$ in $\mathbb{L}^q((\Omega, \mathcal{F}, \mathbf{P}_{x_0}^{\infty,k}); \mathsf{H}_r^w(\mathbb{R}^d, \mathrm{d}x))$, for almost every $t \in (0, T]$. In particular, $\bar{p}(t, \omega)$ is the density of $\Lambda_t(\omega)$ $(\mathrm{d}t \otimes \mathrm{d}\mathbf{P}_{x_0}^{\infty,k})$-a.e., and $\bar{p}(t, \cdot)$ satisfies the bound (3.3) for almost every $t \in (0, T]$.*

***Proof*** This follows from Lemma 1.4 and [8, Thm. 1.5]. Indeed, Lemma 1.4 shows that the coefficients in the $\delta_k$-smoothed $n$-particle system in (1.1) and (1.2) satisfy Assumptions 1.2 and 1.3. These are precisely the properties required by [8, Thm. 1.5]. Since by its definition in (1.17), $(\mathbf{P}_{x_0}^{n,k})_{n \in \mathbb{N}}$ is the sequence of empirical measure flow laws, and $\mathbf{P}_{x_0}^{\infty,k}$ is a cluster point, we obtain from that result the emergence of regularity property (3.2) and the regularity estimate (3.3) under $\mathbf{P}_{x_0}^{\infty,k}$. □

We also take note of the following remark which is vital for our developments below.

**Remark 3.2** | The parameters $r$, $w$, $q$, and $\gamma$ in Theorem 3.1 can be taken to be independent of $k$. This is because by Lemma 1.4 we have that $\|b^\delta\|_\infty \leq \|b\|_\infty$ for any $\delta > 0$, and because the estimates used in [8, Proof of Thm 1.5, Step 1] depend next to the dimension $d$ and the time horizon $T$ only on the primitive quantities from Assumptions 1.2 and 1.3(H), namely the bounds for $\|b^\delta\|_\infty$, $\|\bar{\sigma}\|_\infty$, $\|\sigma\|_\infty$, the $\beta$-Hölder-seminorm bounds for $\sup_{t \in [0,T]}[\sigma_t]_\beta$ and $[\bar{\sigma}]_\beta$, and finally the ellipticity bound $\kappa$ from Assumption 1.3(E). Since all these can be controlled uniformly in $k$, so can the parameters $r$, $w$, $q$, $\gamma$. The same is true for the constant $c_{\text{Thm. 3.1}}$ appearing in Theorem 3.1. It depends next to the primitives from Assumptions 1.2 and 1.3(H) only on the values of $r$, $w$, $q$ and $\gamma$; see [8, Rem. 1.6].

## 3.2 A uniform regularity estimate

From the emergence of regularity result in Theorem 3.1, we obtain in (3.3) for each $k \in \mathbb{N}$ a regularity estimate for the density of $\Lambda$ under $\mathbf{P}_{x_0}^{\infty,k}$. We now verify that the regularity estimates we obtain in this way are not only uniform in $k \in \mathbb{N}$ but, in fact, carry over to narrow cluster points $\mathbf{P}_{x_0}^{\infty,\infty}$ of $(\mathbf{P}_{x_0}^{\infty,k})_{k\in\mathbb{N}}$.

**Proposition 3.3** | *Under Assumptions 1.1, 1.2 and 1.3, let $\mathbf{P}_{x_0}^{\infty,\infty}$ be a cluster point of $(\mathbf{P}_{x_0}^{\infty,k})_{k\in\mathbb{N}}$. Then for all $k \in \mathbb{N}\cup\{\infty\}$ and each $t \in (0,T]$, we have $\mathbf{P}_{x_0}^{\infty,k}$-a.s. that $\Lambda_t \ll \mathrm{d}x$. In fact, there exist exponents*

$$1 < q < \infty \text{ and } 1 < r < \infty, \ 0 < w < \infty \text{ as well as } 0 < \gamma < 1\,, \tag{3.4}$$

*a constant $c = c_{\mathrm{Thm.\,3.1}} < \infty$, all independent of $(\mathbf{P}_{x_0}^{\infty,k})_{k\in\mathbb{N}}$ and $\mathbf{P}_{x_0}^{\infty,\infty}$, and for each $k \in \mathbb{N}\cup\{\infty\}$ strongly a measurable function*

$$p^k : (0,T] \to \mathbb{L}^q\Big((\Omega, \mathcal{F}, \mathbf{P}_{x_0}^{\infty,k}); \mathsf{H}_r^w(\mathbb{R}^d)\Big)\,, \tag{3.5}$$

*which jointly satisfy for all $k \in \mathbb{N}\cup\{\infty\}$ the bound*

$$\Big\| \|p^k(t)\|_{\mathsf{H}_r^w(\mathbb{R}^d)} \Big\|_{\mathbb{L}^q(\mathbf{P}_{x_0}^{\infty,k})} \leq c(1 \wedge t)^{-\gamma} \quad \text{for all } t \in (0,T] \tag{3.6}$$

*and for which for each $t \in (0,T]$, we have $\mathbf{P}_{x_0}^{\infty,k}$-a.s. that*

$$\mathrm{d}\Lambda_t(\omega)/\mathrm{d}x = p^k(t)(\omega) \quad \text{in } \mathsf{H}_r^w(\mathbb{R}^d)\,. \tag{3.7}$$

*As a consequence, $\mathrm{d}\Lambda_t(\omega)/\mathrm{d}x = p^k(t)(\omega)$ in $\mathsf{H}_r^w(\mathbb{R}^d)$ $(\mathrm{d}t\otimes\mathrm{d}\mathbf{P}_{x_0}^{\infty,k})$-a.e. for all $k \in \mathbb{N}\cup\{\infty\}$.*

**Proof** We first consider the case $k \in \mathbb{N}$ and then argue that the result for $k = \infty$ follows from it.

**Step 1** First, for each $k \in \mathbb{N}$, the measure $\mathbf{P}_{x_0}^{\infty,k}$ chosen in (1.18) is a cluster point of $(\mathbf{P}_{x_0}^{n,k})_{n\in\mathbb{N}}$. By Theorem 3.1, we can find exponents $q$, $r$, $w$, $\gamma$ as in (3.4), a constant $c_{\mathrm{Thm.\,3.1}} < \infty$ and for each $k \in \mathbb{N}$ a strongly measurable function $p^k$ as in (3.5) satisfying the regularity estimate (3.6) and for each $t \in (0,T]$ the equality in (3.7). Note that by Remark 3.2, the values of $q$, $r$, $w$, $\gamma$ and $c_{\mathrm{Thm.\,3.1}}$ can be taken to be independent of $k \in \mathbb{N}$.

**Step 2** For $k = \infty$, we proceed as follows. The representation

$$\|\Lambda_t(\omega)\|_{\mathsf{H}_r^w(\mathbb{R}^d)} = \sup\{|\Lambda_t(\omega)[\phi]| \ : \ \phi \in \mathcal{S} \text{ with } \|\phi\|_{\mathsf{H}_{r'}^{-w}(\mathbb{R}^d)} = 1\}$$

gives for $a \in \mathbb{R}$ that

$$\begin{aligned}
&\{\omega \in \Omega : \|\Lambda_t(\omega)\|_{\mathsf{H}_r^w(\mathbb{R}^d)} > a\} \\
&= \Big\{\omega \in \Omega : \sup\{|\Lambda_t(\omega)[\phi]| \ : \ \phi \in \mathcal{S} \text{ with } \|\phi\|_{\mathsf{H}_{r'}^{-w}(\mathbb{R}^d)} = 1\} > a\Big\} \\
&= \{\omega \in \Omega : \text{there exists } \phi \in \mathcal{S} \text{ with } \|\phi\|_{\mathsf{H}_{r'}^{-w}(\mathbb{R}^d)} = 1 \text{ and } |\Lambda_t(\omega)[\phi]| > a\}\,.
\end{aligned}$$

Since the last set is open in $\Omega$, the function $\omega \mapsto \|\Lambda_t(\omega)\|_{\mathsf{H}_r^w(\mathbb{R}^d)}$ is lower-semicontinuous. Since also $\|\Lambda(\omega)\|_{\mathsf{H}_r^w(\mathbb{R}^d)} \geq 0$ for all $\omega \in \Omega$, we deduce from the Portmanteau theorem that

$$\begin{aligned}
\Big\| \|\Lambda_t\|_{\mathsf{H}_r^w(\mathbb{R}^d)} \Big\|_{\mathbb{L}^q(\mathbf{P}_{x_0}^{\infty,\infty})} &\leq \liminf_{k\to\infty} \Big\| \|\Lambda_t\|_{\mathsf{H}_r^w(\mathbb{R}^d)} \Big\|_{\mathbb{L}^q(\mathbf{P}_{x_0}^{\infty,k})} \\
&\leq \liminf_{k\to\infty} \Big\| \|p^k(t)\|_{\mathsf{H}_r^w(\mathbb{R}^d)} \Big\|_{\mathbb{L}^q(\mathbf{P}_{x_0}^{\infty,k})} \\
&\leq c_{\mathrm{Thm.\,3.1}}(1 \wedge t)^{-\gamma}\,.
\end{aligned}$$

This last estimate corresponds to the estimate in [8, Eqn. (6.14)], so proceeding as in [8, Proof of Thm. 1.5, Steps 3–5], we deduce (3.6) and (3.7) also for the case $k = \infty$. $\qquad\square$

# 4   A REFINED TIGHTNESS RESULT

In Lemma 2.9, we obtained a basic tightness result for the sequence $(\mathbf{P}_{x_0}^{\infty,k})_{k\in\mathbb{N}}$ on the space $\Omega$. This result, however, is not sufficient for our purposes below. We therefore establish a second, more refined compactness result. Define for $1 \leq p, \xi \leq \infty$ the space

$$\Omega_{p,\xi} := \mathbb{L}^p\Big([0,T]; \mathbb{L}^\xi(\mathrm{d}x)\Big) \tag{4.1}$$

Note that the topologies on $\Omega$ and on $\Omega_{p,\xi}$ are not related.

If $0 < \gamma < 1$, $r > 1$ and $s = 0$ are as in Theorem 3.1, and $1 \leq p \leq \infty$ is such that $p\gamma < 1$, then we have for each $k \in \mathbb{N}$ that

$$\mathbf{P}_{x_0}^{\infty,k}[\Lambda \in \Omega_{p,r}] = 1\,.$$

Indeed, this directly follows from the estimate in (3.3) by which we have that

$$
\begin{aligned}
\mathbf{E}^{\mathbf{P}_{x_0}^{\infty,k}}\left[\int_0^T \|\Lambda_t\|_{\mathbb{L}^r(\mathrm{d}x)}^p \,\mathrm{d}t\right] &= \mathbf{E}^{\mathbf{P}_{x_0}^{\infty,k}}\left[\int_0^T \|p^k(t)\|_{\mathbb{L}^r(\mathrm{d}x)}^p \,\mathrm{d}t\right] \\
&= \int_0^T \mathbf{E}^{\mathbf{P}_{x_0}^{\infty,k}}[\|p^k(t)\|_{\mathbb{L}^r(\mathrm{d}x)}^p]\,\mathrm{d}t \\
&\leq c_{Thm.\,3.1} \int_0^T (1 \wedge t)^{-\gamma p}\,\mathrm{d}t < \infty\,,
\end{aligned}
$$

which in turn implies for each $k \in \mathbb{N}$ that $\mathbf{P}_{x_0}^{\infty,k}$-a.s. $\|\|\Lambda_t\|_{\mathbb{L}^r(\mathrm{d}x)}\|_{\mathbb{L}^q(\mathrm{d}t)} < \infty$.

In this section, we aim to establish for appropriate choices of $p$ and $r$ tightness in $\Omega_{p,r}$. In other words, for each $\varepsilon > 0$ we aim to exhibit a compact set $K \subseteq \Omega_{p,r}$ with $\sup_{k\in\mathbb{N}} \mathbf{P}_{x_0}^{\infty,k}[\Lambda \in K] \geq 1 - \varepsilon$. This is considerably more intricate than showing tightness on $\Omega$.

## 4.1   Time-regularity on compacts in $(0,T]$

For our developments below, the following result is vital. It can be seen as a refinement of (3.2) in Theorem 3.1. To state the result, let $\mathbf{C}((0,T];\mathbb{L}^r(\mathrm{d}x))$ be the Fréchet space obtained by considering the topology of uniform convergence on compacts in $(0,T]$.

**Proposition 4.1** | *Impose Assumptions 1.1, 1.2 and 1.3. For each $k \in \mathbb{N}$ there exists a modification $\Lambda^k = (\Lambda_t^k)_{t\in[0,T]}$ of $\Lambda$ such that*

$$\mathbf{P}_{x_0}^{\infty,k}\Big[(\Lambda_t^k)_{t\in(0,T]} \in \mathbf{C}\Big((0,T];\mathbb{L}^r(\mathrm{d}x)\Big)\Big] = 1\,. \tag{4.2}$$

*In fact, there exists $\eta_0 > 0$ such that for any $\delta > 0$, any $q \in [1,\infty)$ and any $\eta \in [0,\eta_0)$ the modification $\Lambda^k$ can be chosen such that $\mathbf{P}_{x_0}^{\infty,k}$-a.s. we have $(\Lambda^k)_{t\in[\delta,T]} \in \mathsf{H\ddot{o}l}_\eta([\delta,T];\mathbb{L}^r(\mathrm{d}x))$. More precisely, there exists a random variable*

$$C_{\eta,q,\delta} \in \mathbb{L}^q(\Omega, \mathcal{F}, \mathbf{P}_{x_0}^{\infty,k}) \tag{4.3}$$

*satisfying* $\mathbf{P}_{x_0}^{\infty,k}$*-a.s. that*

$$\|\Lambda_t^k - \Lambda_s^k\|_{\mathbb{L}^r(\mathrm{d}x)} \leq C_{\eta,q,\delta}|t - s|^\eta \qquad \text{for all } \delta \leq s \leq t \leq T.$$ (4.4)

*Finally, all these properties also hold with* $k = \infty$ *under any cluster point* $\mathbf{P}_{x_0}^{\infty,\infty}$ *of* $(\mathbf{P}_{x_0}^{\infty,k})_{k\in\mathbb{N}}$.

We will not distinguish between $\Lambda$ and its modifications in the notation. To prepare for the proof, we begin with two auxiliary results.

**Lemma 4.2** | *Under Assumptions 1.1, 1.2 and 1.3, let $\boldsymbol{X}^{n,\delta_k}$ be the weak solution of* (1.1), (1.2). *For any $q \geq 1$ and $n \in \mathbb{N}$ we have that*

$$\|X_t^{1,n} - X_s^{1,n}\|_{\mathbb{L}^q(\mathbb{P}_{x_0}^{n,k})} \leq c|t - s|^{1/2} \qquad \text{for all } 0 \leq s \leq t \leq T,$$ (4.5)

*where the constant $c = c(b, \sigma, \bar{\sigma}, p, T)$ depends neither on $n \in \mathbb{N}$ nor on $k \in \mathbb{N}$.*

**Proof**   The argument is standard and uses the boundedness of $b_k^\delta$, $\sigma$ and $\bar{\sigma}$ and the Burkholder–David–Gundy inequality. □

**Lemma 4.3** | *Under Assumptions 1.1, 1.2 and 1.3, there exists $\eta_0 > 0$ such that for any $\delta > 0$ and $q \geq 1$ we can find a constant $c_{\delta,\eta_0,q} < \infty$ satisfying*

$$\left\| \|\Lambda_t - \Lambda_s\|_{\mathbb{L}^r(\mathrm{d}x)} \right\|_{\mathbb{L}^q(\mathbf{P}_{x_0}^{\infty,k})} \leq c_{\delta,\eta_0,q}|t - s|^{\eta_0} \qquad \text{for all } \delta \leq s \leq t \leq T,$$ (4.6)

*for any cluster point $\mathbf{P}_{x_0}^{\infty,k}$ of $(\mathbf{P}_{x_0}^{n,k})_{n\in\mathbb{N}}$. The constant $c_{\delta,\eta_0,q}$ does not depend on $k \in \mathbb{N}$. Finally, (4.6) remains valid with $k = \infty$ for any cluster point $\mathbf{P}_{x_0}^{\infty,\infty}$ of $(\mathbf{P}_{x_0}^{\infty,k})_{k\in\mathbb{N}}$.*

**Proof**   Fix $\delta > 0$. Let $w > 0$, $r > 1$ and $\gamma < 1$ be as in Theorem 3.1. Let $u = 1 + d/r'$ and choose $\theta \in (0,1)$ such that $0 \leq \alpha := \theta w - (1 - \theta)u$. From Bergh and Löfström [5, Thm. 4.1.2 and 6.4.5(7)], we have that

$$\|\Lambda_t - \Lambda_s\|_{\mathsf{H}_r^\alpha} \leq \|\Lambda_t - \Lambda_s\|_{\mathsf{H}_r^w}^\theta \|\Lambda_t - \Lambda_s\|_{\mathsf{H}_r^{-u}}^{1-\theta}.$$

For $q > 1$, we now apply Hölder's inequality with conjugate exponents $1/q = 1/q_1 + 1/q_2$ and then Jensen's inequality to show that

$$\left\| \|\Lambda_t - \Lambda_s\|_{\mathsf{H}_r^\alpha} \right\|_{\mathbb{L}^q(\mathbf{P}_{x_0}^{\infty,k})} \leq \left\| \|\Lambda_t - \Lambda_s\|_{\mathsf{H}_r^w}^\theta \right\|_{\mathbb{L}^{q_1}(\mathbf{P}_{x_0}^{\infty,k})} \left\| \|\Lambda_t - \Lambda_s\|_{\mathsf{H}_r^{-u}}^{1-\theta} \right\|_{\mathbb{L}^{q_2}(\mathbf{P}_{x_0}^{\infty,k})}$$
$$\leq \left\| \|\Lambda_t - \Lambda_s\|_{\mathsf{H}_r^w} \right\|_{\mathbb{L}^{q_1}(\mathbf{P}_{x_0}^{\infty,k})}^\theta \left\| \|\Lambda_t - \Lambda_s\|_{\mathsf{H}_r^{-u}} \right\|_{\mathbb{L}^{q_2}(\mathbf{P}_{x_0}^{\infty,k})}^{1-\theta}.$$ (4.7)

To bound the first factor on the right-hand side of (4.7), we use (3.3) in Theorem 3.1. By that result and the subadditivity of $x \mapsto x^\theta$ we get for $\delta \leq s \leq t \leq T$ that

$$\left\| \|\Lambda_t - \Lambda_s\|_{\mathsf{H}_r^w} \right\|_{\mathbb{L}^{q_1}(\mathbf{P}_{x_0}^{\infty,k})}^\theta \leq \left\| \|\Lambda_t\|_{\mathsf{H}_r^w} \right\|_{\mathbb{L}^{q_1}(\mathbf{P}_{x_0}^{\infty,k})}^\theta + \left\| \|\Lambda_s\|_{\mathsf{H}_r^w} \right\|_{\mathbb{L}^{q_1}(\mathbf{P}_{x_0}^{\infty,k})}^\theta$$
$$\leq 2c_{Thm.\,3.1}^\theta (1 \wedge \delta)^{-\gamma\theta}.$$ (4.8)

To control the second factor on the right-hand side of (4.7), we follow similar steps as in [8, Sec. 2.2]. First note that for any $0 < \rho < 1$ we have that $\mathsf{H}_r^u \hookrightarrow \mathsf{H\ddot{o}l}_\rho(\mathbb{R}^d)$. The duality of $\mathsf{H}_r^{-u}$ and $\mathsf{H}_{r'}^u$, the definition of the empirical measure $\mu^n$ in (1.2) and the continuous embedding give

$$
\begin{aligned}
\|\mu_t^n - \mu_s^n\|_{\mathsf{H}_r^{-u}} &= \sup_{\|\phi\|_{\mathsf{H}_{r'}^u} \le 1} |\mu_t^n[\phi] - \mu_s^n[\phi]| \\
&= \sup_{\|\phi\|_{\mathsf{H}_{r'}^u} \le 1} \left| \frac{1}{n} \sum_{i=1}^n \phi(X_t^{i,n}) - \phi(X_s^{i,n}) \right| \\
&\le c_\rho \frac{1}{n} \sum_{i=1}^n |X_t^{i,n} - X_s^{i,n}|^\rho \,.
\end{aligned}
$$

The definition (1.17) of the measure $\mathbf{P}_{x_0}^{n,k}$, the exchangeability of $X^n$ under $\mathbb{P}_{x_0}^{n,\delta_k}$ from Lemma 1.5, and the bound (4.5) in Lemma 4.2 combine to give the estimate

$$
\begin{aligned}
\left\| \|\Lambda_t - \Lambda_s\|_{\mathsf{H}_r^{-u}}^{1-\theta} \right\|_{\mathbb{L}^{p_2}(\mathbf{P}_{x_0}^{n,k})} &\le \left\| \|\mu_t^n - \mu_s^n\|_{\mathsf{H}_r^{-u}} \right\|_{\mathbb{L}^{q_2}(\mathbb{P}_{x_0}^{n,k})}^{1-\theta} \\
&\le \left\| |X_t^{i,n} - X_s^{i,n}|^\rho \right\|_{\mathbb{L}^{q_2}(\mathbb{P}_{x_0}^{n,k})}^{1-\theta} \\
&\le c_{b,\sigma,\tilde{\sigma},q_2,\rho}^{1-\theta} |t-s|^{\rho(1-\theta)/2} \,.
\end{aligned}
\tag{4.9}
$$

We now set $\eta_0 := \rho(1-\theta)/2$ and $c_{\delta,\eta_0,q} := 2c_{Thm.\ 3.1}^\theta (1 \wedge \delta)^{-\gamma\theta} c_{b,\sigma,\tilde{\sigma},q_2,\rho}^{1-\theta}$ and combine (4.7), (4.8) and (4.10) to arrive at

$$
\left\| \|\Lambda_t - \Lambda_s\|_{\mathsf{H}_r^\rho} \right\|_{\mathbb{L}^q(\mathbf{P}_{x_0}^{n,k})} \le c_{\delta,\eta_0,q} |t-s|^{\eta_0} \qquad \text{for all } \delta \le s \le t \le T \,.
$$

Just as in Step 2 of the proof of Proposition 3.3, we obtain for any cluster point $\mathbf{P}_{x_0}^{\infty,k}$ of $(\mathbf{P}_{x_0}^{n,k})_{n \in \mathbb{N}}$ for all $\delta \le s \le t \le T$ that

$$
\left\| \|\Lambda_t - \Lambda_s\|_{\mathsf{H}_r^\rho} \right\|_{\mathbb{L}^q(\mathbf{P}_{x_0}^{\infty,k})} \le \liminf_{n \to \infty} \left\| \|\Lambda_t - \Lambda_s\|_{\mathsf{H}_r^\rho} \right\|_{\mathbb{L}^q(\mathbf{P}_{x_0}^{n,k})} \le c_{\delta,\eta_0,q} |t-s|^{\eta_0} \,.
\tag{4.10}
$$

Note that the choice of $q$ is independent of $\eta_0$ and $\delta$. Finally, since $\mathsf{H}_r^0 = \mathbb{L}^r(\mathrm{d}x)$ with equivalent norms, (4.6) follows. Finally, (4.10) remains valid when we take $\liminf_{k\to\infty}$, so that we get

$$
\left\| \|\Lambda_t - \Lambda_s\|_{\mathsf{H}_r^\rho} \right\|_{\mathbb{L}^q(\mathbf{P}_{x_0}^{\infty,\infty})} \le c_{\delta,\eta_0,q} |t-s|^{\eta_0}
$$

for any cluster point $\mathbf{P}_{x_0}^{\infty,\infty}$ of $(\mathbf{P}_{x_0}^{\infty,k})_{k \in \mathbb{N}}$. This concludes the proof. $\qquad \square$

**Proof of Proposition 4.1**  Let $\delta > 0$. By Lemma 4.3, we find $q_0 > 1/\eta_0 > 1$ and $c_{\delta,\eta_0,q_0}$ such that (4.6) holds for all $\delta \le s \le t \le T$. Now apply the Kolmogorov continuity criterion from e.g. Friz and Victoir [15, Thm. A.10] to deduce that $(\Lambda_t)_{t \in [\delta,T]}$ possesses a modification $(\Lambda_t^k)_{t \in [\delta,T]}$ that is Hölder-continuous of order $\eta$ with $\eta \in [0, \eta_0 - 1/q_0)$. In addition, the random Hölder-constant $C_{\eta,q_0,\delta}$ satisfies $\|C_{\eta,q_0,\delta}\|_{\mathbb{L}^{q_0}(\mathbf{P}_{x_0}^{\infty,k})} \le c_{\delta,\eta_0,q_0}$. Since $\mathbb{L}^{q_0}(\mathbf{P}_{x_0}^{\infty,k}) \hookrightarrow \mathbb{L}^q(\mathbf{P}_{x_0}^{\infty,k})$ if $q \le p_0$, we get (4.3) and (4.4). By pasting together these continuous modifications for arbitrary compact subsets of $(0,T]$, we deduce (4.2) and the proof is complete. $\qquad \square$

## 4.2 Tightness

We now come to our second tightness result.

**Lemma 4.4** | *Under Assumptions 1.1, 1.2 and 1.3 we can find for any $\varepsilon > 0$ a compact set $K \subseteq \Omega_{p,1}$ such that*

$$\inf_{k \in \mathbb{N}} \mathbf{P}_{x_0}^{\infty,k}[\Lambda \in K] > 1 - \varepsilon. \qquad (4.11)$$

*In other words, for any $p \geq 1$ the sequence $(\mathbf{P}_{x_0}^{\infty,k})_{k \in \mathbb{N}}$ is tight on $\Omega_{p,1}$, and any cluster point $\mathbf{P}_{x_0}^{\infty,\infty}$ of $(\mathbf{P}_{x_0}^{\infty,k})_{k \in \mathbb{N}}$ relative to narrow convergence on $\Omega_{p,1}$ satisfies $\mathbf{P}_{x_0}^{\infty,\infty}[\Lambda \in \Omega_{p,1}] = 1$.*

**Proof** We start with some auxiliary estimates to apply the Kolmogorov–M. Riesz–Fréchet compactness ccriterion for $\Omega_{p,1} \coloneqq \mathbb{L}^p([0,T]; \mathbb{L}^1(\mathrm{d}x))$; see Brézis [6, Th. 4.26]. Let $\varepsilon > 0$.

Theorem 3.1 shows that for any $k \in \mathbb{N}$ we have $\mathrm{d}t \otimes \mathrm{d}\mathbf{P}_{x_0}^{\infty,k}$-a.e. that $\|\Lambda_t\|_{\mathbb{L}^1(\mathrm{d}x)} = 1$. Therefore Markov's inequality shows that for some $c > 0$ we have

$$\sup_{k \in \mathbb{N}} \mathbf{P}_{x_0}^{\infty,k}\left[\int_0^T \|\Lambda_t\|_{\mathbb{L}^1(\mathrm{d}x)}^q \, \mathrm{d}t \geq c\right] \leq \frac{1}{c} \sup_{k \in \mathbb{N}} \mathbf{E}^{\mathbf{P}_{x_0}^{\infty,k}}\left[\int_0^T \|\Lambda_t\|_{\mathbb{L}^1(\mathrm{d}x)}^q \, \mathrm{d}t\right] \leq \frac{T}{c} \leq \varepsilon. \qquad (4.12)$$

In fact, the left-hand side in the above display is identically zero whenever $c > T$.

We next show that for any $c > 0$ we can find $0 < \delta < 1$ such that

$$\sup_{k \in \mathbb{N}} \mathbf{P}_{x_0}^{\infty,k}\left[\sup_{|h| < \delta} \int_0^T \|\Lambda_{t+h} - \Lambda_t\|_{\mathbb{L}^1(\mathrm{d}x)}^q \, \mathrm{d}t \geq c\right] \leq \varepsilon. \qquad (4.13)$$

This is more intricate. Fix $c > 0$ and use Markov's inequality to see that for any intermediate $t_0 \in (0, T)$, to be chosen later, we have that

$$\begin{aligned}
\mathbf{P}_{x_0}^{\infty,k}\left[\sup_{|h| < \delta} \int_0^T \|\Lambda_{t+h} - \Lambda_t\|_{\mathbb{L}^1(\mathrm{d}x)}^q \, \mathrm{d}t \geq c\right] &\leq \frac{1}{c} \mathbf{E}^{\mathbf{P}_{x_0}^{\infty,k}}\left[\sup_{|h| < \delta} \int_0^T \|\Lambda_{t+h} - \Lambda_t\|_{\mathbb{L}^1(\mathrm{d}x)}^q \, \mathrm{d}t\right] \\
&\leq \frac{1}{c} \mathbf{E}^{\mathbf{P}_{x_0}^{\infty,k}}\left[\sup_{|h| < \delta} \int_0^{t_0} \|\Lambda_{t+h} - \Lambda_t\|_{\mathbb{L}^1(\mathrm{d}x)}^q \, \mathrm{d}t\right] \\
&\quad + \frac{1}{c} \mathbf{E}^{\mathbf{P}_{x_0}^{\infty,k}}\left[\sup_{|h| < \delta} \int_{t_0}^T \|\Lambda_{t+h} - \Lambda_t\|_{\mathbb{L}^1(\mathrm{d}x)}^q \, \mathrm{d}t\right] \\
&=: \frac{1}{c} A_{t_0,\delta} + \frac{1}{c} B_{t_0,\delta}. \qquad (4.14)
\end{aligned}$$

Note first that $\|\Lambda_0\|_{\mathbb{L}^1(\mathrm{d}x)} = 1$. Note next that by Theorem 3.1 and Proposition 4.1, more specifically, by combining (3.2) and (4.2), we have $\mathbf{P}_{x_0}^{\infty,k}$-a.s. that $\|\Lambda_t\|_{\mathbb{L}^1(\mathrm{d}x)} = 1$ simultaneously for all $t \in (0, T]$. Therefore,

$$A_{t_0,\delta} = \mathbf{E}^{\mathbf{P}_{x_0}^{\infty,k}}\left[\sup_{|h| < \delta} \int_0^{t_0} \|\Lambda_{t+h} - \Lambda_t\|_{\mathbb{L}^1(\mathrm{d}x)}^q \, \mathrm{d}t\right] \leq 2^q t_0. \qquad (4.15)$$

This bound is independent of $\delta$. We can thus choose $t_0 > 0$ such that $2^p t_0 / c < \varepsilon/2$. For $B_{t_0,\delta}$, let us first fix $0 < \delta_0 < t_0$. By Proposition 4.1, $\Lambda$ can be taken to be Hölder-continuous of some order $\eta > 0$ on $[t_0 - \delta_0, T]$. For any $0 < \delta \leq \delta_0$ we thus get that

$$
\begin{aligned}
B_{t_0, \delta} &= \mathbf{E}^{\mathbf{P}_{x_0}^{\infty, k}} \left[ \sup_{|h| < \delta} \int_{t_0}^T \| \Lambda_{t+h} - \Lambda_t \|_{\mathbb{L}^1(\mathrm{d}x)}^q \, \mathrm{d}t \right] \\
&\leq \mathbf{E}^{\mathbf{P}_{x_0}^{\infty, k}} \left[ \sup_{|h| < \delta} \int_{t_0}^T C_{\eta, q, t_0 - \delta_0}^q |h|^{q\eta} \, \mathrm{d}t \right] \\
&\leq \| C_{\eta, q, t_0 - \delta_0} \|_{\mathbb{L}^q(\mathbf{P}_{x_0}^{\infty, k})}^q (T - t_0) \delta^{q\eta} \,.
\end{aligned} \tag{4.16}
$$

By (4.3), we can choose $\delta < \delta_0$ such that $\| C_{\eta, q, t_0 - \delta_0} \|_{\mathbb{L}^q(\mathbf{P}_{x_0}^{\infty, k})}^q (T - t_0) \delta^{q\eta_0} / \eta < \varepsilon/2$. We now combine the bounds in (4.15) and (4.16) via (4.14). Since by Proposition 4.1 all constants are independent of $k \in \mathbb{N}$, we arrive at (4.13).

To deduce (4.11), we can now follow standard arguments. Use (4.12) and (4.13) to choose $c > 0$ and $(\delta_\ell)_{\ell \in \mathbb{N}}$ with $\delta_\ell \to 0$ as $\ell \to \infty$ such that for

$$
\begin{aligned}
B_0 &:= \left\{ \omega \in \Omega_{p,1} \,:\, \int_0^T \| \Lambda_t \|_{\mathbb{L}^1(\mathrm{d}x)}^p \, \mathrm{d}t < c \right\} \\
B_\ell &:= \left\{ \omega \in \Omega_{p,1} \,:\, \sup_{|h| < \delta_\ell} \int_0^T \| \Lambda_{t+h} - \Lambda_t \|_{\mathbb{L}^1(\mathrm{d}x)}^p \, \mathrm{d}t < \frac{1}{\ell} \right\}
\end{aligned}
$$

we have for each $\ell \in \mathbb{N}_0$ that $\mathbf{P}_{x_0}^{\infty, k}[B_\ell] \geq 1 - \varepsilon / 2^{\ell+1}$ for all $k \in \mathbb{N}$. Let $K$ be the closure of $B := \bigcap_{\ell \in \mathbb{N}_0} B_\ell$. By the Kolmogorov–M. Riesz–Fréchet Theorem the set $K \subseteq \Omega_{p,1}$ is compact and by the choices above we have that $\inf_{k \in \mathbb{N}} \mathbf{P}_{x_0}^{\infty, k}[\Lambda \in K] > 1 - \varepsilon$. We thus have (4.11) and therefore tightness. □

**Proof of Lemma 1.8** By Lemma 2.9, $(\mathbf{P}_{x_0}^{\infty, k})_{k \in \mathbb{N}}$ is tight in $\Omega$. We can therefore find a subsequence $(\mathbf{P}_{x_0}^{\infty, k_\ell})_{\ell \in \mathbb{N}}$ that is narrowly convergent in $\Omega$. But by Lemma 4.4, the sequence $(\mathbf{P}_{x_0}^{\infty, k_\ell})_{\ell \in \mathbb{N}}$ is tight in $\Omega_{p,1}$. We therefore find a subsequence $(\mathbf{P}_{x_0}^{\infty, k_m})_{m \in \mathbb{N}}$ of $(\mathbf{P}_{x_0}^{\infty, k_\ell})_{\ell \in \mathbb{N}}$ that is in addition narrowly convergent in $\Omega_{p,1}$. □

### 4.3 Improved convergence result

Starting from the sequence $(\mathbf{P}_{x_0}^{\infty, k})_{k \in \mathbb{N}}$ fixed in (1.18), we obtain via Lemmas 2.9 and 4.4 a narrow cluster point $\mathbf{P}_{x_0}^{\infty, \infty}$ for both $\Omega := \mathbf{C}([0, T]; \mathcal{S}')$ and $\Omega_{p,1} := \mathbb{L}^p([0, T]; \mathbb{L}^1(\mathrm{d}x))$. In fact, for $p_2 > p_1$ with $p_1$ as in Lemma 2.9, we know that $\mathbf{P}_{x_0}^{\infty, m}$ for $m \in \mathbb{N} \cup \{\infty\}$ concentrate all their mass on the Polish space $\mathbf{C}([0, T]; \mathcal{H}_{-p_2})$. Hence there is a subsequence $(\mathbf{P}_{x_0}^{\infty, k_m})_{m \in \mathbb{N}}$ with $\mathbf{P}_{x_0}^{\infty, k_m} \to \mathbf{P}_{x_0}^{\infty, \infty}$ narrowly in $\mathbf{C}([0, T]; \mathcal{H}_{-p_2})$ and in $\Omega_{p,1}$ as $m \to \infty$. To simplify notation we relabel $\mathbf{P}_{x_0}^{\infty, m} := \mathbf{P}_{x_0}^{\infty, k_m}$ and $\delta_m := \delta_{k_m}$ for all $m \in \mathbb{N}$.

Consider next the diagonal embedding $\iota : \Omega \cap \Omega_{q,1} \to \Omega \times \Omega_{q,1}$, $\iota(u) = (u, u)$, and set $\mu_m := \iota_\# \mathbf{P}_{x_0}^{\infty, m}$ and $\mu := \iota_\# \mathbf{P}_{x_0}^{\infty, \infty}$. Since $\Omega$ and $\Omega_{q,1}$ are Polish spaces, the product $\Omega \times \Omega_{q,1}$ is Polish; by tightness and the marginal convergences we have $\mu_m \to \mu$ narrowly with all $\mu_m, \mu$ supported on the diagonal. Applying the Skorohod representation theorem on $\Omega \times \Omega_{q,1}$, see, e.g., Pollard [31, Thm. IV.13], we obtain a probability space $(\hat{\Omega}, \hat{\mathcal{F}}, \hat{\mathbf{P}})$ and random variables $\hat{\Lambda}^m$ for $m \in \mathbb{N} \cup \{\infty\}$, valued in both $\Omega$ and $\Omega_{q,1}$, such that

$$
\mathrm{Law}_{\hat{\mathbf{P}}}(\hat{\Lambda}^m) = \mathrm{Law}_{\mathbf{P}_{x_0}^{\infty, m}}(\Lambda) \qquad \text{for all } m \in \mathbb{N} \cup \{\infty\}. \tag{4.17}
$$

and $\hat{\mathbf{P}}$-a.s. we have

$$\lim_{m\to\infty} \hat{\Lambda}^m = \hat{\Lambda}^\infty \qquad \text{in } \mathbf{C}([0,T]; \mathscr{H}_{-p}), \tag{4.18}$$

$$\lim_{m\to\infty} \hat{\Lambda}^m = \hat{\Lambda}^\infty \qquad \text{in } \mathbb{L}^p\big([0,T]; \mathbb{L}^1(\mathrm{d}x)\big). \tag{4.19}$$

Note that by (1.9) and (4.17), $\hat{\Lambda}^m$ is $\hat{\mathbf{P}}$-a.s. $\mathbf{M}_1^+$-valued for all $m \in \mathbb{N} \cup \{\infty\}$.

With the help of Proposition 4.1, the mode of convergence in (4.19) can be improved.

**Proposition 4.5** | *Let* $1 \le \xi < r$. *Under Assumptions 1.1, 1.2 and 1.3, we have* $\hat{\mathbf{P}}$-*a.s. that*

$$\lim_{m\to\infty} (\hat{\Lambda}_t^m)_{t\in(0,T]} = (\hat{\Lambda}_t^\infty)_{t\in(0,T]} \qquad \text{in } \mathbf{C}\big((0,T]; \mathbb{L}^\xi(\mathrm{d}x)\big), \tag{4.20}$$

*where* $\mathbf{C}((0,T]; \mathbb{L}^\xi(\mathrm{d}x))$ *carries the topology of uniform convergence on compact sets.*

**Proof**    Let $1 \le \xi < r$. For some $\theta \in (0,1)$ we then have the interpolation bound

$$\|\hat{\Lambda}_t^m(\hat{\omega}) - \hat{\Lambda}_t^\infty(\hat{\omega})\|_{\mathbb{L}^\xi(\mathrm{d}x)} \le \|\hat{\Lambda}_t^m(\hat{\omega}) - \hat{\Lambda}_t^\infty(\hat{\omega})\|_{\mathbb{L}^1(\mathrm{d}x)}^\theta \|\hat{\Lambda}_t^m(\hat{\omega}) - \hat{\Lambda}_t^\infty(\hat{\omega})\|_{\mathbb{L}^r(\mathrm{d}x)}^{1-\theta}.$$

Now let $t_0 \in (0,T]$. We then get

$$\int_{t_0}^T \|\hat{\Lambda}_t^k(\hat{\omega}) - \hat{\Lambda}_t^\infty(\hat{\omega})\|_{\mathbb{L}^\xi(\mathrm{d}x)}^p \, \mathrm{d}t \le \int_{t_0}^T \|\hat{\Lambda}_t^k(\hat{\omega}) - \hat{\Lambda}_t^\infty(\hat{\omega})\|_{\mathbb{L}^1(\mathrm{d}x)}^{\theta p} \|\hat{\Lambda}_t^k(\hat{\omega}) - \hat{\Lambda}_t^\infty(\hat{\omega})\|_{\mathbb{L}^r(\mathrm{d}x)}^{(1-\theta)p} \, \mathrm{d}t$$

$$\le c_{r,p,\theta,t_0}(\hat{\omega}) \int_{t_0}^T \|\hat{\Lambda}_t^m(\hat{\omega}) - \hat{\Lambda}_t^\infty(\hat{\omega})\|_{\mathbb{L}^1(\mathrm{d}x)}^{p\theta} \, \mathrm{d}t, \tag{4.21}$$

where by (4.4) we have $\hat{\mathbf{P}}$-a.s. that

$$c_{r,p,\theta,t_0}(\hat{\omega}) := \sup_{t\in[t_0,T]} \|\hat{\Lambda}_t^m(\hat{\omega}) - \hat{\Lambda}_t^\infty(\hat{\omega})\|_{\mathbb{L}^r(\mathrm{d}x)}^{p(1-\theta)} < \infty. \tag{4.22}$$

Next let $\varepsilon > 0$ and define

$$A_k(\hat{\omega}) := \{t \in [t_0,T] : \|\hat{\Lambda}_t^m(\hat{\omega}) - \hat{\Lambda}_t^\infty(\hat{\omega})\|_{\mathbb{L}^\xi(\mathrm{d}x)} > \varepsilon\}.$$

By the definition of $A_k(\hat{\omega})$, the interpolation bound (4.21) and (4.22) we obtain

$$|A_k(\hat{\omega})| \le \frac{1}{\varepsilon^p} \int_{A_k(\hat{\omega})} \|\hat{\Lambda}_t^m(\hat{\omega}) - \hat{\Lambda}_t^\infty(\hat{\omega})\|_{\mathbb{L}^\xi(\mathrm{d}x)}^p \, \mathrm{d}t$$

$$\le \frac{1}{\varepsilon^p} \int_{t_0}^T \|\hat{\Lambda}_t^m(\hat{\omega}) - \hat{\Lambda}_t^\infty(\hat{\omega})\|_{\mathbb{L}^\xi(\mathrm{d}x)}^p \, \mathrm{d}t$$

$$\le \frac{c_{r,p,\theta,t_0}(\hat{\omega})}{\varepsilon^p} \int_{t_0}^T \|\hat{\Lambda}_t^m(\hat{\omega}) - \hat{\Lambda}_t^\infty(\hat{\omega})\|_{\mathbb{L}^1(\mathrm{d}x)}^{p\theta} \, \mathrm{d}t$$

$$\le \frac{c_{r,p,\theta,t_0}(\hat{\omega})}{\varepsilon^p} \Big\| \|\hat{\Lambda}_t^m(\hat{\omega}) - \hat{\Lambda}_t^\infty(\hat{\omega})\|_{\mathbb{L}^1(\mathrm{d}x)} \Big\|_{\mathbb{L}^p(\mathrm{d}t)}^{p\theta},$$

and so (4.18) gives that $\hat{\mathbf{P}}$-a.s. we have

$$|A_k(\hat{\omega})| \to 0 \qquad \text{as } k \to \infty. \tag{4.23}$$

Suppose now that for some $F \in \hat{\mathcal{F}}$ of positive $\hat{\mathbf{P}}$-measure we can find for each $\hat{\omega} \in F$ and $K \in \mathbb{N}$ an integer $k > K$ and a time $t_k \in [t_0, T]$, both of which may depend on $\hat{\omega}$, for with we have that

$$\|\hat{\Lambda}_{t_k}^m(\hat{\omega}) - \hat{\Lambda}_{t_k}^\infty(\hat{\omega})\|_{\mathbb{L}^\xi(\mathrm{d}x)} > 2\varepsilon \,.$$

By the Hölder-continuity established in Proposition 4.1, we find a $\hat{\mathbf{P}}$-a.s. finite constant $c(\hat{\omega})$ that is independent of $k \in \mathbb{N} \cup \{\infty\}$ for which we have

$$\|\hat{\Lambda}_t^m(\hat{\omega}) - \hat{\Lambda}_{t'}^m(\hat{\omega})\|_{\mathbb{L}^\xi(\mathrm{d}x)} \le c(\hat{\omega})|t - t'|^\alpha$$

for all $t_0 \le t \le t' \le T$. Therefore, there exists $\delta > 0$ which may depend on $\hat{\omega}$ but which is independent of $k \in \mathbb{N} \cup \{\infty\}$ such that for all $t' \in [\varepsilon, T]$ with $|t_k - t'| < \delta$ we get

$$\left\|\hat{\Lambda}_{t_k}^m(\hat{\omega}) - \hat{\Lambda}_{t_k}^\infty(\hat{\omega}) - \left(\hat{\Lambda}_{t'}^m(\hat{\omega}) - \hat{\Lambda}_{t'}^\infty(\hat{\omega})\right)\right\|_{\mathbb{L}^\xi(\mathrm{d}x)} < \varepsilon \,,$$

which in turn implies that

$$\left\|\hat{\Lambda}_{t'}^m(\hat{\omega}) - \hat{\Lambda}_{t'}^\infty(\hat{\omega})\right\|_{\mathbb{L}^\xi(\mathrm{d}x)} > 2\varepsilon - \varepsilon > \varepsilon \,.$$

By the definition of $A_k(\hat{\omega})$, we have $\{t' \in [t_0, T] : |t_k - t'| < \delta\} \subseteq A_k(\hat{\omega})$. Since $\delta > 0$ is independent of $k \in \mathbb{N}$, this implies that the set $A_k(\hat{\omega})$ must have Lebesgue measure bounded from below for infinitely many $k \in \mathbb{N}$. Since this contradicts (4.23), we conclude $\hat{\mathbf{P}}$-a.s. that $\sup_{t \in [t_0, T]} \|\hat{\Lambda}_t^m(\hat{\omega}) - \hat{\Lambda}_t^\infty(\hat{\omega})\|_{\mathbb{L}^\xi(\mathrm{d}x)} \to 0$ as $k \to \infty$, and since $t_0 \in (0, T]$ is arbitrary, this establishes (4.20) and concludes the proof. $\qquad \square$

## 5 FORMALLY LIMITING MARTINGALE PROBLEM

In this section we prove our main result which we repeat here for sake of reference.

**Theorem 1.9** | *Under Assumptions 1.1, 1.3 and 1.7, let $(\mathbf{P}_{x_0}^{\infty,k})_{k \in \mathbb{N}}$ be chosen as in (1.16), (1.18). For $p \ge 1$, let $\mathbf{P}_{x_0}^{\infty,\infty}$ be a successive cluster point as in Lemma 1.8. Then $\mathbf{P}_{x_0}^{\infty,\infty}$ solves the martingale problem $\mathrm{MP}(\nu, \mathscr{A}, \mathscr{Q})$. In particular, the coordinate process $\Lambda$ on $(\Omega, \mathcal{F}, \mathbf{P}_{x_0}^{\infty,\infty})$ possesses the $\mathcal{S}'$-valued semimartingale decomposition*

$$\Lambda_t = \Lambda_0 + \int_0^t \mathscr{A}_s(\Lambda_s)\,\mathrm{d}s + M_t, \qquad t \in [0, T]\,, \tag{1.19}$$

*where $\Lambda_0 \stackrel{d}{=} \nu$ under $\mathbf{P}_{x_0}^{\infty,\infty}$ and where $M = (M_t)_{t \in [0,T]}$ is in $\mathscr{M}_2^c(\mathbb{F}, \mathbf{P}_{x_0}^{\infty,\infty}; \mathcal{S}')$ with quadratic variation*

$$\langle\!\langle M \rangle\!\rangle_t = \int_0^t \mathscr{Q}_s(\Lambda_s)\,\mathrm{d}s, \qquad t \in [0, T]\,. \tag{1.20}$$

*Moreover, for all $1 \le q < \infty$, there exist exponents $1 < r < \infty$, $0 < w < \infty$ as well as $0 < \gamma < 1$, and a constant $c = c_{\mathrm{Thm.\ 3.1}} < \infty$, all independent of $\mathbf{P}_{x_0}^{\infty,\infty}$, such that the function $t \mapsto (\omega \mapsto \mathrm{d}\Lambda_t(\omega)/\mathrm{d}x)$ defines a strongly measurable map*

$$p^\infty : (0, T] \to \mathbb{L}^q\Big((\Omega, \mathcal{F}, \mathbf{P}_{x_0}^{\infty,\infty}); \mathsf{H}_r^w(\mathbb{R}^d)\Big)$$

*which satisfies for some constant the bound*

$$\left\| \|p^\infty(t)\|_{\mathsf{H}^w_r(\mathbb{R}^d)} \right\|_{\mathbb{L}^q(\mathbf{P}^{\infty,\infty}_{x_0})} \le c_{\text{Thm. 3.1}} (1 \wedge t)^{-\gamma}$$

*for all $t \in (0, T]$. In particular, for each $t \in (0, T]$, we have $\mathbf{P}^{\infty,\infty}_{x_0}$-a.s. that $\Lambda_t \ll \mathrm{d}x$ and $\mathrm{d}\Lambda_t(\omega)/\mathrm{d}x = p^\infty(t)(\omega)$ in $\mathsf{H}^w_r(\mathbb{R}^d)$ $(\mathrm{d}t \otimes \mathrm{d}\mathbf{P}^{\infty,\infty}_{x_0})$-a.e.*

The strategy of proof was already touched upon in Section 1.4. There we noted that in order to prove the $(\mathbb{F}, \mathbf{P}^{\infty,\infty}_{x_0})$-local martingale property of $M$ in (1.19), or equivalently (2.18), we should aim to establish (1.30). The following technical proposition gives a more general result from which (1.30) is obtained; see (5.2) below.

**Proposition 5.2** | *Under Assumptions 1.1, 1.3 and 1.7, let $(k_m)_{m \in \mathbb{N}}$ be a sequence such that $(\mathbf{P}^{\infty,k_m}_{x_0})_{m \in \mathbb{N}}$ converges to $\mathbf{P}^{\infty,\infty}_{x_0}$ in the narrow topology. Then*

$$\lim_{m \to \infty} \mathbf{E}^{\mathbf{P}^{\infty,k_m}_{x_0}} \left[ \int_s^t g_u(\Lambda) \mathscr{A}^{\delta_{k_m}}_u(\Lambda_u)[\phi] \, \mathrm{d}u \right] = \mathbf{E}^{\mathbf{P}^{\infty,\infty}_{x_0}} \left[ \int_s^t g_u(\Lambda) \mathscr{A}_u(\Lambda_u)[\phi] \, \mathrm{d}u \right] \quad (5.1)$$

*for any $\phi \in \mathcal{S}$, all $s, t \in [0, T]$ with $s \le t$, and all bounded and measurable functions $g : [0, T] \times \Omega \to \mathbb{R}$ such that for each $u \in [0, T]$, the real-valued function $g_u(\,\cdot\,) := g(u, \cdot\,)$ on $\Omega$ is continuous.*

The proof of this result will occupy us for the better half of this section.

**Almost sure representation**    Let $(\mathbf{P}^{\infty,m}_{x_0})_{m \in \mathbb{N}}$ and $\mathbf{P}^{\infty,\infty}_{x_0}$ be as in Section 4.3. Recall from Section 4.3 that $(\hat{\Lambda}_m)_{m \in \mathbb{N} \cup \{\infty\}}$ is a Skorohod representation. With this representation, we have for all $m \in \mathbb{N} \cup \{\infty\}$ and $\phi \in \mathcal{S}$ the identity

$$\mathbf{E}^{\mathbf{P}^{\infty,m}_{x_0}} \left[ \int_s^t g_u(\Lambda) \mathscr{A}^{\delta_m}_u(\Lambda_u)[\phi] \, \mathrm{d}u \right] = \mathbf{E}^{\hat{\mathbf{P}}} \left[ \int_s^t g_u(\hat{\Lambda}^m) \mathscr{A}^{\delta_m}_u(\hat{\Lambda}^m_u)[\phi] \, \mathrm{d}u \right].$$

Instead of (5.1) in Proposition 5.2, we thus aim to prove

$$\lim_{m \to \infty} \mathbf{E}^{\hat{\mathbf{P}}} \left[ \int_s^t g_u(\hat{\Lambda}^m) \mathscr{A}^{\delta_m}_u(\hat{\Lambda}^m_u)[\phi] \, \mathrm{d}u \right] = \mathbf{E}^{\hat{\mathbf{P}}} \left[ \int_s^t g_u(\hat{\Lambda}^\infty) \mathscr{A}^{\delta_\infty}_u(\hat{\Lambda}^\infty_u)[\phi] \, \mathrm{d}u \right]. \quad (5.2)$$

**Comparing terms**    To establish (5.2), recall from (2.13) the definition $\mathscr{A}^\delta_t(\lambda)[\phi] = \lambda[\mathscr{L}^\delta_t(\lambda)\phi]$ for $\lambda \in \mathbf{M}^+_1$ to write

$$\int_s^t g_u(\hat{\Lambda}^\infty) \mathscr{A}^{\delta_\infty}_u(\hat{\Lambda}^\infty_u)[\phi] \, \mathrm{d}u - \int_s^t g_u(\hat{\Lambda}^m) \mathscr{A}^{\delta_m}_u(\hat{\Lambda}^m_u)[\phi] \, \mathrm{d}u \quad (5.3)$$

$$= \int_s^t \left( g_u(\hat{\Lambda}^\infty) \hat{\Lambda}^\infty_u [\mathscr{L}^{\delta_\infty}_u(\hat{\Lambda}^\infty_u)\phi] - g_u(\hat{\Lambda}^m) \hat{\Lambda}^m_u [\mathscr{L}^{\delta_m}_u(\hat{\Lambda}^m_u)\phi] \right) \mathrm{d}u$$

$$= \hat{E}_m + \hat{F}_m + \hat{G}_m \,,$$

where for each $m \in \mathbb{N}$,

$$\hat{E}_m := \int_s^t g_u(\hat{\Lambda}^\infty) \Big( \hat{\Lambda}_u^\infty [\mathscr{L}_u^{\delta_\infty} (\hat{\Lambda}_u^\infty) \phi] - \hat{\Lambda}_u^m [\mathscr{L}_u^{\delta_\infty} (\hat{\Lambda}_u^\infty) \phi] \Big) \, \mathrm{d}u \,, \tag{5.4}$$

$$\hat{F}_m := \int_s^t \Big( g_u(\hat{\Lambda}^\infty) \hat{\Lambda}_u^m [\mathscr{L}_u^{\delta_\infty} (\hat{\Lambda}_u^\infty) \phi] - g_u(\hat{\Lambda}^m) \hat{\Lambda}_u^m [\mathscr{L}_u^{\delta_\infty} (\hat{\Lambda}_u^\infty) \phi] \Big) \, \mathrm{d}u \,, \tag{5.5}$$

$$\hat{G}_m := \int_s^t g_u(\hat{\Lambda}^m) \hat{\Lambda}_u^m \Big[ \big( \mathscr{L}_u^{\delta_m} (\hat{\Lambda}_u^\infty) - \mathscr{L}_u^{\delta_m} (\hat{\Lambda}_u^m) \big) \phi \Big] \, \mathrm{d}u \,. \tag{5.6}$$

Let us also define

$$E_m := \mathbf{E}^{\hat{\mathbf{P}}}[\hat{E}_m] \,, \qquad F_m := \mathbf{E}^{\hat{\mathbf{P}}}[\hat{F}_m] \,, \qquad G_m := \mathbf{E}^{\hat{\mathbf{P}}}[\hat{G}_m] \,. \tag{5.7}$$

Via (5.3), we see that to establish (5.2), it is sufficient to show that

$$\lim_{m \to \infty} |E_m| + \lim_{m \to \infty} |F_m| + \lim_{m \to \infty} |G_m| = 0 \,.$$

Below we bound the individual contributions of $E_m$, $F_m$ and $G_m$. The strategy we use broadly follows [10]. Before doing this we transfer the findings from Proposition 3.3 under the sequence $(\mathbf{P}_{x_0}^{\infty, m})_{m \in \mathbb{N}}$ with narrow limit $\mathbf{P}_{x_0}^{\infty, \infty}$ to the sequence of random variables $(\hat{\Lambda}^m)_{m \in \mathbb{N}}$ with $\hat{\mathbf{P}}$-a.s. limit $\hat{\Lambda}^\infty$.

## 5.1 Bounding the individual contributions

For later reference, note that the emergence of regularity in Theorem 3.1 remains valid for the random variables $\hat{\Lambda}^m$ with $m \in \mathbb{N} \cup \{\infty\}$.

**Lemma 5.3** | *Impose Assumptions 1.1, 1.2 and 1.3. Then for all $m \in \mathbb{N} \cup \{\infty\}$ and each $t \in (0, T]$, we have $\hat{\mathbf{P}}$-a.s. that $\hat{\Lambda}_t^m \ll \mathrm{d}x$. In fact, there exist exponents*

$$1 < q < \infty \text{ and } 1 < r < \infty, \ 0 < \gamma < 1, \tag{5.8}$$

*a constant $c = c_{\text{Thm. 3.1}} < \infty$ and for each $m \in \mathbb{N} \cup \{\infty\}$ strongly a measurable function*

$$\hat{p}^m : (0, T] \to \mathbb{L}^q \Big( (\hat{\Omega}, \hat{\mathcal{F}}, \hat{\mathbf{P}}); \mathbb{L}^r(\mathbb{R}^d, \mathrm{d}x) \Big) \,, \tag{5.9}$$

*which jointly satisfy for all $m \in \mathbb{N} \cup \{\infty\}$ the bound*

$$\Big\| \|\hat{p}^m(t)\|_{\mathbb{L}^r(\mathrm{d}x)} \Big\|_{\mathbb{L}^q(\hat{\mathbf{P}})} \le c(1 \wedge t)^{-\gamma} \quad \text{for all } t \in (0, T] \tag{5.10}$$

*and for which for each $t \in (0, T]$, we have $\hat{\mathbf{P}}$-a.s. that*

$$\mathrm{d}\hat{\Lambda}_t(\hat{\omega})/\mathrm{d}x = \hat{p}^m(t)(\hat{\omega}) \quad \text{in } \mathbb{L}^r(\mathrm{d}x) \,. \tag{5.11}$$

*As a consequence, $\mathrm{d}\hat{\Lambda}_t^m(\hat{\omega})/\mathrm{d}x = \hat{p}^m(t)(\hat{\omega})$ in $\mathbb{L}^r(\mathrm{d}x)$ $(\mathrm{d}t \otimes \mathrm{d}\hat{\mathbf{P}})$-a.e.*

**Proof** Because $\mathrm{Law}_{\hat{\mathbf{P}}}(\hat{\Lambda}^m) = \mathrm{Law}_{\mathbf{P}_{x_0}^{\infty,m}}(\Lambda)$ for all $m \in \mathbb{N}\cup\{\infty\}$ by the Skorohod representation in (4.17), and because by (3.1) we have that $\mathsf{H}_r^w(\mathbb{R}^d) \hookrightarrow \mathbb{L}^r(\mathrm{d}x)$, the result from Proposition 3.3 transfers to $\hat{\Lambda}^m$ under $\hat{\mathbf{P}}$. □

To streamline notation, write $\hat{\Lambda}^{m,\delta} := (\hat{\Lambda}^m)^\delta$ for $\delta > 0$ and $m \in \mathbb{N}\cup\{\infty\}$, where $(\hat{\Lambda}^m)^\delta = \hat{\Lambda}^m * h_\delta$ is defined in (1.3). This makes sense because $\hat{\mathbf{P}}$-a.s. we have that $\hat{\Lambda}_t^m \in \mathbf{M}_1^+(\mathbb{R}^d)$ simultaneously for all $t \in [0, T]$ and all $m \in \mathbb{N}\cup\{\infty\}$.

**Lemma 5.4** | *Under Assumptions 1.1, 1.3 and 1.7, fix $\phi \in \mathbf{C}_c^\infty(\mathbb{R}^d)$ and a bounded measurable $g : [0, T] \times \Omega \to \mathbb{R}$. Then for any $\varepsilon > 0$, there exists $m_E(\varepsilon) \in \mathbb{N}$ such that $|E_m| < \varepsilon$ for all $m > m_E(\varepsilon)$. In other words, $\lim_{m\to\infty} |E_m| = 0$.*

**Proof** Fix $\varepsilon > 0$ and let $0 \leq s \leq t \leq T$. Recall from (2.9) that $\mathcal{L}^{\delta\infty} = \mathcal{L}^0 = \mathcal{L}$. With $c_g := \|g\|_\infty$ we can estimate the term $\hat{E}_m$ from (5.4) by

$$|\hat{E}_m| \leq c_g \int_s^t |(\hat{\Lambda}_u^\infty - \hat{\Lambda}_u^m)[\mathcal{L}_u(\hat{\Lambda}_u^\infty)\phi]|\,\mathrm{d}u\,. \tag{5.12}$$

Taking expectations under $\hat{\mathbf{P}}$ in (5.12), we use (5.11) to represent the $\mathrm{d}(\hat{\Lambda}_t^\infty - \hat{\Lambda}_u^m)$-integral in terms of an integral against $(\hat{p}^\infty(t) - \hat{p}^m(t))\,\mathrm{d}x$ and apply Hölder's inequality to find that

$$\mathbf{E}^{\hat{\mathbf{P}}}[\hat{E}_m] \leq c_g\,\mathbf{E}^{\hat{\mathbf{P}}}\left[\int_s^t \|\mathcal{L}_u(\hat{\Lambda}_u^\infty)\phi\|_{\mathbb{L}^\infty(\mathrm{d}x)}\|\hat{p}_u^\infty - \hat{p}_u^m\|_{\mathbb{L}^1(\mathrm{d}x)}\,\mathrm{d}u\right]\,.$$

The estimates to control the term on the right-hand side are familiar from the proof of [10, Lem. 4.4]. First consider the case $s > 0$. We use Tonelli's theorem to interchange the expectation with the $\mathrm{d}u$-integral and apply Hölder's inequality with conjugate exponents $1 < q, q' < \infty$ for the $\hat{\mathbf{P}}$-expectation. This yields

$$\mathbf{E}^{\hat{\mathbf{P}}}[\hat{E}_{m,\delta',1}] \leq c_g \int_s^t \left\|\|\mathcal{L}_u(\hat{\Lambda}_u^\infty)\phi\|_{\mathbb{L}^\infty(\mathrm{d}x)}\right\|_{\mathbb{L}^{q'}(\hat{\mathbf{P}})}\left\|\|\hat{p}_u^\infty - \hat{p}_u^m\|_{\mathbb{L}^1(\mathrm{d}x)}\right\|_{\mathbb{L}^q(\hat{\mathbf{P}})}\,\mathrm{d}u\,.$$

The estimate (2.21) in Remark 2.10 shows that $\left\|\|\mathcal{L}_u(\hat{\Lambda}_u^\infty)\phi\|_{\mathbb{L}^\infty(\mathrm{d}x)}\right\|_{\mathbb{L}^{q'}(\hat{\mathbf{P}})} \leq c_{b,\sigma,\bar{\sigma}}\|\phi\|_2^*$. Since $\|\hat{p}_u^m\|_{\mathbb{L}^1(\mathrm{d}x)} = 1$ for all $m \in \mathbb{N}\cup\{\infty\}$ and $u \in (0, T]$, Proposition 4.5 and dominated convergence imply that we can find $m_0(\varepsilon) \in \mathbb{N}$ such that for all $m > m_0(\varepsilon)$ we have

$$E_m = \mathbf{E}^{\hat{\mathbf{P}}}[\hat{E}_m] < \varepsilon\,.$$

Now take $s = 0$. Splitting for an intermediate $s' \in (0, t]$, to be chosen below, in (5.12) the integral $\int_0^t = \int_0^{s'} + \int_{s'}^t$ gives the bound

$$|\hat{E}_m| \leq c_g \int_0^{s'} |(\hat{\Lambda}_u^\infty - \hat{\Lambda}_u^m)[\mathcal{L}_u(\hat{\Lambda}_u^\infty)\phi]|\,\mathrm{d}u + c_g \int_{s'}^t |(\hat{\Lambda}_u^\infty - \hat{\Lambda}_u^m)[\mathcal{L}_u(\hat{\Lambda}_u^\infty)\phi]|\,\mathrm{d}u\,. \tag{5.13}$$

For the first term in (5.13) we can use estimate (2.21) to bound it by $c_g c_{b,\sigma,\bar{\sigma}}\|\phi\|_2^* s'$. Choosing $s' > 0$ sufficiently small, we can make the contribution of the first term in (5.13) less than $\varepsilon/2$. For the second term in (5.13) we now proceed as in the case $s > 0$ above to make it less than $\varepsilon/2$, too. In this way we see that $E_m < \varepsilon$ and this concludes the proof. □

Note that for all $(t, x, \lambda) \in [0, T] \times \mathbb{R}^d \times \mathbf{M}_1^+$, $\delta > 0$ and $\phi \in \mathcal{S}$, we have

$$\left( \mathscr{L}_t(\lambda) - \mathscr{L}_t^\delta(\lambda) \right) \phi(x) = \left( b_t(x, \lambda) - b_t^\delta(x, \lambda^\delta) \right) \cdot \nabla \phi(x). \tag{5.14}$$

This is immediately seen by using the definitions of $\mathscr{L}$ and $\mathscr{L}^\delta$ in (2.9).

**Lemma 5.5** | *Under Assumptions 1.1, 1.3 and 1.7, fix $\phi \in \mathbf{C}_c^\infty(\mathbb{R}^d)$ and a bounded measurable $g : [0, T] \times \Omega \to \mathbb{R}$ such that for each $u \in [0, T]$, the real-valued function $g_u(\,\cdot\,) := g(u, \cdot)$ on $\Omega$ is continuous. Then for any $\varepsilon > 0$, there exists $m_F(\varepsilon) \in \mathbb{N}$ such that $|F_m| < \varepsilon$ for all $m > m_F(\varepsilon)$. In other words, $\lim_{m \to \infty} |F_m| = 0$.*

**Proof** Upon adding and subtracting $\int_s^t g_u(\hat{\Lambda}^m) \hat{\Lambda}_u^m [\mathscr{L}_u(\hat{\Lambda}_u^\infty) \phi] \, \mathrm{d}u$ to and from $\hat{F}_m$ in (5.5), we find by the triangle inequality that

$$\begin{aligned}
|\hat{F}_m| &\leq \left| \int_s^t \left( g_u(\hat{\Lambda}^\infty) - g_u(\hat{\Lambda}^m) \right) \hat{\Lambda}_u^m [\mathscr{L}_u(\hat{\Lambda}_u^\infty) \phi] \, \mathrm{d}u \right| \\
&\quad + \left| \int_s^t g_u(\hat{\Lambda}^m) \hat{\Lambda}_u^m \left[ \left( \mathscr{L}_u(\hat{\Lambda}_u^\infty) - \mathscr{L}_u^{\delta_m}(\hat{\Lambda}_u^\infty) \right) \phi \right] \, \mathrm{d}u \right| \\
&=: \hat{F}_{m,1} + \hat{F}_{m,2} .
\end{aligned} \tag{5.15}$$

For $\hat{F}_{m,1}$, we use the estimate in (2.21) from Remark 2.10 to see that

$$\hat{F}_{m,1} \leq c_{\phi, b, \sigma, \bar{\sigma}} \int_s^t |g_u(\hat{\Lambda}^\infty) - g_u(\hat{\Lambda}^m)| \, \mathrm{d}u , \tag{5.16}$$

where $c_{\phi, b, \sigma, \bar{\sigma}} := \|\phi\|_2^* c_{b, \sigma, \bar{\sigma}}$. We note from (4.18) that $\hat{\mathbf{P}}$-a.s. $\hat{\Lambda}^m \to \hat{\Lambda}^\infty$ in $\mathbf{C}([0, T]; \mathscr{H}_{-p})$ as $m \to \infty$. Since that the inclusion $\mathscr{H}_{-p} \hookrightarrow \mathcal{S}'$ is continuous, it follows that $\hat{\mathbf{P}}$-a.s. $\hat{\Lambda}^m \to \hat{\Lambda}^\infty$ in $\mathbf{C}([0, T]; \mathcal{S}')$ as $m \to \infty$. Since for all $u \in [s, t]$, the real-valued function $g_u(\,\cdot\,)$ is continuous on $\Omega = \mathbf{C}([0, T]; \mathcal{S}')$, we get $\hat{\mathbf{P}}$-a.s. that $g_u(\hat{\Lambda}^m) \to g_u(\hat{\Lambda}^\infty)$ in $\mathbb{R}$ as $m \to \infty$ simultaneously for all $u \in [s, t]$. As $g$ is uniformly bounded, dominated convergence thus yields $\int_s^t |g_u(\hat{\Lambda}^m) - g_u(\hat{\Lambda}^\infty)| \, \mathrm{d}u \to 0$ $\hat{\mathbf{P}}$-a.s. as $m \to \infty$. Taking expectations and using again dominated convergence shows via (5.16) that

$$\mathbf{E}^{\hat{\mathbf{P}}}[\hat{F}_{m,1}] \leq c_{\phi, b, \sigma, \bar{\sigma}} \, \mathbf{E}^{\hat{\mathbf{P}}} \left[ \int_s^t |g_u(\hat{\Lambda}^\infty) - g_u(\hat{\Lambda}^m)| \, \mathrm{d}u \right] < \varepsilon/2 \tag{5.17}$$

for all $m > m_{F,1}(\varepsilon)$.

For $\hat{F}_{m,2}$ in (5.15), we use (5.14) to see that

$$\begin{aligned}
\hat{F}_{m,2} &\leq \|g\|_\infty \int_s^t \left| \hat{\Lambda}_u^m \left[ \left( \mathscr{L}_u(\hat{\Lambda}_u^\infty) - \mathscr{L}_u^m(\hat{\Lambda}_u^{\infty, \delta_m}) \right) \phi \right] \right| \, \mathrm{d}u \\
&= \|g\|_\infty \int_s^t \left| \hat{\Lambda}_u^m \left[ \left( b_u(\hat{\Lambda}_u^\infty) - b_u^{\delta_m}(\hat{\Lambda}_u^{\infty, \delta_m}) \right) \cdot \nabla \phi \right] \right| \, \mathrm{d}u .
\end{aligned} \tag{5.18}$$

By Lemma 5.3, we can fix exponents

$$1 < q < \infty \text{ and } 1 < r < \infty, \ 0 < \gamma < 1 \tag{5.19}$$

as in (5.8) and a function $\hat{p}^m$ which satisfies the bound (5.10) and for which we have for all $m \in \mathbb{N} \cup \{\infty\}$ that

$$\mathrm{d}\hat{\Lambda}^m_t(\hat{\omega})/\mathrm{d}x = \hat{p}^m(t)(\hat{\omega}) \qquad \text{in } \mathbb{L}^r(\mathrm{d}x) \ (\mathrm{d}t \otimes \mathrm{d}\hat{\mathbf{P}})\text{-a.e.} \tag{5.20}$$

We denote by

$$K_\phi := \mathrm{cl}_{(\mathbb{R}^d, |\cdot|)}\{x \in \mathbb{R}^d : |\nabla\phi(x)| > 0\} \tag{5.21}$$

the support of the gradient of $\phi$.

Taking expectations under $\hat{\mathbf{P}}$ in (5.18), we use (5.20) and apply Hölder's inequality to the $\mathrm{d}\hat{\Lambda}^m_u$-integral with conjugate exponents $r, r'$ to find with $c_{g,\phi} := c_g\|\phi\|^*_1$ that

$$\mathbf{E}^{\hat{\mathbf{P}}}[\hat{F}_{m,2}] \le c_{g,\phi}\mathbf{E}^{\hat{\mathbf{P}}}\left[\int_s^t \|b_u(\hat{\Lambda}^\infty_u) - b_u^{\delta_m}(\hat{\Lambda}^{\infty,\delta_m}_u)\|_{\mathbb{L}^{r'}(K_\phi,\mathrm{d}x)}\|\hat{p}^m_u\|_{\mathbb{L}^r(\mathrm{d}x)} \,\mathrm{d}u\right].$$

The estimates to control the term on the right-hand side are familiar from the proof of [10, Lem. 4.4]. We use Tonelli's theorem to interchange the expectation with the $\mathrm{d}u$-integral and apply Hölder's inequality twice, with conjugate exponents $q, q'$ for the $\hat{\mathbf{P}}$-expectation and $w, w'$, for $1 < w < \infty$ which we choose below, for the $\mathrm{d}u$-integral. This, combined with the uniform norm bound in (3.6) for $\hat{p}^m_u$ yields

$$\begin{aligned}
\mathbf{E}^{\hat{\mathbf{P}}}[\hat{F}_{m,2}] &\le c_{g,\phi}\left(\int_s^t \left\|\|b_u(\hat{\Lambda}^\infty_u) - b_u^{\delta_m}(\hat{\Lambda}^{\infty,\delta_m}_u)\|_{\mathbb{L}^{r'}(K_\phi,\mathrm{d}x)}\right\|^{w'}_{\mathbb{L}^{q'}(\hat{\mathbf{P}})} \mathrm{d}u\right)^{\frac{1}{w'}} \\
&\quad \times \left(\int_s^t \left\|\|\hat{p}^m_u\|_{\mathbb{L}^r(\mathrm{d}x)}\right\|^w_{\mathbb{L}^q(\hat{\mathbf{P}})} \mathrm{d}u\right)^{\frac{1}{w}} \\
&\le c_{g,\phi}\left(\int_s^t \left\|\|b_u(\hat{\Lambda}^\infty_u) - b_u^{\delta_m}(\hat{\Lambda}^{\infty,\delta_m}_u)\|_{\mathbb{L}^{r'}(K_\phi,\mathrm{d}x)}\right\|^{w'}_{\mathbb{L}^{q'}(\hat{\mathbf{P}})} \mathrm{d}u\right)^{\frac{1}{w'}} \\
&\quad \times \left(\int_s^t c^w(1 \wedge u)^{-\gamma w} \,\mathrm{d}u\right)^{\frac{1}{w}}. 
\end{aligned} \tag{5.22}$$

Since $\gamma < 1$ by (5.19), we can choose $w > 1$ small enough to have

$$\int_s^t (1 \wedge u)^{-\gamma w} \,\mathrm{d}u \le \int_0^T (1 \wedge u)^{-\gamma w} \,\mathrm{d}u =: c_{w,\gamma} < \infty. \tag{5.23}$$

Then (5.22) and (5.23) imply with $c_1 := cc^{1/w}_{w,\gamma}c_{g,\phi}$ and an application of the triangle inequality that

$$\mathbf{E}^{\hat{\mathbf{P}}}[F_{m,2}] \le c_1\left(\int_s^t \left\|\|b_u(\hat{\Lambda}^\infty_u) - b_u^{\delta_m}(\hat{\Lambda}^{\infty,\delta_m}_u)\|_{\mathbb{L}^{r'}(K_\phi,\mathrm{d}x)}\right\|^{w'}_{\mathbb{L}^{q'}(\hat{\mathbf{P}})} \mathrm{d}u\right)^{\frac{1}{w'}} \tag{5.24}$$

$$\le c_1\left(\int_s^t \left\|\|b_u(\hat{\Lambda}^\infty_u) - b_u^{\delta_m}(\hat{\Lambda}^\infty_u)\|_{\mathbb{L}^{r'}(K_\phi,\mathrm{d}x)}\right\|^{w'}_{\mathbb{L}^{q'}(\hat{\mathbf{P}})} \mathrm{d}u\right)^{\frac{1}{w'}}$$

$$+ c_1\left(\int_s^t \left\|\|b_u^{\delta_m}(\hat{\Lambda}^\infty_u) - b_u^{\delta_m}(\hat{\Lambda}^{\infty,\delta_m}_u)\|_{\mathbb{L}^{r'}(K_\phi,\mathrm{d}x)}\right\|^{w'}_{\mathbb{L}^{q'}(\hat{\mathbf{P}})} \mathrm{d}u\right)^{\frac{1}{w'}}. \tag{5.25}$$

Since $\hat{\Lambda}^\infty$ is $\hat{\mathbf{P}}$-a.s. $\mathbf{M}^+_1$-valued, part 2) of Lemma B.2 implies $\hat{\mathbf{P}}$-a.s. that

$$\lim_{m\to 0} \|b_u(\hat{\Lambda}^\infty_u) - b_u^{\delta_m}(\hat{\Lambda}^\infty_u)\|_{\mathbb{L}^{r'}(K_\phi,\mathrm{d}x)} = 0$$

simultaneously for all $u \in [s, t]$. Next define $K'_\phi := K_\phi - K_h = \{x - x' : x \in K_\phi, x' \in K_h\}$, where $K_\phi$ denotes as in (5.21) the support of $\nabla \phi$ and $K_h := \mathrm{cl}_{(\mathbb{R}^d, |\cdot|)}\{x \in \mathbb{R}^d : |h(x)| > 0\}$; see Lemma B.2. By part 4) of Lemma B.2, the set $K'_\phi$ is compact and the linearity of the convolution operation shows that

$$\|b_u^{\delta_m}(\hat{\Lambda}_u^\infty) - b_u^{\delta_m}(\hat{\Lambda}_u^{\infty, \delta_m})\|_{\mathbb{L}^{r'}(K_\phi, \mathrm{d}x)} \leq \|b_u(\hat{\Lambda}_u^\infty) - b_u(\hat{\Lambda}_u^{\infty, \delta_m})\|_{\mathbb{L}^{r'}(K'_\phi, \mathrm{d}x)}. \quad (5.26)$$

To bound the right-hand side, we first use the uniform boundedness of $b$ from Assumption 1.7 to get the estimate

$$\int_{K'_\phi} |b_u(\hat{\Lambda}_u^\infty) - b_u(\hat{\Lambda}_u^{\infty, \delta_m})|^{r'} \, \mathrm{d}x \leq 2\|b\|_\infty^{r'-1} \int_{K_{\phi'}} |b_u(\hat{\Lambda}_u^\infty) - b_u(\hat{\Lambda}_u^{\infty, \delta_m})| \, \mathrm{d}x,$$

so that with $c_{b, r'} := \|b\|_\infty^{r'-1}$, we get

$$\|b_u(\hat{\Lambda}_u^\infty) - b_u(\hat{\Lambda}_u^{\infty, \delta_m})\|_{\mathbb{L}^{r'}(K_\phi, \mathrm{d}x)}^{r'} \leq 2c_{b, r'}\|b_u(\hat{\Lambda}_u^\infty) - b_u(\hat{\Lambda}_u^{\infty, \delta_m})\|_{\mathbb{L}^1(K'_\phi, \mathrm{d}x)}. \quad (5.27)$$

Since by Lemma 5.3 we have $(\mathrm{d}t \otimes \mathrm{d}\hat{\mathbf{P}})$-a.e. that $\hat{\Lambda}_t^\infty$ is in $\mathbb{L}^r(\mathrm{d}x)$, it follows from standard results that $\hat{\Lambda}_t^{\infty, \delta_m} \to \hat{\Lambda}_t^\infty$ in $\mathbb{L}^r(\mathrm{d}x)$ as $m \to \infty$, so Assumption 1.7 and (5.27) give

$$\lim_{m \to \infty} \|b_u(\hat{\Lambda}_u^\infty) - b_u(\hat{\Lambda}_u^{\infty, \delta_m})\|_{\mathbb{L}^{r'}(K_\phi, \mathrm{d}x)} = 0. \quad (5.28)$$

The integrable bound $2\|b\|_\infty \mathbb{1}_{[s,t] \times \hat{\Omega} \times K_\phi}$ lets us appeal to dominated convergence to conclude that there exists $m_{F,2}(\varepsilon)$ such that for all $m > m_{F,2}(\varepsilon)$,

$$\mathbf{E}^{\hat{\mathbf{P}}}[\hat{F}_{m,2}] < \varepsilon/2. \quad (5.29)$$

Now combining (5.15), (5.17) and (5.29) implies that

$$F_m = \mathbf{E}^{\hat{\mathbf{P}}}[\hat{F}_m] \leq \mathbf{E}^{\hat{\mathbf{P}}}[\hat{F}_{m,1}] + \mathbf{E}^{\hat{\mathbf{P}}}[\hat{F}_{m,2}] < \varepsilon$$

for all $m > m_G := \max\{m_{G,1}(\varepsilon), m_{G,2}(\varepsilon)\}$. $\qquad\square$

**Lemma 5.6** | *Under Assumptions 1.1, 1.3 and 1.7, fix $\phi \in \mathbf{C}_c^\infty(\mathbb{R}^d)$ and a bounded measurable $g : [0, T] \times \Omega \to \mathbb{R}$. Then for any $\varepsilon > 0$, there exists $m_G(\varepsilon) \in \mathbb{N}$ such that $|G_m| < \varepsilon$ for all $m > m_G(\varepsilon)$. In other words, $\lim_{m \to \infty} |G_m| = 0$.*

**Proof** From the definitions of $\hat{G}_m$ and $\mathscr{L}^{\delta_m}$ in (5.6) and (2.9), we get by the triangle inequality that

$$|\hat{G}_m| \leq \|g\|_\infty \int_s^t \left|\hat{\Lambda}_u^m\left[\left(\mathscr{L}_u^{\delta_m}(\hat{\Lambda}_u^\infty) - \mathscr{L}_u^{\delta_m}(\hat{\Lambda}_u^m)\right)\phi\right]\right| \mathrm{d}u \quad (5.30)$$

$$\leq \|g\|_\infty \int_s^t \left|\hat{\Lambda}_u^m\left[\left(b_u^{\delta_m}(\hat{\Lambda}_u^{\infty, \delta_m}) - b_u^{\delta_m}(\hat{\Lambda}_u^{m, \delta_m})\right) \cdot \nabla\phi\right]\right| \mathrm{d}u$$

$$+ \frac{1}{2}\|g\|_\infty \int_s^t \left|\hat{\Lambda}_u^m\left[\left(a_u(\hat{\Lambda}_u^\infty) - a_u(\hat{\Lambda}_u^m)\right) : \nabla^2\phi\right]\right| \mathrm{d}u$$

$$=: \hat{G}_{m,1} + \hat{G}_{m,2}.$$

We start with the term $\hat{G}_{m,1}$. It is similar to (5.18); the only difference is that here we have $b_u^{\delta_m}(\hat{\Lambda}_u^{\infty,\delta_m}) - b_u^{\delta_m}(\hat{\Lambda}_u^{m,\delta_m})$, while in (5.18) we have $b_u(\hat{\Lambda}_u^{\infty}) - b_u^{\delta_m}(\hat{\Lambda}_u^{\infty,\delta_m})$. With this replacement, we can follow the same steps as those leading from (5.18) to (5.24) and use the bound $\|b_u^{\delta_m}(\hat{\Lambda}_u^{\infty,\delta_m}) - b_u^{\delta_m}(\hat{\Lambda}_u^{m,\delta_m})\|_{\mathbb{L}^{r'}(K_\phi,\mathrm{d}x)} \leq \|b_u(\hat{\Lambda}_u^{\infty,\delta_m}) - b_u(\hat{\Lambda}_u^{m,\delta_m})\|_{\mathbb{L}^{r'}(K_\phi',\mathrm{d}x)}$ to obtain the estimate

$$\mathbf{E}^{\hat{\mathbf{P}}}[\hat{G}_{m,1}] \leq c_1 \left( \int_s^t \left\| \|b_u^{\delta_m}(\hat{\Lambda}_u^{\infty,\delta_m}) - b_u^{\delta_m}(\hat{\Lambda}_u^{m,\delta_m})\|_{\mathbb{L}^{r'}(K_\phi,\mathrm{d}x)} \right\|_{\mathbb{L}^q(\hat{\mathbf{P}})}^{w'} \mathrm{d}u \right)^{\frac{1}{w'}}$$

$$\leq c_1 \left( \int_s^t \left\| \|b_u(\hat{\Lambda}_u^{\infty,\delta_m}) - b_u(\hat{\Lambda}_u^{m,\delta_m})\|_{\mathbb{L}^{r'}(K_\phi,\mathrm{d}x)} \right\|_{\mathbb{L}^q(\hat{\mathbf{P}})}^{w'} \mathrm{d}u \right)^{\frac{1}{w'}},$$

where $c_1 = cc_{w,r}^{1/w}c_{g,\phi}$ and $K_\phi'$ are as defined in Lemma 5.5 and the Hölder conjugate pair $1 < w, w' < \infty$ is chosen as in (5.23).

Consider the case $s > 0$. Young's inequality implies for $1 \leq \xi < r$ and $m \in \mathbb{N}$ that

$$\|\hat{\Lambda}_t^{\infty,\delta_m} - \hat{\Lambda}_t^{m,\delta_m}\|_{\mathbb{L}^\xi(\mathrm{d}x)} \leq \|\hat{\Lambda}_t^{\infty} - \hat{\Lambda}_t^{m}\|_{\mathbb{L}^\xi(\mathrm{d}x)}. \tag{5.31}$$

Therefore, by (4.20) we get that $\hat{\mathbf{P}}$-a.s. simultaneously for all $u \in [s,t]$ that $\hat{\Lambda}_u^m \to \hat{\Lambda}_u^m$ in $\mathbb{L}^\xi(\mathrm{d}x)$ as $m \to \infty$ so that Assumption 1.7 and the arguments in (5.26)–(5.28) imply $\hat{\mathbf{P}}$-a.s. simultaneously for all $u \in [s,t]$ that $\|b_u(\hat{\Lambda}_u^m) - b_u(\hat{\Lambda}_u^{\infty})\|_{\mathbb{L}^1(K_\phi',\mathrm{d}x)} \to 0$ as $m \to \infty$. Therefore, dominated convergence the integrable bound $4\|b\|_\infty \mathbb{1}_{[s,t] \times \hat{\Omega} \times K_\phi'}$ shows that we can find $m_{G,1}(\varepsilon) \in \mathbb{N}$ such that

$$\mathbf{E}^{\hat{\mathbf{P}}}[\hat{G}_{m,1}] < \varepsilon/2 \tag{5.32}$$

whenever $m > m_{G,1}(\varepsilon)$.

We continue with the term $\hat{G}_{m,2}$ from (5.30). Let

$$K_\phi'' := \mathrm{cl}_{(\mathbb{R}^d, |\cdot|)}\{x \in \mathbb{R}^d : |\nabla^2\phi(x)| > 0\}$$

denote the support of the Hessian of $\phi$. Recall that $a = \sigma\sigma^{\mathsf{T}} + \bar{\sigma}\bar{\sigma}^{\mathsf{T}}$ as defined following (2.9). Using first that for all $m \in \mathbb{N}$, the process $\hat{\Lambda}^m$ is $\hat{\mathbf{P}}$-a.s. $\mathbf{M}_1^+$-valued, and then the Hölder-continuity of $\sigma, \bar{\sigma}$ from Assumption 1.3(H), we get

$$\lim_{m \to \infty} \left| \hat{\Lambda}_u^m \left[ \left( a_u(\hat{\Lambda}_u^{\infty}) - a_u(\hat{\Lambda}_u^m) \right) : \nabla^2\phi \right] \right| \leq \|\phi\|_2^* \lim_{m \to \infty} \sup_{x \in K_\phi''} |a_u(x, \hat{\Lambda}_u^{\infty}) - a_u(x, \hat{\Lambda}_u^m)| = 0.$$

So the integrable bound $(\|\sigma\|_\infty^2 + \|\bar{\sigma}\|_\infty^2)\mathbb{1}_{[s,t] \times \hat{\Omega} \times K_\phi''}$ lets us apply dominated convergence. There thus exists $m_{2,G}(\varepsilon) \in \mathbb{N}$ such that

$$\mathbf{E}^{\hat{\mathbf{P}}}[G_{m,2}] = \frac{1}{2}\|g\|_\infty \int_s^t \mathbf{E}^{\hat{\mathbf{P}}}\left[ \left| \hat{\Lambda}_u^m \left[ \left( a_u(\hat{\Lambda}_u^{\infty}) - a_u(\hat{\Lambda}_u^m) \right) : \nabla^2\phi \right] \right| \right] \mathrm{d}u < \varepsilon/2 \tag{5.33}$$

for all $m > m_{G,2}(\varepsilon)$.

Combining (5.30)–(5.33) implies for the case $s > 0$ that

$$G_m = \mathbf{E}^{\hat{\mathbf{P}}}[\hat{G}_m] \leq \mathbf{E}^{\hat{\mathbf{P}}}[\hat{G}_{m,1}] + \mathbf{E}^{\hat{\mathbf{P}}}[\hat{G}_{m,2}] < \varepsilon$$

for all $m > m_G := \max\{m_{G,1}(\varepsilon), m_{G,2}(\varepsilon)\}$.

Finally, consider the case $s = 0$. Splitting in (5.30) the integral $\int_0^t = \int_0^{s'} + \int_{s'}^t$ for an appropriately chosen $s' > 0$, lets us proceed as in the proof of Lemma 5.5 to control the contribution on the interval $[0, s')$ using the estimate in Remark 2.10, and controlling the contribution on $[s', t]$ as in the case $s > 0$ above. □

## 5.2 Proof of Proposition 5.2 and Theorem 1.9

We are ready to prove Proposition 5.2 which, to recall, claims that

$$\lim_{m \to \infty} \mathbf{E}^{\mathbf{P}_{x_0}^{\infty, k_m}} \left[ \int_s^t g_u(\Lambda) \mathscr{A}_u^{\delta_{k_m}}(\Lambda_u)[\phi] \, \mathrm{d}u \right] = \mathbf{E}^{\mathbf{P}_{x_0}^{\infty, \infty}} \left[ \int_s^t g_u(\Lambda) \mathscr{A}_u(\Lambda_u)[\phi] \, \mathrm{d}u \right] \tag{5.1}$$

for any $\phi \in \mathcal{S}$, all $s, t \in [0, T]$ with $s \leq t$, and all bounded measurable functions $g : [0, T] \times \Omega \to \mathbb{R}$ such that for each $u \in [0, T]$, the real-valued function $g_u(\cdot) = g(u, \cdot)$ on $\Omega$ is continuous.

**Proof of Proposition 5.2**  To avoid double subscripts, we continue using the notation from the almost sure representation in the beginning of Section 5 and write $\delta_m$ instead of $\delta_{k_m}$, and recall from before (1.18) that $\delta_\infty = 0$. Fix $\phi \in \mathcal{S}$ and $\varepsilon > 0$, and recall from (2.13) that $\mathscr{A}_u^{\delta_m}(\Lambda_u)[\phi] = \Lambda_u[\mathscr{L}_u^{\delta_m}(\Lambda_u)\phi]$ for all $m \in \mathbb{N} \cup \{\infty\}$. As seen after Proposition 5.2, in order to establish (5.1), we only need to prove (5.2), i.e. that there exists $m_0(\varepsilon)$ such that $m > m_0(\varepsilon)$ implies

$$\left| \mathbf{E}^{\hat{\mathbf{P}}} \left[ \int_s^t g_u(\hat{\Lambda}^m) \hat{\Lambda}_u^m [\mathscr{L}_u^{\delta_m}(\hat{\Lambda}_u^m)\phi] \, \mathrm{d}u - \int_s^t g_u(\hat{\Lambda}^\infty) \hat{\Lambda}_u^\infty [\mathscr{L}_u^{\delta_\infty}(\hat{\Lambda}_u^\infty)\phi] \, \mathrm{d}u \right] \right| < \varepsilon. \tag{5.34}$$

We first claim that there exists a function $\psi \in \mathbf{C}_c^\infty$ such that for all $m \in \mathbb{N} \cup \{\infty\}$, we have

$$\left| \mathbf{E}^{\hat{\mathbf{P}}} \left[ \int_s^t g_u(\hat{\Lambda}^m) \hat{\Lambda}_u^m [\mathscr{L}_u^{\delta_m}(\hat{\Lambda}_u^m)(\phi - \psi)] \, \mathrm{d}u \right] \right| < \varepsilon/3. \tag{5.35}$$

To see this, recall from Remark 2.10 the estimate (2.21), from which we deduce the bound

$$\left| \mathbf{E}^{\hat{\mathbf{P}}} \left[ \int_s^t g_u(\hat{\Lambda}^m) \hat{\Lambda}_u^m [\mathscr{L}_u^{\delta_m}(\hat{\Lambda}_u^m)(\phi - \psi)] \, \mathrm{d}u \right] \right| \leq \|g\|_\infty c_{b,\sigma,\bar{\sigma}}(t-s) \|\phi - \psi\|_2^*. \tag{5.36}$$

Since $\mathbf{C}_c^\infty$ is dense in $\mathcal{S}$, we can choose $\psi \in \mathbf{C}_c^\infty$ to make $\|\phi - \psi\|_2^*$ and hence (5.36) arbitrarily small.

Now fix $\psi \in \mathbf{C}_c^\infty$ satisfying (5.36) for all $m \in \mathbb{N} \cup \{\infty\}$. Then linearity and twice the triangle inequality show that to prove (5.34), it is sufficient to establish that

$$\left| \mathbf{E}^{\hat{\mathbf{P}}} \left[ \int_s^t \left( g_u(\hat{\Lambda}^{\delta_\infty}) \hat{\Lambda}_u^\infty [\mathscr{L}_u(\hat{\Lambda}_u^\infty)\psi] - g_u(\hat{\Lambda}^m) \hat{\Lambda}_u^m [\mathscr{L}_u^{\delta_m}(\hat{\Lambda}_u^m)\psi] \right) \mathrm{d}u \right] \right| \leq E_m + F_m + G_m < \varepsilon/3,$$

where $E_m$, $F_m$ and $G_m$ are defined in (5.4)–(5.6). But this is an immediate consequence of Lemmas 5.4–5.6. This shows (5.34) and concludes the proof. □

We can now proceed to the proof of our main result in Theorem 1.9. To make its proof more transparent, we recall an easy auxiliary lemma.

**Lemma 5.7** | *Impose Assumption 1.3. For all $\phi, \psi \in \mathcal{S}$ and $t \in [0, T]$, the mapping*

$$\lambda \mapsto \mathcal{Q}_t(\lambda)[\phi, \psi] = \lambda[\mathcal{R}_t(\lambda)\phi] \cdot \lambda[\mathcal{R}_t(\lambda)\psi] = \lambda[\bar{\sigma}_t^{\mathsf{T}}(\lambda)\nabla\phi] \cdot \lambda[\bar{\sigma}_t^{\mathsf{T}}(\lambda)\nabla\psi]$$

*defined in (2.14) is continuous on $\mathcal{P}_{\mathrm{wk}^*}$, the space of probability measures on $\mathcal{B}(\mathbb{R}^d)$ equipped with the narrow topology.*

**Proof** See [10, Lem. 4.8]. $\qquad\square$

The proof of Theorem 1.9 runs in parallel to that of [10, Thm. 1.6]. We must show that $\mathbf{P}_{x_0}^{\infty,\infty}$ solves the martingale problem $\mathrm{MP}(\nu, \mathscr{A}, \mathscr{Q})$, i.e. satisfies the requirements of Definition 2.3 with $\nu_0 = \nu$, $A = \mathscr{A}$ and $Q = \mathscr{Q}$.

**Proof of Theorem 1.9** By the construction of the sequence $(\mathbf{P}_{x_0}^{\infty,k})_{k\in\mathbb{N}}$ in (1.18) and by Proposition 2.8, $\mathbf{P}_{x_0}^{\infty,k}$ solves for each $k \in \mathbb{N}$ the martingale problem $\mathrm{MP}(\nu, \mathscr{A}^{\delta_k}, \mathscr{Q})$ with $\delta_k$ from (1.16), $\mathscr{A}^{\delta_k}$ from (2.13) and $\mathscr{Q}$ from (2.14). Thus the process $\Lambda$ is an $\mathcal{S}'$-valued $(\mathbb{F}, \mathbf{P}_{x_0}^{\infty,k})$-semimartingale with continuous trajectories and canonical decomposition given by (2.18) in terms of $\mathscr{A}^{\delta_k}$, and with martingale part $M \in \mathscr{M}_2^c(\mathbf{P}_{x_0}^{\infty,k}; \mathcal{S}')$ having the quadratic variation process (2.19). Hence for each $\phi \in \mathcal{S}$, the process $\Lambda[\phi] = (\Lambda_t[\phi])_{t\in[0,T]}$ is a continuous real-valued $(\mathbb{F}, \mathbf{P}_{x_0}^{\infty,k})$-semimartingale with canonical decomposition given for all $t \in [0, T]$ by

$$\Lambda_t[\phi] = \Lambda_0[\phi] + \left(\int_0^t \mathscr{A}_u^{\delta_k}(\Lambda_u)\,\mathrm{d}u\right)[\phi] + M_t^k[\phi] \tag{5.37}$$
$$= \Lambda_0[\phi] + \int_0^t \mathscr{A}_u^{\delta_k}(\Lambda_u)[\phi]\,\mathrm{d}u + M_t^k[\phi]\,,$$

where $M^k[\phi] = (M_t^k[\phi])_{t\in[0,T]} \in \mathscr{M}_2^c(\mathbf{P}_{x_0}^{\infty,k}; \mathbb{R})$ has the quadratic variation process

$$\langle M^k[\phi]\rangle_t = \left(\int_0^t \mathscr{Q}_u(\Lambda_u)\,\mathrm{d}u\right)[\phi, \phi] \tag{5.38}$$
$$= \int_0^t \mathscr{Q}_u(\Lambda_u)[\phi, \phi]\,\mathrm{d}u\,, \qquad t \in [0, T]\,.$$

Note that (5.38) takes a form slightly different from [10, Eqn. (4.47)] for the following reason. Since each $\mathbf{P}_{x_0}^{\infty,k}$ is the law of a limiting particle system (with smoothed drift), there is no contribution from $\mathscr{C}_t(\Lambda_t)[\phi, \phi]$. Indeed, this term is present only for the $n$-particle systems, i.e. under the measure $\mathbf{P}_{x_0}^{n,k}$ from Lemma 2.9; it vanishes in the narrow limit $\mathbf{P}_{x_0}^{n,k} \to \mathbf{P}_{x_0}^{\infty,k}$ as $n \to \infty$; compare (2.17) and (2.19).

**Step 1** Let $\mathbf{P}_{x_0}^{\infty,\infty}$ be a cluster point of $(\mathbf{P}_{x_0}^{\infty,k})_{k\in\mathbb{N}}$ and $(\mathbf{P}_{x_0}^{\infty,k_m})_{m\in\mathbb{N}}$ a subsequence with $\mathbf{P}_{x_0}^{\infty,k_m} \to \mathbf{P}_{x_0}^{\infty,\infty}$ narrowly in $\mathbf{M}_1^+(\mathbf{C}([0,T]; \mathcal{S}'))$ as $m \to \infty$. To avoid double subscripts, we write $m$ instead of $k_m$; in particular, we continue using the notation from the almost sure representation in the beginning of Section 5 and write $\mathbf{P}_{x_0}^{\infty,m}$ and $\delta_m$ instead of $\mathbf{P}_{x_0}^{\infty,k_m}$ and $\delta_{k_m}$. So $\delta_\infty = \delta_{\infty\infty} = 0$. Recall from (2.13) that $\mathscr{A}_s^{\delta_m}(\Lambda_s)[\phi] = \Lambda_s[\mathscr{L}_s^{\delta_m}(\Lambda_s)\phi]$, and that $\mathscr{A}_s^0(\Lambda_s)[\phi] = \mathscr{A}_s(\Lambda_s)[\phi]$. We need to verify that $\mathbf{P}_{x_0}^{\infty,\infty}$ is a solution of the martingale problem $\mathrm{MP}(\nu, \mathscr{A}, \mathscr{Q})$.

Since by Assumption 1.1, $X_0^{1,n} \stackrel{d}{=} \nu$ and $(X_0^{i,n})_{i=1,\ldots,n}$ are i.d.d. for all $n \in \mathbb{N}$, the law of large numbers gives that $\Lambda_0[\phi] \stackrel{d}{=} \nu[\phi]$ under $\mathbf{P}_{x_0}^{\infty,\infty}$ for all $\phi \in \mathcal{S}$. We therefore have part 1) of Definition 2.3.

By Lemma 2.9, the process $\Lambda$ is $\mathbf{M}_1^+$-valued under $\mathbf{P}_{x_0}^{\infty,\infty}$. Hence the estimate (2.21) from Remark 2.10 gives the bound

$$\int_0^T |\mathscr{A}_t(\Lambda_t)[\phi]| \, \mathrm{d}t = \int_0^T |\Lambda_t[\mathscr{L}_t(\Lambda_t)\phi]| \, \mathrm{d}t \leq c_{b,\sigma,\tilde{\sigma}} T \|\phi\|_{2*} < \infty \, ,$$

and so part 2) of Definition 2.3 is satisfied. For parts 3) and 4), let $F \in \mathcal{S}$ and $0 \leq s \leq t \leq T$. By Itô's formula and (2.16), we have for $m \in \mathbb{N}$

$$\begin{aligned}
&F(\Lambda_t[\phi]) - F(\Lambda_s[\phi]) - \int_s^t F'(\Lambda_u[\phi]) \mathscr{A}_u^{\delta_{m_k}}(\Lambda_u)[\phi] \, \mathrm{d}u \\
&\quad - \frac{1}{2} \int_s^t F''(\Lambda_u[\phi]) \mathscr{Q}_u(\Lambda)[\phi, \phi] \, \mathrm{d}u \\
&= \int_s^t F'(\Lambda_u[\phi]) \, \mathrm{d}M_u^{k_m}[\phi] \, .
\end{aligned} \tag{5.39}$$

The derivative $F'$ of $F \in \mathcal{S}$ is bounded so that the process $\int F'(\Lambda_u[\phi]) \mathrm{d}M_u^{k_m}[\phi]$ is an $(\mathbb{F}, \mathbf{P}_{x_0}^{\infty,k_m})$-martingale for each $m \in \mathbb{N}$. We claim that we can follow the same reasoning as in Step 1 of the proof of [10, Thm. 1.6] to show that

$$\begin{aligned}
\mathbf{E}^{\mathbf{P}_{x_0}^{\infty,\infty}} &\left[ \left( F(\Lambda_t[\phi]) - F(\Lambda_s[\phi]) - \int_s^t F'(\Lambda_u[\phi]) \mathscr{A}_u(\Lambda_u)[\phi] \, \mathrm{d}u \right.\right. \\
&\left.\left. - \frac{1}{2} \int_s^t F''(\Lambda_u[\phi]) \mathscr{Q}_u(\Lambda)[\phi, \phi] \, \mathrm{d}u \right) G(\Lambda) \right] = 0
\end{aligned} \tag{5.40}$$

for all continuous, bounded and $\mathcal{F}_s$-measurable $G : \Omega \to \mathbb{R}$.

**Step 2**  Conditioning (5.39) on $\mathcal{F}_s$ thus gives for all continuous, bounded, $\mathcal{F}_s$-measurable functions $G : \Omega \to \mathbb{R}$ that

$$\mathbf{E}^{\mathbf{P}_{x_0}^{\infty,k_m}} \left[ \left( \int_s^t F'(\Lambda_u[\phi]) \mathrm{d}M_u^{\delta_m}[\phi] \right) G(\Lambda) \right] = 0 \qquad \text{for } m \in \mathbb{N} \, . \tag{5.41}$$

For the individual terms on the left-hand side of (5.39), we note the following. Clearly, by the narrow convergence $\mathbf{P}_{x_0}^{\infty,k_m} \to \mathbf{P}_{x_0}^{\infty,\infty}$ and the continuity of $G$ and of the evaluation map $\phi \mapsto \Lambda[\phi]$, we get

$$\begin{aligned}
\lim_{m \to \infty} \mathbf{E}^{\mathbf{P}_{x_0}^{\infty,k_m}} &\left[ \left( F(\Lambda_t[\phi]) - F(\Lambda_s[\phi]) \right) G(\Lambda) \right] \\
&= \mathbf{E}^{\mathbf{P}_{x_0}^{\infty,\infty}} \left[ \left( F(\Lambda_t[\phi]) - F(\Lambda_s[\phi]) \right) G(\Lambda) \right] \, .
\end{aligned} \tag{5.42}$$

Consider now the map $g : [0,T] \times \Omega \to \mathbb{R}$ defined by $g_u(\Lambda) := g(u, \Lambda) := G(\Lambda) F'(\Lambda_u[\phi])$. It is clearly bounded and measurable. In addition, $\lambda \mapsto g_u(\lambda)$ is continuous on $\Omega$ for each $u \in [0, T]$. So $g$ satisfies the hypotheses of Proposition 5.2, and we get that

$$
\begin{aligned}
\lim_{m \to \infty} \mathbf{E}^{\mathbf{P}^{\infty, k_m}_{x_0}} & \left[ \left( \int_s^t F'(\Lambda_u[\phi]) \mathscr{A}_u^{\delta_m}(\Lambda_u)[\phi] \, \mathrm{d}u \right) G(\Lambda) \right] \\
&= \lim_{m \to \infty} \mathbf{E}^{\mathbf{P}^{\infty, k_m}_{x_0}} \left[ \int_s^t g_u(\Lambda) \mathscr{A}_u^{\delta_m}(\Lambda_u)[\phi] \, \mathrm{d}u \right] \\
&= \mathbf{E}^{\mathbf{P}^{\infty, \infty}_{x_0}} \left[ \int_s^t g_u(\Lambda) \mathscr{A}_u(\Lambda_u)[\phi] \, \mathrm{d}u \right] \\
&= \mathbf{E}^{\mathbf{P}^{\infty, \infty}_{x_0}} \left[ \left( \int_s^t F'(\Lambda_u[\phi]) \mathscr{A}_u(\Lambda_u)[\phi] \, \mathrm{d}u \right) G(\Lambda) \right].
\end{aligned}
\tag{5.43}
$$

In addition, the function $\lambda \mapsto F''(\lambda_u[\phi]) G(\lambda)$ is bounded, measurable and continuous on $\Omega$ for each $u \in [0, T]$. Hence applying Lemma 5.7 shows that

$$
\begin{aligned}
\lim_{m \to \infty} \mathbf{E}^{\mathbf{P}^{\infty, k_m}_{x_0}} & \left[ \left( \frac{1}{2} \int_s^t F''(\Lambda_u[\phi]) \mathscr{Q}_u(\Lambda_u)[\phi, \phi] \, \mathrm{d}u \right) G(\Lambda) \right] \\
&= \mathbf{E}^{\mathbf{P}^{\infty, \infty}_{x_0}} \left[ \left( \frac{1}{2} \int_s^t F''(\Lambda_u[\phi]) \mathscr{Q}_u(\Lambda_u)[\phi, \phi] \, \mathrm{d}u \right) G(\Lambda) \right].
\end{aligned}
\tag{5.44}
$$

Starting from (5.39), multiplying both sides by $G(\Lambda)$, taking expectations with respect to $\mathbf{P}^{\infty, k_m}_{x_0}$, passing to the limit $m \to \infty$ and using (5.41)–(5.44), we therefore see that for all $F \in \mathcal{S}$ and all continuous, bounded and $\mathcal{F}_s$-measurable $G : \Omega \to \mathbb{R}$, we have that

$$
\begin{aligned}
\mathbf{E}^{\mathbf{P}^{\infty, \infty}_{x_0}} & \left[ \left( F(\Lambda_t[\phi]) - F(\Lambda_s[\phi]) - \int_s^t F'(\Lambda_u[\phi]) \mathscr{A}_u(\Lambda_u)[\phi] \, \mathrm{d}u \right. \right. \\
& \left. \left. - \frac{1}{2} \int_s^t F''(\Lambda_u[\phi]) \mathscr{Q}_u(\Lambda_u)[\phi, \phi] \, \mathrm{d}u \right) G(\Lambda) \right] = 0.
\end{aligned}
\tag{5.45}
$$

By a standard functional monotone class argument, see for instance Ethier and Kurtz [13, Sec. 3, Eqn. (3.4) and Thm. A.4.3], (5.45) can be extended from continuous, bounded and $\mathcal{F}_s$-measurable $G$ to general bounded and $\mathcal{F}_s$-measurable $G$, and this shows that the process $N^{F, \phi} := (N_t^{F, \phi})_{t \in [0, T]}$ defined by

$$
\begin{aligned}
N_t^{F, \phi} := F(\Lambda_t[\phi]) - F(\Lambda_0[\phi]) - \int_0^t F'(\Lambda_u[\phi]) \mathscr{A}_u(\Lambda_u)[\phi] \, \mathrm{d}u \\
- \frac{1}{2} \int_0^t F''(\Lambda_u[\phi]) \mathscr{Q}_u(\Lambda_u)[\phi, \phi] \, \mathrm{d}u, \qquad t \in [0, T],
\end{aligned}
\tag{5.46}
$$

is a continuous real-valued $(\mathbb{F}, \mathbf{P}^{\infty, \infty}_{x_0})$-martingale. Since $F \in \mathcal{S}$ is arbitrary, standard results on the characterization of diffusion processes, see e.g. Revuz and Yor [32, Prop. VII.2.4], show that for each $\phi \in \mathcal{S}$, the process $M^\phi = (M_t^\phi)_{t \in [0, T]}$ defined by

$$
M_t^\phi := \Lambda_t[\phi] - \Lambda_0[\phi] - \int_0^t \mathscr{A}_u(\Lambda_u)[\phi] \, \mathrm{d}u, \qquad t \in [0, T],
\tag{5.47}
$$

is a continuous $(\mathbb{F}, \mathbf{P}^{\infty, \infty}_{x_0})$-martingale with quadratic variation process

$$
\langle M^\phi \rangle_t = \int_0^t \mathscr{Q}_u(\Lambda_u)[\phi, \phi] \, \mathrm{d}u, \qquad t \in [0, T].
\tag{5.48}
$$

Since $\phi$ is arbitrary, Lemma 2.2 implies that $\Lambda$ is an $\mathcal{S}'$-valued $(\mathbb{F}, \mathbf{P}_{x_0}^{\infty,\infty})$-semimartingale with continuous (in $\mathcal{S}'$) trajectories and decomposition given by (1.19). Moreover, from (2.4), Remark 2.4 and polarization, we obtain the quadratic variation process $\langle\!\langle M \rangle\!\rangle$ in (1.20) directly from (5.48). This establishes parts 3) and 4) of Definition 2.3 and shows that $\mathbf{P}_{x_0}^{\infty,\infty}$ solves $\mathrm{MP}(\nu, \mathscr{A}, \mathscr{Q})$.

The rest of Theorem 1.9 now follows directly from Proposition 3.3. $\qquad\square$

# 6 EXAMPLES

In this section, we give two motivating examples of drifts that satisfy Assumption 1.7. In Section 7, below, we discuss Assumption 1.7 and thus also these examples can be generalized.

**Convolution with an irregular kernel**　　Let $k : [0, T] \times \mathbb{R}^d \to \mathbb{R}^d$ be a kernel that satisfies

$$\sup_{t \in [0,T]} \|k_t\|_{\mathbb{L}^\infty(\mathrm{d}x)} < \infty \,. \tag{6.1}$$

Consider the drift

$$b_t(x, \mu) \coloneqq (\mu * k_t)(x) \,, \tag{6.2}$$

where $*$ denotes as in (1.3) the usual convolution. We immediately see from Young's inequality that $b$ satisfies Assumption 1.7($\mathrm{C}_w$).

In fact, the boundedness requirement (6.1) can be generalized in several ways. While we didn't make the precise value of $r$ in Theorem 3.1 explicit, this can be done by carefully studying the proofs in [8]. The above argument remains valid for kernels $k$ for which $\sup_{t \in [0,T]} \|k_t\|_{\mathbb{L}^{r'}(\mathrm{d}x)} < \infty$ where $r'$ satisfies the conjugacy relation $1 = 1/r + 1/r'$. Indeed, in this case $b$ in (6.2) is by Young's inequality bounded.

We discuss in Section 7 how the case of unbounded $b$ can be treated.

**Local dependence on the density: Drift coefficients of Nemytskii-type**　　Consider a bounded function $B : [0, T] \times \mathbb{R}^d \times \mathbb{R}^{\geq 0} \to \mathbb{R}^d$ for which $(x, x') \mapsto B_t(x, x')$ is continuous on $\mathbb{R}^d \times \mathbb{R}_{\geq 0}$ for (almost) every $t \in [0, T]$. We define $b : [0, T] \times \mathbb{R}^d \times \mathbf{M}_1^+(\mathbb{R}^d) \to \mathbb{R}^d$ for $\mu \in \mathbf{M}_1^+(\mathbb{R}^d)$ with $\mu \ll \mathrm{d}x$ by

$$b_t(x, \mu) \coloneqq B_t(x, \delta_x[\mu]) \,, \tag{6.3}$$

where $\delta_x[\mu] = f_\mu(x)$ and $f_\mu \coloneqq \mathrm{d}\mu/\mathrm{d}x$. This is well defined $\mathrm{d}x$-a.e. since the density $f_\mu$ is unique up to $\mathrm{d}x$-a.e. equivalence.

We now show that a drift $b$ as in (6.3) satisfies Assumption 1.7. Let $(\mu_k)_{k \in \mathbb{N}}$ be a sequence in $\mathbf{M}_1^+$ with $\mu_k \ll \mathrm{d}x$ and let $\mu \in \mathbf{M}_1^+$ with $\mathrm{d}\mu \ll \mathrm{d}x$. Suppose that $f_{\mu_k} \to f_\mu$ in $\mathbb{L}^1(\mathrm{d}x)$ as $k \to \infty$. We can then find a subsequence $(f_{\mu^{k_m}})_{m \in \mathbb{N}}$ for which we have that $f_{\mu^{k_m}} \to f_\mu$ $\mathrm{d}x$-a.e. as $m \to \infty$. Now let $K \subseteq \mathbb{R}^d$ be compact. By the continuity of $B_t$, we have that $B_t(x, f_{\mu_{k_m}}(x)) \to b(x, f_\mu(x))$ for almost every $x \in K$ as $m \to \infty$. Let $M = \sup_{(x, x') \in \mathbb{R} \times \mathbb{R}_{\geq 0}} |B_t(x, x')| < \infty$. Then $2M \mathbf{1}_K \in L^1(\mathbb{R}^d, \mathrm{d}x)$ serves as an integrable bound to apply dominated convergence and conclude that

$$\int_K \left| B_t\big(x, f_{\mu_{k_m}}(x)\big) - B_t\big(x, f_\mu(x)\big) \right| \mathrm{d}x \longrightarrow 0 \qquad \text{as } m \to \infty.$$

To pass from the subsequence result to the full sequence, suppose for contradiction that for some $\varepsilon > 0$ we can find a subsequence $(f_{k_j})_{j \in \mathbb{N}}$ such that

$$\int_K \left| B_t\big(x, f_{\mu_{k_j}}(x)\big) - B_t\big(x, f_\mu(x)\big) \right| dx > \varepsilon \qquad \text{for all } j \in \mathbb{N} \,. \tag{6.4}$$

Since $(f_{k_j})_{j \in \mathbb{N}}$ converges to $f$ in $\mathbb{L}^1(\mathrm{d}x)$, we can extract from $(f_{k_j})_{j \in \mathbb{N}}$ a further subsequence $(f_{k_{j_m}})_{m \in \mathbb{N}}$ which converges to $f$ $\mathrm{d}x$-a.e. The previous argument now gives us that $\int_K |B_t(x, f_{\mu_{k_{j_m}}}(x)) - B_t(x, f_\mu(x))| \, \mathrm{d}x \to 0$ as $m \to \infty$. Since this is a contradiction to (6.4), we deduce that

$$\| B_t(x, \delta_x[\mu_k]) - B_t(x, \delta_x[\mu]) \|_{\mathbb{L}^1(K, \mathrm{d}x)} \longrightarrow 0 \qquad \text{as } k \to \infty \,,$$

and this shows that $b$ in (6.3) satisfies Assumption 1.7.

# 7 DISCUSSION AND CONCLUDING REMARKS

To initiate the discussion, let us recapitulate the key steps in the proof of Theorem 1.9.

Via the $\delta_k$-smoothed $n$-particle systems (1.1) and (1.2), we obtained the solution $\mathbf{P}^{n,k}_{x_0}$ of the martingale problem $\mathrm{MP}(\nu_n, \mathscr{A}^{\delta_k}, \mathscr{Q} + \frac{1}{n}\mathscr{C})$. The sequence $(\mathbf{P}^{n,k}_n)_{n \in \mathbb{N}}$ is tight in $\Omega$, and cluster points $\mathbf{P}^{\infty,k}_{x_0}$ solve the $\delta_k$-smoothed martingale problem $\mathrm{MP}(\nu, \mathscr{A}^{\delta_k}, \mathscr{Q})$. This uses the results in [10] and requires that the coefficients $b_t^{\delta_k}(x, \mu^{\delta_k})$, $\sigma_t(x, \mu)$, $\bar{\sigma}_t(x, \mu)$ are continuous in the measure argument $\mu$ relative to narrow convergence. In principle, this type of continuity is restrictive; however for the smoothed coefficients, Lemma 2.7 shows that even if the unsmoothed drift $b_t(x, \mu)$ does not satisfy this type of continuity, the smoothed drift $b_t^{\delta_k}(x, \mu^{\delta_k})$ does under Assumption 1.7.

To remove the smoothing, i.e., to take $\delta_k \to 0$, we study the sequence $(\mathbf{P}^{\infty,k}_{x_0})_{k \in \mathbb{N}}$. To avoid placing more stringent assumptions on the unsmoothed drift $b$, a more careful choice of the topology with respect to which tightness is considered must be made. The key ingredient that facilitates such a choice is the *emergence of regularity*, which shows that $\Lambda$ under $\mathbf{P}^{\infty,k}_{x_0}$ satisfies by Theorem 3.1 a good regularity property uniformly in $k \in \mathbb{N}$. Using it, we first deduce that the sequence $(\mathbf{P}^{\infty,k}_{x_0})_{k \in \mathbb{N}}$ is tight in $\Omega_{p,1} = \mathbb{L}^p([0, T]; \mathbb{L}^1(\mathbb{R}^d))$, and then show that for a Skorohod-representation $(\hat{\Lambda}^k)_{k \in \mathbb{N} \cup \{\infty\}}$ of a convergent subsequence $(\mathbf{P}^{\infty,k_m}_{x_0})_{m \in \mathbb{N}}$ with limit $\mathbf{P}^{\infty,\infty}_{x_0}$, we get the improved mode of convergence $(\Lambda_t^m)_{t \in (0,T]} \to (\Lambda_t^\infty)_{t \in (0,T]}$ in $\mathbf{C}((0, T]; \mathbb{L}^\xi(\mathbb{R}^d))$ with convergence uniformly on compacts as $m \to \infty$, valid for appropriate $\xi > 1$. This strong mode of convergence which we deduce in Proposition 4.5 is enough to let us complete our task. Using it and that for each $k \in \mathbb{N}$ the measure $\mathbf{P}^{\infty,k}_{x_0}$ solves the martingale problem $\mathrm{MP}(\nu, \mathscr{A}^{\delta_k}, \mathscr{Q})$, we finally show that $\mathbf{P}^{\infty,\infty}_{x_0}$ solves $\mathrm{MP}(\nu, \mathscr{A}, \mathscr{Q})$.

Our approach builds on the previous two papers [8, 10] in this series. However, the way in which we obtain Theorem 1.9 is not by a straightforward application of these results. Instead, they can only be combined and via the improved convergence result, which itself is obtained through rather intricate arguments. This synthesis is novel.

In Section 6, we gave two examples of coefficients that appear in the literature but are not continuous relative to narrow convergence. The first is given by convolution with an irregular kernel, the second is that of Nemytskii-type coefficients that depend locally on the density $\mu$. Recent contributions that obtain results for these types of coefficients include Röckner and Zhang [33], Röckner et al. [34] for the case of (singular) interaction kernels, and Barbu and Röckner [3], Barbu and Roöckner [4] as well as Grube [17, 18] for the case Nemitskii-type coefficients. However, we should note that these contributions do not include the common noise $Z$. Its presence not only significantly complicates the analysis but likely precludes the use of their techniques to obtain results. In addition, in some of these works, a functional form of the coefficients is prescribed, and the assumptions placed on them are more stringent than what we demand in Assumption 1.7. In these more restricted settings, some of those papers offer uniqueness results, a question that we do not address here.

## 7.1 Generalizations

To conclude, let us mention two avenues for generalizing Theorem 1.9. The first is weakening the boundedness requirement in Assumption 1.7. Already in [8, Sec. 7] and [10, Sec. 6], it is mentioned that this extension is possible but requires more work. Let us not recapitulate the necessary steps in those works upon which the present paper builds, but instead mention how we use boundedness in the proof of Theorem 1.9. The boundedness assumption simplifies the estimates we use, e.g., for (5.13). There, we split the $\mathrm{d}u$-integral on $[0, t]$ between $[0, s']$ and $[s', t]$, and use boundedness of $b$ to control via Remark 2.10 the contribution on $[0, s']$. Instead, we could use the uniform $\mathbb{L}^r(\mathrm{d}x)$-regularity of $\hat{\Lambda}_t^m$ and $\hat{\Lambda}_t^\infty$ from Lemma 5.3 to obtain a more careful bound via Hölder's inequality.

The second is weakening Assumption 1.7 to consider only tight sequences $(\mu_k)_{k \in \mathbb{N} \cup \{\infty\}}$ of absolutely continuous measures relative to Lebesgue for which $\mathrm{d}\mu_n/\mathrm{d}x \to \mathrm{d}\mu_\infty/\mathrm{d}x$ in $\mathbb{L}^\xi(\mathbb{R}^d)$ for appropriate $\xi > 1$ as $n \to \infty$. Note that for such sequences we have that $\mathrm{d}\mu_n/\mathrm{d}x \to \mathrm{d}\mu_\infty/\mathrm{d}x$ in $\mathbb{L}^1(\mathbb{R}^d)$ as $n \to \infty$, but the converse is not true. This weakening is immediately possible. In fact, in anticipation of it, we invested extra effort into Proposition 4.5, the almost sure convergence of $(\hat{\Lambda}_t^m)_{t \in (0,T]} \to (\hat{\Lambda}_t^\infty)_{t \in (0,T]}$ in $\mathbf{C}((0,T]; \mathbb{L}^\xi(\mathbb{R}^d))$ as $m \to \infty$, for all $\xi$ satisfying $1 \le \xi < r$ with $r$ as in Theorem 3.1. With it, the way we use Assumption 1.7 following (5.31) in the proof of Lemma 5.6 shows that we do not need continuity of $b$ in the measure argument $\mu$ relative to convergence in $\mathbb{L}^1(\mathbb{R}^d)$, but instead only require continuity for tight sequences of measures relative to convergence in $\mathbb{L}^\xi(\mathbb{R}^d)$. The reason we formulated Assumption 1.7 as it is stated in Section 1.3 is in this way we can avoid the mention of the parameter $r$. Since Theorem 3.1 does not specify its possible values, we preferred the more concise yet less general way of presenting Assumption 1.7 and thus Theorem 1.9. However, by studying the proof of Theorem 3.1 in [8], precise values of $r$ in terms of the dimensions $d$ and the Hölder-exponent $\beta$ of $\sigma$, $\bar{\sigma}$ from Assumption 1.3 can be given. With these values of $r$, also $\xi$ can be determined.

## A BACKGROUND

### A.1 Functional analysis

**The spaces of Schwartz functions an tempered distributions.** Let $\alpha \in \mathbb{N}_0^d$ be a multi-index, define $|\alpha| = \alpha_1 + \cdots + \alpha_d$, write $\partial^\alpha := \partial_1^{\alpha_1} \cdots \partial_d^{\alpha_d}$ and set

$$\langle x \rangle := (1 + |x|^2)^{1/2} \quad \text{for } x \in \mathbb{R}^d. \tag{A.1}$$

For each $m \in \mathbb{N}_0$ and $\phi \in \mathbf{C}^\infty(\mathbb{R}^d)$, we define the seminorm

$$\|\phi\|_m^* := \max_{|\alpha| \leq m} \sup_{x \in \mathbb{R}^d} |\langle x \rangle^m \partial^\alpha \phi(x)|. \tag{A.2}$$

The *space of rapidly decreasing functions* or Schwartz functions is the vector space

$$\mathcal{S} := \mathcal{S}(\mathbb{R}^d) := \{\phi \in \mathbf{C}^\infty(\mathbb{R}^d) : \|\phi\|_m^* < \infty \text{ for all } m \in \mathbb{N}_0\}.$$

We endow $\mathcal{S}$ with the Fréchet topology generated by the family $\| \cdot \|_m^*$ for $m \in \mathbb{N}_0$, i.e. $\phi_n \to \phi$ in $\mathcal{S}$ if and only if $\|\phi_n - \phi\|_m^* \to 0$ when $n \to \infty$, for each $m \in \mathbb{N}$. As usual, a set $B \subseteq \mathcal{S}$ is said to be *bounded* if $\sup_{\phi \in B} \|\phi\|_m^* < \infty$ for all $m \in \mathbb{N}_0$.

It is well known that $\mathcal{S}$ is metrizable, complete and separable. Indeed, let $\mathbf{C}_c^\infty = \mathbf{C}_c^\infty(\mathbb{R}^d)$ be the Fréchet space of smooth and compactly supported functions $\mathbb{R}^d \to \mathbb{R}$; see Rudin [35, Sec. 6.2] for details. Then $\mathbf{C}_c^\infty \hookrightarrow \mathcal{S}$, which is to say that $\mathbf{C}_c^\infty$ embeds continuously into $\mathcal{S}$. Moreover, this embedding is dense. As a consequence we also have the continuous embedding $\mathcal{S} \hookrightarrow \mathbb{L}^r(dx) = \mathbb{L}^r(\mathbb{R}^d, dx)$ for all $r \in [1, \infty]$, which is dense if $r \in [1, \infty)$.

Let $\lambda : \mathcal{S} \to \mathbb{R}$ be a linear functional. We denote the duality pairing by $\langle \lambda; \phi \rangle_{\mathcal{S}' \times \mathcal{S}}$ for $\phi \in \mathcal{S}$, and usually write $\lambda[\phi] := \langle \lambda; \phi \rangle_{\mathcal{S}' \times \mathcal{S}}$. The functional $\lambda$ is *continuous* if $\phi_n \to \phi$ in $\mathcal{S}$ implies $\lambda[\phi_n] \to \lambda[\phi]$. Equivalently, there exist $C > 0$ and $m \in \mathbb{N}$ such that

$$|\lambda[\phi]| \leq C \|\phi\|_m^* \quad \text{for all } \phi \in \mathcal{S}(\mathbb{R}^d). \tag{A.3}$$

The minimal index $m$ validating this inequality is called the *order* of $\lambda$. Continuous linear functionals are called *(tempered) distributions*; they form the strong topological dual space

$$\mathcal{S}' := \mathcal{S}'(\mathbb{R}^d) := \{\lambda : \mathcal{S}(\mathbb{R}^d) \to \mathbb{R} : \lambda \text{ is continuous}\}.$$

Its strong topology is generated by uniform convergence on bounded sets, i.e. the topology is generated the family of seminorms $\|\lambda\|_{\mathcal{S}', B} := \sup_{\phi \in B} |\lambda[\phi]|$ as $B$ ranges over the bounded sets of $\mathcal{S}$.

The strong topology generates the Borel-$\sigma$-field on $\mathcal{S}'$, denoted $\mathcal{B}(\mathcal{S}'(\mathbb{R}^d))$ or simply $\mathcal{B}(\mathcal{S}')$. It can be shown that $\mathcal{B}(\mathcal{S}')$ equals the $\sigma$-algebra generated by the collection of sets $\{\lambda \in \mathcal{S}' : \lambda[\phi] < a\}$ for $\phi \in \mathcal{S}$ and $a \in \mathbb{R}$; see e.g. Kallianpur and Xiong [24, Thm. 3.1.1].

Let $f$ be a measurable function. We define

$$\lambda_f[\phi] := \int_{\mathbb{R}^d} \phi(x) f(x) \, dx, \tag{A.4}$$

for $\phi \in \mathcal{S}$, provided the integral exists. Consider $\lambda \in \mathcal{S}'$ and suppose that there exists a $f \in \mathbb{L}_{\text{loc}}^1(dx)$ such that $\lambda = \lambda_f$ in $\mathcal{S}'$, that is to say $\lambda[\phi] = \lambda_f[\phi]$ for all $\phi \in \mathcal{S}$. Then $f$ is determined uniquely up to d$x$-a.e. equivalence; see Hörmander [20, Thm. 1.2.5]. In this case it is customary to say that $\lambda$ *is a function*, or that $\lambda$ is representable by $f$. The *canonical embedding* $\mathcal{S} \hookrightarrow \mathcal{S}'$ mapping $f \in \mathcal{S}$ to $\lambda_f \in \mathcal{S}'$ is thus well-defined and, in fact, continuous. This inclusion allows us to identify $\mathcal{S}$ with a subset of $\mathcal{S}'$: If $\lambda \in \mathcal{S}'$ is represented by a function in $\mathcal{S}$, then we freely write $\lambda \in \mathcal{S}$. Similarly, we have the canonical embedding $\mathbb{L}^r(\mathrm{d}x) \hookrightarrow \mathcal{S}'$ for all $r \in [1, \infty]$, and if $\lambda \in \mathcal{S}'$ is represented by a function in $\mathbb{L}^r(\mathrm{d}x)$, we simply write $\lambda \in \mathbb{L}^r(\mathrm{d}x)$. In the same fashion, if $\mu$ is a measure on the Borel sets of $\mathbb{R}^d$, we may define

$$\lambda_\mu[\phi] := \int_{\mathbb{R}^d} \phi(x)\, \mu(\mathrm{d}x) \tag{A.5}$$

for $\phi \in \mathcal{S}$, provided the integral exists. In analogy to the above, we say that $\lambda \in \mathcal{S}'$ is a measure if there is a measure $\mu$ such that $\lambda = \lambda_\mu$ in $\mathcal{S}'$. Let $\mathbf{M}_1^+ = \mathbf{M}_1^+(\mathbb{R}^d)$ denote the set of probability measures on the Borel sets of $\mathbb{R}^d$. It is standard to verify that $\mathbf{M}_1^+ \subseteq \mathcal{S}'$ by virtue of (A.5). In fact, the inclusion is sequentially continuous if $\mathbf{M}_1^+$ carries the narrow topology; see Lemma A.3 below.

We recall some basics from functional analysis, notably the space of tempered distributions. The material in this section is standard and can be found, e.g., in Rudin [35, Ch. 6] and Simon [37, Ch. 6].

**Identification of $\mathcal{S}'$ with a sequence space**    The canonical embedding (A.4) lets us write

$$\mathcal{S}(\mathbb{R}^d) \hookrightarrow \mathbb{L}^2(\mathbb{R}^d) \hookrightarrow \mathcal{S}'(\mathbb{R}^d)\,.$$

This allows to give an elegant and useful identification of $\mathcal{S}(\mathbb{R}^d)$ with rapidly decreasing sequences, and dually of $\mathcal{S}'(\mathbb{R}^d)$ with slowly increasing sequences, which we describe by using Hermite functions. An efficient presentation of this classical identification is found in Simon [36], to which we also refer for more background.

Starting from the standard Gaussian density

$$g(x) := \frac{1}{\sqrt{2\pi}} \exp\left(-\frac{1}{2}x^2\right) \qquad \text{for } x \in \mathbb{R}\,,$$

we define the *Hermite functions* on $\mathbb{R}$ by

$$h_{k+1}(x) := \frac{(-1)^k}{\sqrt{k!}}\, g(x)^{-\frac{1}{2}} \left(\frac{\mathrm{d}}{\mathrm{d}x}\right)^k g(x) \qquad \text{for } k = 0, 1, 2, \dots$$

The above definition extends to $\mathbb{R}^d$ as follows. We define the Hermite functions on $\mathbb{R}^d$ by $h_k(x) := h_{k_1}(x_1) \cdots h_{k_d}(x_d)$ for each multi-index $k \in \mathbb{N}^d$ and $x = (x_1, \dots, x_d) \in \mathbb{R}^d$. The Hermite functions are elements of $\mathcal{S}$, i.e. $h_k \in \mathcal{S}$ for any $k \in \mathbb{N}^d$. Moreover, they are uniformly bounded, i.e. $\sup_{k \in \mathbb{N}^d} \|h_k\|_\infty \leq (2\pi)^{d/4}$; see Indritz [22, p. 981] and Simon [36, Sec. 3, Note 12].

Let $(\,\cdot\,,\cdot\,)$ denote the usual inner product on $\mathbb{L}^2(\mathrm{d}x)$, which we use to identify $\mathbb{L}^2(\mathrm{d}x)$ with its dual via the Riesz isomorphism. It is well known that the family $(h_k)_{k \in \mathbb{N}^d}$ is a complete orthonormal system in $\mathbb{L}^2(\mathbb{R}^d, \mathrm{d}x)$ with respect to the usual inner product; see e.g. Simon [37, Thm. 6.4.3]

Since $\mathcal{S} \subseteq \mathbb{L}^2(dx)$, the basis expansion

$$\phi = \sum_{k \in \mathbb{N}^d} (\phi; h_k) \, h_k \qquad \text{in } \mathbb{L}^2(dx)$$

is well defined for any $\phi \in \mathcal{S}$, and we write $\phi_k^\# := (\phi; h_k)$ with $k \in \mathbb{N}^d$ for the *Hermite coefficients* of $\phi$. Similarly, if $\lambda \in \mathcal{S}'$, we define the Hermite coefficient by $\lambda_k^\# := \lambda[h_k]$, which makes sense since $h_k \in \mathcal{S}$.

We are now ready to discuss the identification of $\mathcal{S}$ and $\mathcal{S}'$ with certain sequences as promised at the beginning of the subsection. We call a sequence $(a_k)_{k \in \mathbb{N}^d}$ of real numbers *rapidly decreasing* if $(\langle k \rangle^m a_k)_{k \in \mathbb{N}^d} \in \ell^2(\mathbb{N}^d)$ for each $m \in \mathbb{N}_0$, with $\langle k \rangle = (1 + |k|^2)^{1/2}$ as in (A.1). We call $(a_k)_{k \in \mathbb{N}^d}$ *slowly increasing* if there exist $m \in \mathbb{N}_0$ and $c_a < \infty$ such that $|a_k| \leq c_a \langle k \rangle^m$ for all $k \in \mathbb{N}^d$. Clearly, a rapidly decreasing sequence is slowly increasing.

**Lemma A.1** | 1) *Let $\phi \in \mathcal{S}(\mathbb{R}^d)$. Then $(\phi_k^\#)_{k \in \mathbb{N}^d}$ is rapidly decreasing, i.e. $(\langle k \rangle^m \phi_k^\#)_{k \in \mathbb{N}^d} \in \ell^2(\mathbb{N}^d)$ for any $m \in \mathbb{N}_0$. Conversely, any rapidly decreasing sequence $(a_k)_{k \in \mathbb{N}^d}$ defines a function $\psi \in \mathcal{S}(\mathbb{R}^d)$ by virtue of $\psi = \sum_{k \in \mathbb{N}^d} a_k h_k$ with convergence in the topology of $\mathcal{S}(\mathbb{R}^d)$.*

*2) Let $\lambda \in \mathcal{S}'(\mathbb{R}^d)$. Then $(\lambda_k^\#)_{k \in \mathbb{N}^d}$ is slowly increasing, i.e. there are $m \in \mathbb{N}_0$ and $c_\lambda < \infty$ such that $|\lambda_k^\#| \leq c_\lambda \langle k \rangle^m$ for $k \in \mathbb{N}^d$. Conversely, any slowly increasing sequence $(b_k)_{k \in \mathbb{N}^d}$ defines a distribution $\tau \in \mathcal{S}'$ by virtue of $\tau[\phi] := \sum_{k \in \mathbb{N}^d} \phi_k^\# b_k$ in $\mathbb{R}$, for each $\phi \in \mathcal{S}$.*

**Proof**    See e.g. [36, Thms. 1 and 2]. ∎

By combining parts 1) and 2) of Lemma A.1, it can be shown that for any $\lambda \in \mathcal{S}'$ and $\phi \in \mathcal{S}$, we have

$$\lambda[\phi] = \sum_{k \in \mathbb{N}^d} \lambda_k^\# \phi_k^\# \,. \tag{A.6}$$

Indeed, the sum in (A.6) converges because $(\phi_k^\#)_{k \in \mathbb{N}^d}$ is rapidly decreasing and $(\lambda_k^\#)_{k \in \mathbb{N}^d}$ is slowly increasing; see [36, Thm. 3].

**Intermediate Hilbert spaces**    Let $\lambda \in \mathcal{S}'$. By Lemma A.1, the Hermite coefficients $(\lambda_k^\#)_{k \in \mathbb{N}^d}$ are slowly increasing. There thus exists a $p \in \mathbb{R}$ such that

$$\|\lambda\|_{\mathcal{H}_p} := \left( \sum_{k \in \mathbb{N}^d} \left( \langle k \rangle^{p/d} \lambda_k^\# \right)^2 \right)^{\frac{1}{2}} \tag{A.7}$$

is finite. Moreover, if $\phi \in \mathcal{S}$, then $\|\phi\|_{\mathcal{H}_p}$ is finite for any $p \in \mathbb{R}$, since the Hermite coefficients $(\phi_k^\#)_{k \in \mathbb{N}^d}$ are rapidly decreasing.

For each $p \in \mathbb{R}$, we thus define

$$\mathcal{H}_p := \mathcal{H}_p(\mathbb{R}^d) := \left\{ f \in \mathcal{S}' : \|f\|_{\mathcal{H}_p} < \infty \right\}.$$

Since $\mathcal{S}$ is $\sigma(\mathcal{S}', \mathcal{S})$-dense in $\mathcal{S}'$, the space $(\mathcal{H}_p, \|\cdot\|_{\mathcal{H}_p})$ can equivalently be defined as the Hilbert space completion of $\mathcal{S}$ under the semi-norm $\|\cdot\|_{\mathcal{H}_p}$; so $\mathcal{S}$ is dense in $\mathcal{H}_p$ for any $p \in \mathbb{R}$.

On $\mathcal{H}_p$, (A.7) specifies a Hilbertian semi-norm, that is, a semi-norm satisfying the parallelogram identity. The *Hermite–Fourier space* $(\mathcal{H}_p, \|\cdot\|_{\mathcal{H}_p})$ is thus a Hilbert subspace of $\mathcal{S}'$, with semi-inner product $(\,\cdot\,;\cdot\,)_p$ obtained from $\|\cdot\|_{\mathcal{H}_p}$ by polarization. With this notation, we observe that $(\,\cdot\,;\cdot\,) = (\,\cdot\,;\cdot\,)_0$.

Let $\mathcal{H}'_p$ denote the (strong) topological dual of $\mathcal{H}_p$. By definition, it is equipped with the norm

$$\|\lambda\|_{\mathcal{H}'_p} := \sup\{|\lambda[\phi]| : \|\phi\|_{\mathcal{H}_p} \leq 1\}. \tag{A.8}$$

The pair $(\mathcal{H}'_p, \|\cdot\|_{\mathcal{H}'_p})$ is again a Hilbert space.

It is classical that the canonical inclusion $\mathcal{S}(\mathbb{R}^d) \hookrightarrow \mathcal{H}'_p$ defined by $f \mapsto \lambda_f[\,\cdot\,] = (\,\cdot\,,f)_0$ extends uniquely to a continuous linear map $J : \mathcal{H}_{-p}(\mathbb{R}^d) \to \mathcal{H}'_p(\mathbb{R}^d)$ which is an isometric isomorphism so that $\mathcal{H}_{-p} \simeq \mathcal{H}'_p$; see for instance Kallianpur and Xiong [24, Thm. 1.3.2(c)] and Demengel and Demengel [12, Prop. 4.10].

An important consequence of the above is the chain of continuously embedded spaces

$$\mathcal{S}(\mathbb{R}^d) \hookrightarrow \mathcal{H}_p(\mathbb{R}^d) \hookrightarrow \mathbb{L}^2(\mathbb{R}^d) \hookrightarrow \mathcal{H}_{-p}(\mathbb{R}^d) \hookrightarrow \mathcal{S}'(\mathbb{R}^d). \tag{A.9}$$

We call the parameter $p \geq 0$ the *regularity*, and note that $\mathbb{L}^2(\mathrm{d}x) = \mathcal{H}_0(\mathbb{R}^d)$ with equal norms. The terminology we employ is suggestive: It indicates that objects of negative regularity act on objects of positive regularity by *canonical duality*

From the definition (A.7), we see that $\|\cdot\|_{\mathcal{H}_{p_1}} \leq \|\cdot\|_{\mathcal{H}_{p_2}}$ for $p_1 \leq p_2$. Thus $\mathcal{H}_{p_2} \hookrightarrow \mathcal{H}_{p_1}$ whenever $p_1 \leq p_2$, and dually, $\mathcal{H}_{-p_1} \hookrightarrow \mathcal{H}_{-p_2}$; so the chain (A.9) can be expanded indefinitely. In fact, from the abstract definitions, we can obtain

$$\mathcal{S}(\mathbb{R}^d) = \bigcap_{p \in \mathbb{R}} \mathcal{H}_p(\mathbb{R}^d) \quad \text{and} \quad \mathcal{S}'(\mathbb{R}^d) = \bigcup_{p \in \mathbb{R}} \mathcal{H}'_p(\mathbb{R}^d). \tag{A.10}$$

Both equalities in (A.10) are to be understood not just as equalities of sets, but as equalities of topological spaces. For this, we endow $\bigcap_{p \in \mathbb{R}} \mathcal{H}_p(\mathbb{R}^d)$ with the limit topology and $\bigcup_{p \in \mathbb{R}} \mathcal{H}_p(\mathbb{R}^d)$ with the colimit topology. In other words, if $(\phi_n)_{n \in \mathbb{N}}$ is a sequence in $\mathcal{S}$, then $\phi_n \to \phi$ in $\mathcal{S}$ as $n \to \infty$ if and only if $\phi_n \to \phi$ as $n \to \infty$ in $\mathcal{H}_p$ for any $p \in \mathbb{R}$. Dually, if $(\lambda_n)_{n \in \mathbb{N}}$ is a sequence in $\mathcal{S}'$, then $\lambda_n \to \lambda$ as $n \to \infty$ in $\mathcal{S}'$ if and only if $\lambda_n \to \lambda$ as $n \to \infty$ in $\mathcal{H}_p$ for some $p \in \mathbb{R}$.

**Lemma A.2** | *The two families of semi-norms $\|\cdot\|_{\mathcal{H}_p}$ in (A.7) for $p \in \mathbb{N}_0$ and $\|\cdot\|^*_m$ in (A.2) for $m \in \mathbb{N}_0$ on $\mathcal{S}(\mathbb{R}^d)$ are equivalent. That is, for each $m$, there are a $p$ and a universal constant $c_m$ such that*

$$\|\phi\|^*_m \leq c_m \|\phi\|_{\mathcal{H}_p},$$

*and for each $p$, there are an $m$ and a universal constant $c_p$ such that*

$$\|\phi\|_{\mathcal{H}_p} \leq c_p \|\phi\|^*_m,$$

*both for all $\phi \in \mathcal{S}(\mathbb{R}^d)$. For each $m$, the number $p$ and the constant $c_m$ depend on $d$, and similarly for each $p$, so do $m$ and $c_p$.*

**Proof** See e.g. Holley and Stroock [19, Appendix, (A.14)–(A.20)]. □

In view of the definition of $\|\cdot\|^*_m$ in (A.2), the above lemma shows in particular that functions in $\mathcal{H}_p$ enjoy certain differentiability properties as well as polynomial decay up to a certain order. An immediate consequence of this is in particular that $\mathcal{H}_p(\mathbb{R}^d) \hookrightarrow \mathsf{Lip}(\mathbb{R}^d)$, the space of Lipschitz functions on $\mathbb{R}^d$, for some $p > 0$ large enough depending on $d$.

### A.2  Probability theory

**Spaces of probability measures.**  A special role is our study is taken by positive measures, and in particular probability measures. We thus take a moment to clarify how they relate to our exposition so far. For a topological space $X$, we write $\mathbf{M}_1^+(X)$ for the set of probability measures on $\mathcal{B}(X)$. If $X = \mathbb{R}^d$, then we oftentimes write $\mathbf{M}_1^+ := \mathbf{M}_1^+(\mathbb{R}^d)$.

The *narrow topology* $\tau_{\mathrm{wk}^*}$ is induced by duality with continuous bounded functions $\mathbf{C}_b(X)$. Specifically, if $(\lambda_m)_{m \in \mathbb{N}}$ is a sequence in $\mathbf{M}_1^+(X)$, then $\lambda_m \to \lambda$ as $m \to \infty$ relative to $\tau_{\mathrm{wk}^*}$ if $\lambda_m[f] \to \lambda[f]$ for all $f \in \mathbf{C}_b(X)$. We write $\mathcal{P}_{\mathrm{wk}^*}(X) := (\mathbf{M}_1^+(X), \tau_{\mathrm{wk}^*})$ for the space $\mathbf{M}_1^+(X)$ endowed with the narrow topology $\tau_{\mathrm{wk}^*}(X)$. If $X$ is Polish, then so is $\mathcal{P}_{\mathrm{wk}^*}(X)$; see, e.g. Aliprantis and Border [2, Thm. 15.15].

Since $\mathbf{M}_1^+ \subseteq \mathcal{S}'$, we can consider the trace topology of $\mathcal{S}'$ on $\mathbf{M}_1^+$. This topology is, however, different from narrow convergence. To discuss this difference, denote by $\mathbf{M}_{\leq 1}^+ = \mathbf{M}_{\leq 1}^+(\mathbb{R}^d)$ the space of subprobability measures on $\mathcal{B}(\mathbb{R}^d)$, i.e. positive measures of no more than unit mass. Evidently $\mathbf{M}_1^+ \subseteq \mathbf{M}_{\leq 1}^+ \subseteq \mathcal{S}'$. The set $\mathbf{M}_1^+$ is *not closed* in $\mathcal{P}_{\leq 1, s'}$; see e.g. [8, Rem. A.2]. The set $\mathbf{M}_{\leq 1}^+$ equipped with the trace topology of $\mathcal{S}'$ is denoted by $\mathcal{P}_{\leq 1, s'} := (\mathbf{M}_{\leq 1}^+, \tau_{s'})$. Specifically, if $(\lambda_m)_{m \in \mathbb{N}}$ is a sequence in $\mathbf{M}_{\leq 1}^+$, then $\lambda_m \to \lambda$ as $m \to \infty$ in $\mathcal{P}_{\leq 1, s'}$ if for all $\phi \in \mathcal{S}$, we have that $\lambda_m[\phi] \to \lambda[\phi]$.

**Lemma A.3** | *If $(\nu_m)_{m \in \mathbb{N}}$ is a sequence in $\mathbf{M}_1^+$, then:*
*1) If $\nu_m \to \nu$ in $\mathcal{P}_{\mathrm{wk}^*}$, then $\nu \in \mathbf{M}_1^+$ and $\nu_m \to \nu$ in $\mathcal{P}_{s'}$.*
*2) If $\nu_m \to \nu$ in $\mathcal{P}_{\leq 1, s'}$ and $\nu[\mathbb{R}^d] = 1$, then $\nu_m \to \nu$ in $\mathcal{P}_{\mathrm{wk}^*}$.*

**Proof**  Part 1) follows from the fact that $\mathcal{S} \subseteq \mathbf{C}_b$ and that the notions of weak and strong convergence in $\mathcal{S}'$ coincide for sequences; see Huang and Yan [21, Thm. 3.12]. Part 2) follows from classical properties of narrow convergence; see Klenke [27, Thm. 13.16].  □

## B  AUXILIARY RESULTS AND OMITTED PROOFS

### B.1  Deferred proofs from Section 2

We begin with an auxiliary result.

**Lemma B.1** | *Let $(\mu_k)_{k \in \mathbb{N} \cup \{\infty\}}$ be a sequence of probability measures on $\mathcal{B}(\mathbb{R}^d)$ such that $\mu_k \to \mu_\infty$ narrowly as $k \to \infty$. Let $h \in \mathbf{C}_c^\infty(\mathbb{R}^d)$ and define $\nu_k := \mu_k * h$ for $k \in \mathbb{N} \cup \{\infty\}$. Then $\nu_k \to \nu_\infty$ in $\mathbb{L}^p(\mathbb{R}^d)$ as $k \to \infty$ for every $1 \leq p < \infty$.*

**Proof**  For each fixed $x \in \mathbb{R}^d$, the function $y \mapsto h(x-y)$ is bounded and continuous. Hence by narrow convergence,

$$\nu_k(x) = \int_{\mathbb{R}^d} h(x-y) \, \mu_k(\mathrm{d}y) \; \longrightarrow \; \int_{\mathbb{R}^d} h(x-y) \, \mu_\infty(\mathrm{d}y) = \nu_\infty(x) \tag{B.1}$$

as $k \to \infty$. Moreover, $|\nu_k(x)| \leq \|h\|_\infty$ and $|\nu_k(x) - \nu_k(x')| \leq \|h\|_1^* |x - x'|$. Thus $(\nu_k)_{k \in \mathbb{N}}$ is equicontinuous and uniformly bounded, and by the Arzelà–Ascoli theorem and uniqueness of the limit in (B.1) we obtain uniform convergence $\nu_k \to \nu_\infty$ on compact sets in $\mathbb{R}^d$ as $k \to \infty$.

Now fix $1 \leq p < \infty$. Jensen's inequality and Tonelli's theorem give

$$\int_{|x|>R} |\nu_k(x)|^p \mathrm{d}x = \int_{|x|>R} |(\mu_k * h)(x)|^p \mathrm{d}x \leq \int \left( \int_{|x|>R} |h(x-y)|^p \mathrm{d}x \right) \mu_k(\mathrm{d}y).$$

If $\mathrm{spt}(h) \subset B_{R_h}$ and $R > 2R_h$, then

$$\int_{|x|>R} |h(x-y)|^p \mathrm{d}x \leq \begin{cases} 0, & |y| \leq R/2, \\ \|h\|_p^p, & |y| > R/2. \end{cases}$$

Thus

$$\int_{|x|>R} |\nu_k(x)|^p \, \mathrm{d}x \; \leq \; \mu_k(\{|y| > R/2\}) \, \|h\|_p^p. \tag{B.2}$$

Since $(\mu_k)_{k\in\mathbb{N}}$ converges narrowly to $\mu_\infty$, the right-hand side can be made small for $R$ large enough uniformly in $k \in \mathbb{N} \cup \{\infty\}$.

Therefore, on $B_R$, uniform convergence implies $\|\nu_k(x) - \nu_\infty(x)\|_{\mathbb{L}^p(B_R, \mathrm{d}x)} \to 0$ as $k \to \infty$. On $B_R^c$, the estimate in (B.2) shows that $\|\nu_k - \nu_\infty\|_{\mathbb{L}^p(B_R^c, \mathrm{d}x)}$ is arbitrarily small for large $R$, uniformly in $k \in \mathbb{N} \cup \{\infty\}$. Combining both parts yields $\|\nu_k - \nu_\infty\|_{\mathbb{L}^p(\mathbb{R}^d)} \to 0$ as $k \to \infty$ for all $1 \leq p < \infty$. □

**Proof of Lemma 2.7** Let $K \subseteq \mathbb{R}^d$ be a compact set and denote the support of $h$ by

$$K_h := \mathrm{cl}_{(\mathbb{R}^d, |\cdot|)}\{x \in \mathbb{R}^d : |h(x)| > 0\}.$$

This set is compact, too, since we chose $h \in \mathbf{C}_c^\infty(\mathbb{R}^d)$. Thus also the set

$$K' := K - K_h = \{x - y : x \in K, y \in K_h\}$$

is compact. Note that for all $\delta > 0$, the support of $h_\delta$ is contained in $K_h$. Now consider a bounded measurable function $f : \mathbb{R}^d \to \mathbb{R}$ and set $f^\delta(x) := \int_{\mathbb{R}^d} f(y) h_\delta(x - y) \, \mathrm{d}y$. Since $x \in K$ implies $x - y \in K'$ if $y \in K_h$, it can be verified directly from the definitions that $f^\delta(x) = (f \mathbb{1}_{K'})^\delta(x)$ for all $x \in K$; see e.g. Makarov and Podkorytov [28, Lem. 7.5.4] for the full argument. Moreover, $|(f\mathbb{1}_{K'})^\delta(x)| \leq \|h_\delta\|_{\mathbb{L}^\infty(\mathrm{d}x)} \|f\mathbb{1}_{K'}\|_{\mathbb{L}^1(\mathrm{d}x)}$ by Hölder's inequality and the translation invariance of the norms on the right-hand side. Since we also have $\|f\mathbb{1}_{K'}\|_{\mathbb{L}^1(\mathrm{d}x)} = \|f\|_{\mathbb{L}^1(K', \mathrm{d}x)}$, we obtain $|f^\delta(x)| \leq \|h_\delta\|_{\mathbb{L}^\infty(\mathrm{d}x)} \|f\|_{\mathbb{L}^1(K', \mathrm{d}x)}$ for all $x \in K$.

Now let $\mu_\infty \in \mathbf{M}_1^+$ be absolutely continuous with $\mathrm{d}\mu/\mathrm{d}x \in \mathbb{L}^1(\mathrm{d}x)$ and $(\mu_k)_{k\in\mathbb{N}}$ a sequence in $\mathbf{M}_1^+$ such that $\mu_k \to \mu_\infty$ narrowly as $k \to \infty$. Applying the preceding estimates to $f(x) = b_t(x, \mu_k^\delta) - b_t(x, \mu_\infty^\delta)$ and noting that by linearity of the convolution we have $(b_t(x, \mu_k^\delta) - b_t(x, \mu_\infty^\delta))^\delta = b_t^\delta(x, \mu_k^\delta) - b_t^\delta(x, \mu_\infty^\delta)$, we get for each $k \in \mathbb{N}$ that

$$\sup_{x\in K} |b_t^\delta(x, \mu_k^\delta) - b_t^\delta(x, \mu_\infty^\delta)| \leq \|h^\delta\|_{\mathbb{L}^\infty(\mathrm{d}x)} \|b_t(\mu_k^\delta) - b_t(\mu_\infty^\delta)\|_{\mathbb{L}^1(K', \mathrm{d}x)}.$$

Since $\mu_k \to \mu_\infty$ narrowly, Lemma B.1 implies that $\mu_k^\delta \to \mu_\infty^\delta$ in $\mathbb{L}^1(\mathrm{d}x)$ as $k \to \infty$. Moreover, since $\mu_k^\delta$ and $\mu_\infty^\delta$ are absolutely continuous and $\mathrm{d}\mu_k^\delta/\mathrm{d}x$ and $\mathrm{d}\mu_\infty^\delta/\mathrm{d}x$ are in $\mathbb{L}^1(\mathrm{d}x)$ for all $k \in \mathbb{N}$, Assumption 1.7 gives us that $\lim_{k\to\infty} \|b_t(\mu_k^\delta) - b_t(\mu_\infty^\delta)\|_{\mathbb{L}^1(K', \mathrm{d}x)} = 0$. This shows that the continuity property $(C_s)$ from Lemma 2.7 is satisfied. Since the continuity property $(C_s)$ implies the property $(C_w)$ from Assumption 1.7, the proof is complete. □

**Proof of Lemma 2.9** The proof is similar to that in [8, Lem. 1.4] and [10, Prop. 2.6]. For each $\phi \in \mathcal{S}$, define on $\mathbf{C}([0, T]; \mathbb{R})$ the probability measure $\nu_k^\phi := \mathbf{P}_{x_0}^{\infty, k} \circ (\Lambda[\phi])^{-1}$. By Kallianpur and Xiong [24, Thm. 2.4.4], it is enough to prove that $(\nu_k^\phi)_{k \in \mathbb{N}}$ is tight for each $\phi \in \mathcal{S}$. It follows from Proposition 2.8, more specifically that $\Lambda[\phi]$ is a real-valued $(\mathbf{P}_{x_0}^{\infty, k}, \mathbb{F})$-semimartingale with characteristics given in terms of $b^{\delta_k}$, $\sigma$, $\bar{\sigma}$. By part 3) of Lemma B.2 and Remark 2.10, these functions and thus the characteristics are bounded uniformly in $k \in \mathbb{N}$. Thus tightness of $(\nu_k^\phi)_{k \in \mathbb{N}}$ follows from standard criteria for real-valued semimartingales; see e.g. Stroock and Varadhan [40, Thm. 1.4.6]. The fact that the process $\Lambda$ is $\mathbf{P}_{x_0}^{\infty, \infty}$-a.s. $\mathbf{M}_1^+$-valued is shown as in [8, Lem. 1.4], using that the bounds there depend only on $\|b^\delta\|_\infty \leq \|b\|_\infty$ and $\|\sigma\|_\infty$, $\|\bar{\sigma}\|_\infty$. Finally, the existence of $p_1 > 0$ such that $\mathbf{P}^{\infty, k}[\mathbf{C}([0, T]; \mathcal{H}_{-p})] = 1$ for all $k \in \mathbb{N} \cup \{\infty\}$ and $p > p_1$ follows directly from the estimate in (2.23) together with [24, Thm. 2.5.2]. $\square$

## B.2 Some auxiliary results for the proofs in Section 5

Let $h \in \mathbf{C}_c^\infty(\mathbb{R}^d)$ be a symmetric function such that $h \geq 0$ and $\int_{\mathbb{R}^d} h(x) \, dx = 1$, and let $h_\delta := \delta^d h(\delta x)$ for $\delta > 0$. Then $(h_\delta)_{\delta > 0}$ is a compactly supported approximate identity; see Rudin [35, Def. 6.31]. Recall from (1.5) that

$$b^\delta : [0, T] \times \mathbb{R}^d \times \mathbf{M}_1^+ \to \mathbb{R}^d$$

is defined for each $\delta > 0$ by

$$b_t^\delta(x, \lambda) := \big( b_t(\lambda) * h_\delta \big)(x) := \int_{\mathbb{R}^d} b_t(y, \lambda) h_\delta(x - y) \, dy \,,$$

where we denote by $b_t(\lambda)$ for $(t, \lambda) \in [0, T] \times \mathbf{M}_1^+$ the function $b_t(\,\cdot\,, \lambda) : \mathbb{R}^d \to \mathbb{R}^d$.

For the next lemma, we recall from before Assumption 1.7 that $\mathcal{E}$ denotes the product-$\sigma$-algebra on $[0, T] \times \mathbb{R}^d \times \mathbf{M}_1^+$ generated by the Lebesgue-measurable sets of $[0, T] \times \mathbb{R}^d$ and the Borel-measurable sets of $\mathcal{P}_{\mathrm{wk}^*}$, the set of probability measures on $\mathcal{B}(\mathbb{R}^d)$ endowed with the narrow topology.

The function $b^\delta$ enjoys the following properties.

**Lemma B.2** | *Let $b : [0, T] \times \mathbb{R}^d \times \mathbf{M}_1^+ \to \mathbb{R}^d$ satisfy Assumption 1.2 and fix $\delta > 0$. Then $b^\delta$ is $\mathcal{E}/\mathcal{B}(\mathbb{R}^d)$-measurable. Moreover:*

*1) For each $\lambda \in \mathbf{M}_1^+$ and $t \in [0, T]$, $x \mapsto b_t^\delta(x, \lambda)$ is continuous.*

*2) If $K \subseteq \mathbb{R}^d$ is compact and $1 \leq r' < \infty$, then for each $\lambda \in \mathbf{M}_1^+$ and $t \in [0, T]$, we have $\lim_{\delta \to 0} \|b_t(\lambda) - b_t^\delta(\lambda)\|_{\mathbb{L}^{r'}(K, \, dx)} = 0$.*

*3) For each $\lambda \in \mathbf{M}_1^+$ and $t \in [0, T]$, $\|b_t^\delta(\lambda)\|_\infty \leq \|b_t(\lambda)\|_{\mathbb{L}^\infty(\mathbb{R}^d, \, dx)}$.*

*4) Let $K_h := \mathrm{cl}_{(\mathbb{R}^d, |\cdot|)}\{x \in \mathbb{R}^d : |h(x)| > 0\}$, and $K \subseteq \mathbb{R}^d$ be compact. Then the set $K' := K - K_h = \{x - x' : x \in K, x' \in K_h\}$ is compact, and for all $\lambda \in \mathbf{M}_1^+$, $t \in [0, T]$ and $r' \in [1, \infty]$, we have $\|b_t^\delta(\lambda)\|_{\mathbb{L}^{r'}(K, \, dx)} \leq \|b_t(\lambda)\|_{\mathbb{L}^{r'}(K', \, dx)}$.*

**Proof** Measurability is established in e.g. Makarov and Podkorytov [28, Lem. 7.5.2]. For part 1), see e.g. Brézis [6, Prop. 4.20].

For parts 2)–4), we set $f^\delta(x) := \int_{\mathbb{R}^d} f(y) h_\delta(x-y) \, \mathrm{d}y$ for bounded and measurable $f : \mathbb{R}^d \to \mathbb{R}$ and observe the two following facts. First, defining $K'$ as in part 4), it can be verified directly from the definitions that $f^\delta(x) = (f \mathbb{1}_{K'})^\delta(x)$ for all $x \in K$; see e.g. [28, Lem. 7.5.4]. Second, the following version of Young's inequality holds (see e.g. Adams and Fournier [1, Thm. 2.24]): For all $r_1, r_2, r_3 \in [1, \infty]$ satisfying $1/r_1 + 1/r_2 + 1/r_3 = 2$ and all $w \in \mathbb{L}^{r_3}(\mathbb{R}^d, \mathrm{d}x)$, we have

$$\left| \int_{\mathbb{R}^d} f^\delta(x) w(x) \, \mathrm{d}x \right| \leq \|f\|_{\mathbb{L}^{r_1}(\mathbb{R}^d, \mathrm{d}x)} \|h\|_{\mathbb{L}^{r_2}(\mathbb{R}^d, \mathrm{d}x)} \|w\|_{\mathbb{L}^{r_3}(\mathbb{R}^d, \mathrm{d}x)}.$$

We can now prove part 2). Indeed, setting $f = b_t(\lambda) \mathbb{1}_{K'}$, which is bounded and measurable by Assumption 1.7, we have by [6, Thm. 4.22] that $\lim_{\delta \to 0} \|f - f^\delta\|_{\mathbb{L}^r(\mathbb{R}^d, \mathrm{d}x)} = 0$ whenever $r' < \infty$. Thus also $\lim_{\delta \to 0} \|f - f^\delta\|_{\mathbb{L}^r(K, \mathrm{d}x)} = 0$. But since $f = b_t(\lambda) \mathbb{1}_{K'}$, we have by the first fact observed above that

$$\left\| b_t(\lambda) \mathbb{1}_{K'} - \left( b_t(\lambda) \mathbb{1}_{K'} \right)^\delta \right\|_{\mathbb{L}^r(K, \mathrm{d}x)} = \|b_t(\lambda) - b_t(\lambda)^\delta\|_{\mathbb{L}^r(K, \mathrm{d}x)}.$$

So part 2) follows.

For part 3), we use Hölder's inequality and the fact that we have $\|h_\delta\|_{\mathbb{L}^1} = 1$ to deduce $|b_t^\delta(x, \lambda)| \leq \|b_t(\lambda)\|_{\mathbb{L}^\infty}$ for all $x \in \mathbb{R}^d$, so $\|b_t(\lambda)\|_\infty \leq \|b_t(\lambda)\|_{\mathbb{L}^\infty}$.

Finally, for part 4), we once more consider $f = b_t(\lambda) \mathbb{1}_{K'}$. Let $w$ be a function that is supported in $K$ with $\|w\|_{\mathbb{L}^{r'}(K, \mathrm{d}x)} = 1$. Young's inequality with $r_1 = r'/(r'-1)$, $r_2 = 1$ and $r_3 = r'$ shows that

$$\left| \int_K f^\delta(x) w(x) \, \mathrm{d}x \right| = \left| \int_{\mathbb{R}^d} f^\delta(x) w(x) \, \mathrm{d}x \right| \leq \|f\|_{\mathbb{L}^r(\mathbb{R}^d, \mathrm{d}x)}.$$

But from the dual characterization of the $\mathbb{L}^r(K, \mathrm{d}x)$-norm, we get

$$\|f^\delta\|_{\mathbb{L}^r(K, \mathrm{d}x)} = \sup \left\{ \left| \int_K f^\delta(x) w(x) \, \mathrm{d}x \right| \, : \, \|w\|_{\mathbb{L}^{r'}(K, \mathrm{d}x)} = 1 \right\}.$$

Using the definition $f = b_t(\lambda) \mathbb{1}_{K'}$, the first fact observed above gives

$$\|f^\delta\|_{\mathbb{L}^r(K, \mathrm{d}x)} = \|(\mathbb{1}_{K'} b_t(\lambda))^\delta\|_{\mathbb{L}^r(K, \mathrm{d}x)} = \|b^\delta\|_{\mathbb{L}^r(K, \mathrm{d}x)},$$

and so

$$\|b^\delta\|_{\mathbb{L}^r(K, \mathrm{d}x)} = \|f^\delta\|_{\mathbb{L}^r(K, \mathrm{d}x)} \leq \|f\|_{\mathbb{L}^r(\mathbb{R}^d, \mathrm{d}x)} = \|b_t(\lambda)\|_{\mathbb{L}^r(K', \mathrm{d}x)},$$

where the last equality uses once more the definition $f = b_t(\lambda) \mathbb{1}_{K'}$. This establishes part 4) and completes the proof. □

August 29, 2025